\documentstyle[amscd,amssymb,verbatim,12pt]{amsart}
\newcommand{\B}{\mbox{\rm B}}
\newcommand{\cen}{\mbox{\rm center}}
\newcommand{\ch}{\mbox{\rm ch}}
\newcommand{\Inv}{\mbox{\rm Inv}}
\newcommand{\const}{\mbox{\rm const.}}
\newcommand{\SpRad}{\mbox{\rm SpRad}}

\newcommand{\Id}{\mbox{\rm Id}}

\newcommand{\Ker}{\mbox{\rm Ker}}

\newcommand{\Image}{\mbox{\rm Im}}
\newcommand{\Aut}{\mbox{\rm Aut}}
\newcommand{\tr}{\mbox{\rm tr}}
\newcommand{\Tr}{\mbox{\rm Tr}}
\newcommand{\TR}{\mbox{\rm TR}}
\newcommand{\Pf}{\mbox{\rm Pf}}

\newcommand{\End}{\mbox{\rm End}}

\newcommand{\rk}{\mbox{\rm rk}}

\newcommand{\HH}{\mbox{\rm H}}

\newcommand{\vol}{\mbox{\rm vol}}
\newcommand{\Hom}{\mbox{\rm Hom}}

\newcommand{\Diff}{\mbox{\rm Diff}}
\newcommand{\Homeo}{\mbox{\rm Homeo}}

\newcommand{\Q}{{\Bbb Q}}
\newcommand{\R}{{\Bbb R}}
\newcommand{\C}{{\Bbb C}}
\newcommand{\N}{{\Bbb N}}
\newcommand{\Z}{{\Bbb Z}}
\theoremstyle{plain}
\newtheorem{definition}{Definition}

\newtheorem{lemma}{Lemma}
\newtheorem{theorem}{Theorem}
\newtheorem{proposition}{Proposition}
\newtheorem{corollary}{Corollary}

\newtheorem{hypothesis}{Hypothesis}
\numberwithin{equation}{section}
\renewcommand{\rm}{\normalshape}
\begin{document}
\title{Diffeomorphisms, Analytic Torsion and Noncommutative Geometry}
\author{John Lott}
\address{Department of Mathematics\\
University of Michigan\\
Ann Arbor, MI  48109\\
USA}
\email{lott@@math.lsa.umich.edu}
\thanks{Research supported by NSF grant DMS-9403652}
\date{July 18, 1996}
\maketitle
\begin{abstract}
We prove an index theorem concerning the pushforward of flat ${\frak B}$-vector
bundles, where ${\frak B}$ is an appropriate algebra.  We construct an
associated analytic torsion form ${\cal T}$. 
If $Z$ is a smooth closed aspherical 
manifold, we show that ${\cal T}$ gives invariants of 
$\pi_*(\Diff(Z))$.
\end{abstract}
\section{Introduction}
Let $Z$ be a smooth connected closed $n$-dimensional $K(\Gamma, 1)$-manifold.
Let $\Diff(Z)$ be the group of diffeomorphisms of $Z$, with its
natural smooth topology \cite{Milnor (1984)}.
What are the rational homotopy groups of $\Diff(Z)$? Farrell and
Hsiang made the following conjecture : \\ \\
{\bf Conjecture } \cite{Farrell-Hsiang (1978)} : \label{conj}
$\pi_1(\Diff(Z)) \otimes_{\Z} \Q = \cen (\Gamma) \otimes_{\Z} \Q$ \\ \\
and if $i > 1$ is sufficiently small compared to $n$,
\begin{equation} \label{pii}
\pi_i(\Diff(Z)) \otimes_{\Z} \Q = 
\begin{cases}
\bigoplus_{j=1}^\infty \HH_{i+1-4j}
(\Gamma; \Q) & \text{ if } n \text{ is odd},\\
0 & \text{ if } n \text{ is even}. \\
\end{cases}
\end{equation}

It follows from the work of Farrell and Jones \cite{Farrell-Jones (1993)} that
the conjecture is true when $n > 10$, $i < \frac{n-7}{3}$ and
$\Gamma$ is a discrete cocompact
subgroup of a Lie group with a finite number of connected components.
(For example, it is true when $Z$ is a torus, something which was already
shown in \cite{Farrell-Hsiang (1978)}.)
The $\pi_1$-result is what one would expect from
homotopy theory. However, (\ref{pii}) is peculiar to the fact that we are
looking at diffeomorphisms; the analogous rational homotopy 
groups of $\Homeo(Z)$ vanish. 
In the cases when the conjecture has been proven, the proofs are 
very impressive but rather indirect, using a great deal of
topological machinery.

From a constructive viewpoint,
suppose that we are given a smooth based map $\alpha : S^{i} \rightarrow
\Diff(Z)$. How could we compute the corresponding rational
homotopy class $[\alpha]_{\Q} \in 
\pi_i(\Diff(Z)) \otimes_{\Z} \Q$? First, let us make an auxiliary fiber 
bundle. Using $\alpha$, we can glue
two copies of $D^{i+1} \times Z$ along their boundaries to obtain a
smooth manifold $M$ which fibers over $S^{i+1}$, with fiber $Z$. 
Any (smooth) topological
invariant of fiber bundles will give an invariant of $\pi_i(\Diff(Z))$. 

Wagoner suggested \cite{Wagoner (1978)} 
that the relevant invariant
is a fiber-bundle extension of the Ray-Singer analytic torsion
\cite{Ray-Singer (1971)}. In \cite{Bismut-Lott (1995)}, J.-M. Bismut and
the author constructed a certain extension of the Ray-Singer analytic torsion
which does give some information about $\pi_*(\Diff(Z))$.
However, that extension is inadequate to capture all of the information in
(\ref{pii}). In this paper, using ideas from noncommutative geometry, we
will construct a ``higher'' analytic torsion which does
potentially detect the right-hand-side of (\ref{pii}).
 
One can think of the 
analytic torsion as arising from the transgression of certain index theorems.
We describe the relevant index theorems.  Suppose that $M 
\stackrel{\pi}{\rightarrow} B$ is a smooth fiber bundle with 
connected closed fibers $Z$.
Let $E$ be a flat complex vector bundle on $M$.   
In \cite{Bismut-Lott (1995)}, certain
characteristic classes $c(E) \in \HH^{odd}(M; \R)$ were defined.
The pushforward of $E$ is defined to be
\begin{equation} \label{push}
\pi_*(E) = \sum_{p=0}^{dim(Z)} (-1)^p \: \HH^p(Z; E \big|_Z),
\end{equation} a formal 
alternating sum of flat vector bundles on $B$ constructed from the
cohomology groups of the fibers. Let
$e(TZ) \in \HH^{dim(Z)}(M; o(TZ))$ be the Euler class of the vertical
tangent bundle $TZ$.
The index theorem of 
\cite[Theorem 0.1]{Bismut-Lott (1995)} 
stated
\begin{equation} \label{BLthm}
c(\pi_* E) = \int_Z e(TZ) \cup c(E) \hspace{.5in} \text{ in } 
\HH^{odd}(B; \R). 
\end{equation}
In its proof, which was analytic
in nature, a certain differential form ${\cal T} \in \Omega^{even}(B)$ 
appeared, called the analytic torsion form.

This index theorem was reproved topologically and extended by Dwyer,
Weiss and Williams
\cite{Dwyer-Weiss-Williams (1995)}. Their setup was a fiber bundle as
above, a ring ${\frak B}$ and a local system ${\cal E}$ of 
finitely-generated projective ${\frak B}$-modules on $M$. The
local system defines a class $[{\cal E}] \in K^{alg}_{\frak B}(M)$ in
a generalized cohomology group of $M$.
One again has
local systems $\{\HH^p(Z; {\cal E} \big|_Z)\}_{p=0}^{dim(Z)}$ 
of finitely-generated
${\frak B}$-modules on $B$. Suppose that they are projective.
Define $\pi_*({\cal E})$ as in (\ref{push}). 
Then \cite[Equation (0-3)]{Dwyer-Weiss-Williams (1995)} stated
\begin{equation} \label{DWWthm}
[\pi_*({\cal E})] = \tr^* [{\cal E}]
\hspace{.5in} \text{ in } K^{alg}_{\frak B}(B),
\end{equation}
where $\tr^*$ is the Becker-Gottlieb-Dold transfer.  When
${\frak B} = \C$, (\ref{BLthm}) is a consequence of
applying the characteristic class $c$ to both sides of (\ref{DWWthm}).

In the present paper, we essentially give an analytic proof of (\ref{DWWthm}).
Provided that ${\frak B}$ is an algebra over $\C$ which satisfies certain
technical conditions,
we define a characteristic class
\begin{equation}
[CS] : K^{alg}_{\frak B}(M) \longrightarrow 
\bigoplus 
\begin{Sb}
p > q \\
p + q \: odd
\end{Sb}
\HH^p(M; \overline{H}_q({\frak B})),
\end{equation}
where $\overline{H}_*({\frak B})$ is the noncommutative de Rham cohomology
of the algebra ${\frak B}$.
Using analytic methods, we prove the following theorem.
\begin{theorem} \label{introprop}
Let the fiber bundle $M \stackrel{\pi}{\rightarrow} B$ be as above.
Let ${\cal E}$ be a local system of finitely-generated projective
${\frak B}$-modules on $M$. Suppose that the
fiberwise differentials $\overline{d^Z}$ have closed image. Then
\begin{equation}
\left[ CS \left( \pi_*{\cal E} )
\right) \right] = \int_Z e(TZ) \cup
\left[CS \left( {\cal E} \right)\right]
\text{ in }
\bigoplus 
\begin{Sb}
p > q \\
p + q \: odd
\end{Sb}
\HH^p(B; \overline{H}_q({\frak B})). 
\end{equation}
\end{theorem}
The condition that $\overline{d^Z}$ have closed image guarantees that
$\HH^* (Z; {\cal E} \big|_{Z} )$ is a local system of  projective
${\frak B}$-modules. If ${\frak B} = \C$ then 
\begin{equation}
\overline{H}_q({\frak B}) =
\begin{cases}
\C & \text{ if } q = 0, \\
0 & \text{ if } q > 0 \\
\end{cases}
\end{equation}
and so we recover (\ref{BLthm}).

The statement of Theorem \ref{introprop} can also be obtained by applying the
characteristic class $[CS]$ to both sides of (\ref{DWWthm}).
As in \cite{Bismut-Lott (1995)}, the interest of the analytic proof is
that it gives a more refined statement at the level of differential forms.
With notation that will be explained later, a certain explicit
differential form ${\cal T} \in 
\overline{\Omega}^{\prime \prime,even}
(B, {\frak B})$ appears naturally in the proof of Theorem \ref{introprop}. 
We call it the analytic torsion form. 

\begin{theorem} \label{introprop2}
With the hypotheses of Theorem \ref{introprop},
\begin{equation} \label{introeqn}
d {\cal T} = \int_Z e \left(TZ, \nabla^{TZ} \right) \wedge
CS \left( \nabla^{\cal E}, h^{\cal E} \right) -
CS \left( \nabla^{\pi_*{\cal E}},
h^{\pi_*{\cal E}} \right)
\text{ in } \overline{\Omega}^{\prime \prime, odd}(B, {\frak B}).
\end{equation}
\end{theorem}
\noindent
Here $h^{\cal E}$ is a ${\frak B}$-valued Hermitian metric on ${\cal E}$ and
$h^{\pi_*{\cal E}}$ is the induced ${\frak B}$-valued Hermitian metric on 
$\pi_*{\cal E}$.

If $\dim(Z)$ is odd and $\HH^*(Z; {\cal E} \big|_{Z} ) = 0$ then the
right-hand-side of (\ref{introeqn}) vanishes automatically, implying that
${\cal T}$ is closed.

\begin{theorem} \label{introprop3}
With the hypotheses of Theorem \ref{introprop},
suppose in addition that $\dim(Z)$ is
odd and $\HH^* (Z; {\cal E} \big|_{Z} ) = 0$.
Then the cohomology class 
\begin{equation}
[{\cal T}] \in \bigoplus
\begin{Sb}
p > q\\
p + q \: even\\
\end{Sb}
\HH^p(B; \overline{H}_q({\frak B}))
\end{equation}
is a (smooth) topological invariant
of the fiber bundle $M \stackrel{\pi}{\rightarrow} B$ and the local
system ${\cal E}$.
\end{theorem}

In particular, $[{\cal T}]$ will give invariants of
$\pi_*(\Diff(Z))$ when $Z$ is a smooth closed aspherical manifold.

We now describe the contents of this paper.
The local system ${\cal E}$ on $M$
can be thought of as a flat
${\frak B}$-vector bundle.  The first order of business is to define
the relevant characteristic classes of ${\cal E}$. Unlike in 
\cite{Bismut-Lott (1995)}, it is not enough to just use the flat connection
on ${\cal E}$. Instead, we will need a connection on ${\cal E}$ which
also differentiates in the ``noncommutative'' directions.  The
correct notion, due to Karoubi, is that of a partially flat connection
(called a ``connexion \`a courbure plate'' in \cite{Karoubi (1987)}).

In Section \ref{Noncommutative Bundle Theory} we briefly review
the geometry of ${\frak B}$-vector bundles. We define certain complexes
of noncommutative differential forms and describe their cohomologies.
We review the notion of a ${\frak B}$-connection on ${\cal E}$ and its Chern
character.  We define the relative
Chern-Simons class of two ${\frak B}$-vector bundles which are topologically 
isomorphic, each having a partially flat connection.
Because our connections are
partially flat and not completely flat, the formalism involved is
fundamentally different than that of \cite{Bismut-Lott (1995)}.

In Section \ref{Groups and Covering Spaces} we look at the case when
$\Gamma$ is a finitely-generated discrete group and 
${\frak B}$ lies between the group algebra $\C \Gamma$ and the group 
$C^*$-algebra $C^*_r \Gamma$. If $M$ is a manifold with a normal 
$\Gamma$-covering $M^\prime \rightarrow M$
then there is a canonical ${\frak B}$-vector bundle 
${\cal E} = {\frak B} \times_\Gamma M^\prime$ on $M$. We describe an
explicit partially flat connection on ${\cal E}$ and compute the pairing
of its Chern character with the group cohomology of $\Gamma$. This
computation is important for applications. 

In Section \ref{Hermitian Metrics and Characteristic Classes} we 
define the notion of a ${\frak B}$-valued Hermitian metric $h^{\cal E}$ on a 
${\frak B}$-vector bundle ${\cal E}$. With our assumptions on ${\frak B}$,
a Hermitian metric on ${\cal E}$ always exists and is unique up to isotopy.
A Hermitian metric gives a topological isomorphism between ${\cal E}$ and
its antidual bundle $\overline{\cal E}^*$. If ${\cal E}$ has a partially
flat connection, we can use this isomorphism to define 
the relative Chern-Simons class $CS \left( {\cal E}, h^{\cal E} \right)
\in \overline{\Omega}^{\prime \prime,odd}
(M, {\frak B})$ of ${\cal E}$ and 
$\overline{\cal E}^*$.  Its cohomology class is independent of
the choice of Hermitian metric, giving the characteristic class
\begin{equation}
\left[ CS({\cal E}) \right] \in 
\bigoplus 
\begin{Sb}
p > q \\
p + q \: odd
\end{Sb}
\HH^p(M; \overline{H}_q({\frak B})).
\end{equation}

In Section \ref{Noncommutative Superconnections} we generalize the preceding
results from connections to superconnections. This generalization will
be crucial for the fiber bundle results.
We define the notion of a partially flat
superconnection on $M$. We construct the relative Chern-Simons class
and analytic torsion form. Using these constructions, we prove a
finite-dimensional analog of (\ref{DWWthm}). This analog is similar
to \cite[Theorem 2.19]{Bismut-Lott (1995)}, but the large-time analysis
requires new techniques.
We then relate the finite-dimensional analytic torsion form to various
versions of the Reidemeister torsion. 

In Section \ref{Fiber Bundles} we extend the methods of Section
\ref{Noncommutative Superconnections} to the setting of a fiber bundle
$Z \rightarrow M \stackrel{\pi}{\rightarrow} B$. 
First, we prove some basic facts
about ${\frak B}$-pseudodifferential operators.
Using heat kernel techniques,
we prove Theorem \ref{introprop}.
We then define the analytic torsion form ${\cal T} \in 
\overline{\Omega}^{\prime \prime,even}
(B, {\frak B})$ and prove Theorems \ref{introprop2} and 
\ref{introprop3}.

Relevant examples of the preceding formalism come from finitely-generated
discrete groups $\Gamma$. Let us introduce
a certain hypothesis on $\Gamma$ :
\begin{hypothesis} \label{introhypo}
There is a Fr\'echet locally $m$-convex algebra ${\frak B}$ containing
$\C \Gamma$ such that\\
1. ${\frak B}$ is dense in $C^*_r \Gamma$ and stable under the 
holomorphic functional calculus in $C^*_r \Gamma$.\\
2. For each $[\tau] \in \HH^q(\Gamma; \C)$, there is a
representative cocycle $\tau \in Z^q(\Gamma; \C)$ such that the ensuing
cyclic cocycle $Z_\tau \in HC^q(\C \Gamma)$ extends to
a continuous cyclic cocycle on ${\frak B}$.
\end{hypothesis}
Hypothesis \ref{introhypo} arises in analytic proofs of the Novikov
Conjecture. It is known to be satisfied by virtually 
nilpotent groups and Gromov-hyperbolic groups
\cite[Section 3.5]{Connes (1994)}.
Using the characteristic class
$[CS]$, in Section 
\ref{Hermitian Metrics and Characteristic Classes} 
we give a simple proof that the algebraic K-theory
assembly map is rationally injective for such groups.

Let $Z$ be a smooth connected closed $n$-dimensional manifold with fundamental
group $\Gamma$.  If the
hypotheses of the preceding theorems are satisfied, we can
define invariants of $\pi_i(\Diff(Z))$, $i > 1$, by constructing
the auxiliary fiber bundle mentioned at the beginning, computing its
analytic torsion form ${\cal T}$ and integrating over $B = S^{i+1}$
to get
\begin{equation}
\int_B [{\cal T}] \in \bigoplus
\begin{Sb}
q < i + 1\\
q \equiv i + 1 \mod 2\\
\end{Sb}
\overline{H}_q({\frak B}).
\end{equation}
By Hypothesis \ref{introhypo}, $\int_B [{\cal T}]$ then
pairs with $\HH^q(\Gamma; \C)$.

To make contact with (\ref{pii}), 
in Section \ref{Diffeomorphism Groups} we assume that 
$Z$ is a $K(\Gamma, 1)$-manifold.
In order to satisfy the hypotheses of the preceding theorems, we
would have to know that the differential form Laplacian on 
the universal cover $\widetilde{Z}$ is
invertible in all degrees, something which is probably never the case.
We present two ways to
get around this problem. First, we consider the case when $\Gamma = \Z^n$.
In this case we can apply ordinary ``commutative'' analysis to study the
problem.  We show that we can define a pairing between
$\int_B [{\cal T}]$ and $\HH^q(\Gamma; \C)$ provided that
$q < \min(i+1, n)$. This pairing vanishes for trivial reasons unless
$n$ is odd and $q \equiv i+1 \mod 4$. Second, we consider
general $\Gamma$ satisfying Hypothesis \ref{introhypo}. Using the
fact that the auxiliary fiber bundle $M \stackrel{\pi}{\rightarrow} S^{i+1}$ is
fiber-homotopically trivial, we construct a relative analytic torsion form
${\cal T}$ such that $\int_B [{\cal T}]$ pairs with $\HH^q(\Gamma; \C)$ 
provided that $q < i+1$. Again, the pairing vanishes for trivial reasons unless
$n$ is odd and $q \equiv i+1 \mod 4$. 

Based on a comparison with (\ref{pii}), we expect that the pairing 
between $\int_B [{\cal T}]$ and $\HH^q(\Gamma; \C)$ will
be nonzero if $n$ is odd and $q \equiv i+1 \mod 4$, at least if $i$ is
sufficiently small with respect to $n$. To show this,
one will probably have to make a direct link between the
analytic constructions of the present paper and the topological machinery.

We note that we do not construct a
``higher'' analytic torsion of a single manifold, in the sense of
Novikov's higher signatures.  In the case of a single manifold,
i.e. if the base $B$ of the fiber bundle is a point,
the analytic torsion that we construct in this paper lies in
${\frak B}/\overline{[{\frak B},{\frak B}]}$, something which pairs with the
zero-dimensional cyclic cohomology of ${\frak B}$. The
higher-dimensional cyclic cohomology of ${\frak B}$ only enters when the
base of the fiber bundle is also higher-dimensional.

So far, our topological applications of the analytic torsion form 
are to the rational homotopy of diffeomorphism groups of
aspherical manifolds.  There are also results in the literature about
the rational homotopy of diffeomorphism groups of simply-connected manifolds
\cite[\S 4]{Burghelea (1989)}. It would be interesting to see if there is
an analog of the analytic torsion form in the simply-connected case.

Finally, let us remark that
in \cite{Lott (1992a)}, we constructed a higher eta-invariant of a manifold
with virtually nilpotent fundamental group.  Using the methods of 
Sections \ref{Noncommutative Superconnections} and \ref{Fiber Bundles} of the
present paper, one can relax this condition to allow for Gromov-hyperbolic
fundamental groups.

I thank Alain Connes and Michael Weiss for helpful discussions. I thank
the Max-Planck-Institut-Bonn for its hospitality while this paper
was written.

\section{Noncommutative Bundle Theory} \label{Noncommutative Bundle Theory}
In this section we review some facts about ${\frak B}$-vector bundles and
their characteristic classes.
The material in this section is taken from \cite{Karoubi (1987)},
along with \cite{Connes-Karoubi (1988)} and
\cite{Lott (1992)}.
\subsection{Noncommutative Differential Forms} 
\label{Noncommutative Differential Forms}

Let ${\frak B}$ be a Fr\'echet locally $m$-convex algebra with unit, 
i.e. the projective limit of a sequence
\begin{equation} \label{projlimit}
\ldots \longrightarrow B_{j+1} \longrightarrow B_j \longrightarrow \ldots
\longrightarrow B_0
\end{equation}
of Banach algebras with unit. (A relevant example is
${\frak B} = C^\infty(S^1)$ and
$B_j = C^j(S^1)$.) We recall some basic facts about such algebras 
\cite{Mallios (1986)}. For $j \ge 0$, let $i_j : {\frak B} \rightarrow
B_j$ be the obvious homomorphism. The Banach norm $| \cdot |_j$ on $B_j$ 
induces a submultiplicative seminorm $\parallel \cdot \parallel_j$ 
on ${\frak B}$ by $\parallel b \parallel_j = | i_j(b)|_j$.
Given $b \in {\frak B}$, its
spectrum $\sigma(b) \subseteq \C$ is given by
\begin{equation} \label{specunion}
\sigma(b) = \bigcup_{j=0}^\infty \: \sigma(i_j(b)).
\end{equation}
As each Banach algebra $B_j$ has a holomorphic functional calculus,
it follows from (\ref{specunion}) that ${\frak B}$ also has a holomorphic
functional calculus.

We assume that $B_0$ is a $C^*$-algebra
$\Lambda$, that $i_0$ is injective with dense image and that 
${\frak B}$ is stable under the
holomorphic functional calculus in $\Lambda$. A consequence is that the
invertible elements $\Inv({\frak B})$ are open in ${\frak B}$, as
$\Inv(\Lambda)$ is open in $\Lambda$ and $\Inv({\frak B}) = 
i_0^{-1}(\Inv(\Lambda))$. Furthermore, $\sigma(b) = \sigma(i_j(b))$ for
all $j \ge 0$.

Let us ignore the topology of ${\frak B}$ for a moment.
The universal graded differential algebra (GDA) of ${\frak B}$ is 
\begin{equation}
\Omega_*({\frak B}) = \bigoplus_{k=0}^\infty \Omega_k({\frak B})
\end{equation}
where as a vector space, $\Omega_k({\frak B}) = {\frak B} \otimes
(\otimes^k ({\frak B}/\C))$. As a GDA, $\Omega_*({\frak B})$ is generated
by ${\frak B} = \Omega_0({\frak B})$ and $d{\frak B} \subset 
\Omega_1({\frak B})$ with the relations
\begin{equation}
d1=0, \: \: \: d^2 = 0, \: \: \: 
d(\omega_k \omega_l) = (d \omega_k) \omega_l + (-1)^k \omega_k
d\omega_l
\end{equation}
for $\omega_k \in \Omega_k({\frak B})$, $\omega_l \in \Omega_l({\frak B})$.
It will be convenient to write an element $\omega_k$ of $\Omega_k({\frak B})$
as a finite sum $\sum b_0 db_1 \ldots db_k$. There is a differential
complex 
\begin{equation} \label{diffcomplex}
\overline{\Omega}_*({\frak B}) = \Omega_*({\frak B})/[
\Omega_*({\frak B}),\Omega_*({\frak B})].
\end{equation}
Let $\overline{Z}_*({\frak B})$, $\overline{B}_*({\frak B})$ and
$\overline{H}_*({\frak B})$ denote its cocycles, coboundaries and cohomology,
respectively. The latter
is given by 
\begin{equation} \label{cycliciso}
\overline{H}_*({\frak B})  =
\begin{cases}
\Ker \left(B : HC_0({\frak B}) \: (= {\frak B}/
[{\frak B}, {\frak B}]) \longrightarrow H_1({\frak B}, {\frak B}) \right)
& \text{ if } * = 0,\\
\Ker \left(B : \overline{HC}_*({\frak B})
\longrightarrow H_{*+1}({\frak B}, {\frak B}) \right) & \text{ if } *>0.
\end{cases}
\end{equation}
Here $\overline{HC}_*({\frak B})$ is the reduced cyclic homology of
${\frak B}$ and $H_{*}({\frak B}, {\frak B})$ is the Hochschild homology.
In particular, there is a pairing between $\overline{H}_*({\frak B})$ and
the (reduced) cyclic cohomology of ${\frak B}$.

Taking the topology on ${\frak B}$ into consideration, 
there is a Fr\'echet completion of $\Omega_*({\frak B})$, which we
again denote by $\Omega_*({\frak B})$. Furthermore, there is a Fr\'echet
space $\overline{\Omega}_*({\frak B})$ defined as in (\ref{diffcomplex}), 
except quotienting by the closure of
the commutator. Hereafter, when we refer to spaces of
differential forms we will always mean these Fr\'echet spaces.
Furthermore, all tensor products of Fr\'echet spaces will implicitly be
projective tensor products. We again denote the (separable)
homology of $\overline{{\Omega}}_*({\frak B})$ by
$\overline{H}_*({\frak B})$. It pairs with the
(reduced) topological cyclic cohomology of ${\frak B}$.

Let ${\frak E}$ be a Fr\'echet left ${\frak B}$-module, meaning a Fr\'echet
space which is a continuous left
${\frak B}$-module. Hereafter, we assume that ${\frak E}$ is a 
finitely-generated projective ${\frak B}$-module. If ${\frak F}$ is a Fr\'echet right ${\frak B}$-module 
then there is a Fr\'echet space
${\frak F} \otimes_{\frak B} {\frak E}$. 
If ${\frak F}$ is
a Fr\'echet ${\frak B}$-bimodule then there is a ``trace map''
\begin{equation}
\Tr : \Hom_{\frak B} ({\frak E}, {\frak F} \otimes_{\frak B}
{\frak E}) \longrightarrow {\frak F}/\overline{[{\frak B}, {\frak F}]}.
\end{equation}
(We quotient by the closure of $[{\frak B}, {\frak F}]$ to ensure that the
result lies in a Fr\'echet space.)
If ${\frak F}$ is a Fr\'echet algebra containing ${\frak B}$ then $\Tr$ 
gives a trace
\begin{equation} \label{tracemap}
\Tr : \Hom_{\frak B} ({\frak E}, {\frak F} \otimes_{\frak B}
{\frak E}) \longrightarrow {\frak F}/\overline{[{\frak F}, {\frak F}]}.
\end{equation}
In the case that ${\frak E}$ is $\Z_2$-graded by an operator
$\Gamma_{\frak E} \in \End_{\frak B}({\frak E})$ satisfying 
$\Gamma_{\frak E}^2 = 1$, we can extend $\Tr$ to a supertrace by
\begin{equation}
\Tr_s(T) = \Tr(\Gamma_{\frak E} T).
\end{equation}

Let $M$ be a smooth connected manifold.  
Put 
\begin{alignat}{2}
\Omega^{p,q}(M, {\frak B}) & = \Omega^p(M; \Omega_q({\frak B})), &
\hspace{.25in}
\overline{\Omega}^{p,q}(M, {\frak B}) & = \Omega^p(M; 
\overline{\Omega}_q({\frak B})) \\
\Omega^{k}(M, {\frak B}) & = 
\bigoplus_{p+q=k} \Omega^{p,q}(M, {\frak B}), &
\hspace{.25in}
\overline{\Omega}^{k}(M, {\frak B}) & = 
\bigoplus_{p+q=k} \overline{\Omega}^{p,q}(M, {\frak B}). \notag
\end{alignat}
We also write $C^\infty(M; {\frak B})$ for $\Omega^0(M, {\frak B})$.
There is a total differential $d$ on
$\overline{\Omega}^{*}(M, {\frak B})$ which decomposes as the sum
of two differentials $d = d^{1,0} + d^{0,1}$. 
Put
\begin{align} \label{formspaces}
\overline{\Omega}^{\prime, 2k}(M, {\frak B}) & = 
Z^k(M; \overline{\Omega}_k({\frak B})) \oplus \left( \bigoplus
\begin{Sb}
p + q = 2k\\
p < q
\end{Sb}
\overline{\Omega}^{p,q}(M, {\frak B}) \right), \\
\overline{\Omega}^{\prime, 2k+1}(M, {\frak B}) & = 
\bigoplus
\begin{Sb}
p + q = 2k+1\\
p < q
\end{Sb}
\overline{\Omega}^{p,q}(M, {\frak B}), \notag \\
\overline{\Omega}^{\prime \prime, *}(M, {\frak B}) & = 
\overline{\Omega}^{*}(M, {\frak B})/
\overline{\Omega}^{\prime, *}(M, {\frak B}). \notag
\end{align}
Then $\overline{\Omega}^{\prime, *}(M, {\frak B})$ and 
$\overline{\Omega}^{\prime \prime, *}(M, {\frak B})$
are also differential complexes.
Let $H_{\frak B}^*(M)$, $H_{\frak B}^{\prime, *}(M)$ and 
$H_{\frak B}^{\prime \prime, *}(M)$ denote
the cohomology groups of $\overline{\Omega}^{*}(M, 
{\frak B})$, $\overline{\Omega}^{\prime, *}(M, {\frak B})$ and 
$\overline{\Omega}^{\prime \prime, *}(M, {\frak B})$, respectively.
Then
\begin{align} \label{isos}
H_{\frak B}^k(M) & \cong
\bigoplus_{p+q=k} \HH^p \left(M; \overline{H}_q({\frak B})\right), \\
H_{\frak B}^{\prime, 2k}(M) & \cong
\HH^k(M; \overline{Z}_k({\frak B})) \oplus \left( \bigoplus
\begin{Sb}
p + q = 2k\\
p < q
\end{Sb}
\HH^p(M; \overline{H}_q({\frak B})) \right), \notag \\
H_{\frak B}^{\prime, 2k+1}(M) & \cong
\bigoplus
\begin{Sb}
p + q = 2k+1\\
p < q
\end{Sb}
\HH^p(M; \overline{H}_q({\frak B})), \notag \\
H_{\frak B}^{\prime \prime, 2k}(M) & \cong
\bigoplus
\begin{Sb}
p + q = 2k\\
p > q
\end{Sb}
\HH^p(M; \overline{H}_q({\frak B})), \notag \\
H_{\frak B}^{\prime \prime, 2k+1}(M) & \cong
\HH^{k+1}\left(M; 
\frac{\overline{\Omega}_k({\frak B})}{\overline{B}_k({\frak B})}\right) 
\oplus \left( \bigoplus
\begin{Sb}
p + q = 2k+1\\
p > q + 1
\end{Sb}
\HH^p(M; \overline{H}_q({\frak B})) \right). \notag
\end{align}
To realize the first isomorphism in (\ref{isos}) explicitly, 
if $\omega \in \overline{\Omega}^{k}(M, 
{\frak B})$ is $d$-closed and $z \in Z_p(M; \C)$ then $\int_z \omega \in 
\overline{Z}_{k-p}({\frak B})$. The other isomorphisms can be realized
similarly.
\subsection{Noncommutative Connections and Chern Character}
\label{Connections and Chern Character}

Let ${\cal E}$ be a smooth
${\frak B}$-vector bundle on $M$ with fibers isomorphic to ${\frak E}$. 
This means that if ${\cal E}$ is defined using charts $\{U_\alpha\}$ then
a transition function is a smooth map $\phi_{\alpha \beta} : U_\alpha
\cap U_\beta \rightarrow \Aut_{\frak B} ({\frak E})$. 
There is a corresponding element $[{\cal E}]$ in the topological K-group
$K^{top}_{\frak B}(M) = [M, K_0({\frak B}) \times BGL({\frak B})]$. 

We will denote the
fiber of ${\cal E}$ over $m \in M$ by ${\cal E}_m$. 
If ${\frak F}$ is a Fr\'echet
${\frak B}$-bimodule, let ${\frak F} \otimes_{\frak B}
{\cal E}$ denote the ${\frak B}$-vector bundle on $M$ with fibers
$({\frak F} \otimes_{\frak B} {\cal E})_m = 
{\frak F} \otimes_{\frak B} {\cal E}_m$ and transition functions
$\Id_{\frak F} \otimes_{\frak B} \phi_{\alpha \beta}$.
 
Let $C^\infty(M; {\cal E})$ denote the left ${\frak B}$-module of smooth
sections of ${\cal E}$ and let $\Omega (M; {\cal E})$ denote the 
left ${\frak B}$-module of smooth
sections of $\Lambda(T^*M) \otimes {\cal E}$. We put 
\begin{equation}
\Omega^{p,q}(M, {\frak B}; {\cal E}) = \Omega^p \left(
M; \Omega_q({\frak B}) \otimes_{\frak B} {\cal E} \right)
\end{equation}
and
\begin{equation}
\Omega^{k}(M, {\frak B}; {\cal E}) = 
\bigoplus_{p+q=k} \Omega^{p,q}(M, {\frak B}; {\cal E}).
\end{equation}
\begin{definition}
A connection on ${\cal E}$ is a $\C$-linear map
\begin{equation}
\nabla^{\cal E} : C^\infty(M; {\cal E}) \rightarrow \Omega^{1}(M, {\frak B}; {\cal E}) 
\end{equation}
such that for all $f \in C^\infty(M; {\frak B})$ and $s \in
C^\infty(M; {\cal E})$, 
\begin{equation}
\nabla^{\cal E}(fs) = f \: \nabla^{\cal E} s + df 
\otimes_{C^\infty(M; {\frak B})} s. 
\end{equation}
\end{definition}
We can decompose $\nabla^{\cal E}$ as 
\begin{equation} \label{conndecomp}
\nabla^{\cal E} = 
\nabla^{{\cal E},1,0} \oplus 
\nabla^{{\cal E},0,1},
\end{equation}
 where
\begin{equation}
\nabla^{{\cal E},1,0} : 
C^\infty(M; {\cal E}) \rightarrow \Omega^{1}(M; {\cal E}) 
\end{equation}
is a connection on ${\cal E}$ in the usual sense which happens to be
${\frak B}$-linear, and 
\begin{equation}
\nabla^{{\cal E},0,1} : C^\infty(M; {\cal E}) \rightarrow
C^\infty(M; \Omega_{1}({\frak B}) \otimes_{\frak B} 
{\cal E} )
\end{equation}
is a $C^\infty(M)$-linear map which
comes from a $\C$-linear bundle homomorphism
\begin{equation} \label{partialdef}
\partial^{\cal E} : 
{\cal E} \rightarrow \Omega_{1}({\frak B}) \otimes_{\frak B} 
{\cal E} 
\end{equation}
satisfying
\begin{equation}
\partial^{\cal E}(b s_m) = b \: \partial^{\cal E} s_m + 
db \otimes_{\frak B} s_m
\end{equation}
for all $m \in M$, $s_m \in {\cal E}_m$ and $b \in {\frak B}$.
One can consider $\nabla^{{\cal E},0,1}$ to be the part of 
$\nabla^{{\cal E}}$
which involves differentiation in the ``noncommutative'' direction.

Extend $\nabla^{\cal E}$ to a $\C$-linear map
\begin{equation} \label{connection}
\nabla^{\cal E} : \Omega^{*}(M, {\frak B}; {\cal E}) \rightarrow 
\Omega^{*+1}(M, {\frak B}; {\cal E}) 
\end{equation}
by requiring that for all $\omega \in \Omega^k(M, {\frak B})$ and $s \in
\Omega^{l}(M, {\frak B}; {\cal E})$, 
\begin{equation}
\nabla^{\cal E}(\omega s) = (-1)^k \: \omega \wedge 
\nabla^{\cal E} s + d\omega 
\otimes_{C^\infty(M; {\frak B})} s. 
\end{equation}
Similarly, extend $\nabla^{{\cal E},1,0}$ to
\begin{equation}
\nabla^{{\cal E},1,0} :
\Omega^{*}(M; {\cal E}) \rightarrow \Omega^{*+1}(M; {\cal E}).
\end{equation}

Now $\left(\nabla^{\cal E}\right)^2$ is multiplication by an element of 
$$\bigoplus_{p+q=2} \Omega^p(M; \Hom_{\frak B}({\cal E}, \Omega_q({\frak B})
\otimes_{\frak B} {\cal E})),$$ 
which we also denote by $\left( \nabla^{\cal E}\right)^2$.
\begin{definition}
The Chern character of $\nabla^{\cal E}$ is 
\begin{equation}
\ch(\nabla^{\cal E}) = \Tr \left( e^{- 
\left( \nabla^{\cal E}\right)^2} \right) \in 
\overline{\Omega}^{even}(M, {\frak B}).
\end{equation}
\end{definition}

As usual, $\ch(\nabla^{\cal E})$ is $d$-closed and its cohomology class 
$\left[ \ch(\nabla^{\cal E}) \right] \in {H}^{even}_{\frak B}(M)$ only
depends on $[{\cal E}] \in K^{top}_{\frak B}(M)$. If
${\cal E}$ and ${\cal E}^\prime$ are ${\frak B}$ and 
${\frak B}^\prime$-vector
bundles on $M$, respectively, then
\begin{equation} \label{chmull}
\ch (\nabla^{{\cal E} \otimes_\C {\cal E}^\prime} ) = 
\ch(\nabla^{\cal E}) \cdot \ch(\nabla^{{\cal E}^\prime}) \in
\overline{\Omega}^{even}(M, {\frak B} \otimes_\C {\frak B}^\prime).
\end{equation} 
\begin{definition} \label{partiallyflat}
A connection $\nabla^{\cal E}$ is partially flat if its component
$\nabla^{{\cal E}, 1, 0}$ is flat, meaning
$\left(\nabla^{{\cal E}, 1, 0}\right)^2 = 0$.
\end{definition} 
\begin{definition} A flat structure on ${\cal E}$ is given by a connection
\begin{equation}
\nabla^{{\cal E},flat} : C^\infty(M; {\cal E}) \rightarrow
\Omega^1(M; {\cal E})
\end{equation}
which is ${\frak B}$-linear and whose extension to $\Omega^*(M; {\cal E})$
satisfies $\left( \nabla^{{\cal E},flat} \right)^2 = 0$.
\end{definition}

Clearly a partially flat connection on ${\cal E}$ determines a flat
structure on ${\cal E}$ through its $(1,0)$-part. 
Conversely, given a flat structure on ${\cal E}$,
there is a partially flat connection on ${\cal E}$ which is compatible
with the flat structure, although generally not a unique one.

The flat structure $\left({\cal E}, \nabla^{{\cal E},flat} \right)$ is
classified by a map $M \rightarrow B\Aut_{\frak B}({\frak E})_\delta$, where
$\delta$ denotes the discrete topology. Then there is a composite map
\begin{equation}
M \rightarrow B\Aut_{\frak B}({\frak E})_\delta \rightarrow 
BGL({\frak B})_\delta \rightarrow BGL({\frak B})_\delta^+,
\end{equation}
where $+$ denotes Quillen's plus construction.  Thus the pair $({\cal E},
\nabla^{{\cal E},flat})$ gives an 
element $[{\cal E},\nabla^{{\cal E},flat}] \in K^{alg}_{\frak B}(M) = 
[M, K_0({\frak B}) \times BGL({\frak B})_\delta^+]$,
the $K_0({\frak B})$ factor simply representing the K-theory class of 
the fiber ${\frak E}$.

If $\nabla^{\cal E}$ is partially flat then
\begin{equation}
\left( \nabla^{\cal E} \right)^2 \in
\Omega^1(M; \Hom_{\frak B}({\cal E}, \Omega_1({\frak B})
\otimes_{\frak B} {\cal E})) \oplus
\Omega^0(M; \Hom_{\frak B}({\cal E}, \Omega_2({\frak B})
\otimes_{\frak B} {\cal E})). 
\end{equation}
Thus
$\ch(\nabla^{\cal E}) \in \bigoplus
\begin{Sb}
p \le q \\
p + q \: even
\end{Sb}
\overline{\Omega}^{p,q}
(M, {\frak B})$. As $\ch(\nabla^{\cal E})$ is $d$-closed, its
$(p,p)$-component $\ch^{p,p}$ must satisfy $d^{1,0} \ch^{p,p} = 0$. Hence 
$\ch(\nabla^{\cal E}) \in \overline{\Omega}^{\prime, even}(M, {\frak B})$
and
$[\ch(\nabla^{\cal E})] \in H^{\prime, even}_{\frak B}(M)$.
There is a commutative diagram
\begin{equation}
\begin{array}{ccc}
K^{alg}_{\frak B}(M) & \longrightarrow & K^{top}_{\frak B}(M) \\
\ch \downarrow & & \ch \downarrow \\
H^{\prime, even}_{\frak B}(M) & \longrightarrow &H^{even}_{\frak B}(M). \\
\end{array}
\end{equation}
\noindent
{\bf Example 1 :} If ${\frak B} = \C$ then $\Omega_0({\frak B}) =
\overline{Z}_0({\frak B}) =
\overline{H}_0({\frak B}) = \C$ and $\Omega_*({\frak B}) =
\overline{Z}_*({\frak B}) =
\overline{H}_*({\frak B}) = 0$ for $* > 0$. Then if ${\cal E}$ has a
flat structure, by (\ref{isos}) we have that
$[\ch(\nabla^{{\cal E}})]$ 
lies in $\HH^0(M; \C)$ 
and simply represents $\rk({\cal E})$. On the other hand,
$K^{alg}_{\C}(M)$ can be very rich. Thus the Chern character does not see
the interesting part of $K^{alg}_{\C}$.  
We now give another construction which will be used in
Section \ref{Hermitian Metrics and Characteristic Classes} 
to see more of $K^{alg}_{\frak B}$.

\subsection{Chern-Simons Classes of Partially Flat Connections}
\label{Chern-Simons Classes}

Let ${\cal E}_1$ and ${\cal E}_2$ be smooth ${\frak B}$-vector bundles on $M$
with flat structures.  Suppose that there is a smooth isomorphism
$\alpha : {\cal E}_1 \rightarrow {\cal E}_2$ of 
${\cal E}_1$ and ${\cal E}_2$ as topological ${\frak B}$-vector bundles.
The triple $({\cal E}_1, {\cal E}_2, \alpha)$ defines an element of
Karoubi's relative K-group $K^{rel}_{\frak B}(M)$, which fits into an
exact sequence
\begin{equation}
K^{alg, -1}_{\frak B}(M) \longrightarrow
K^{top, -1}_{\frak B}(M) \longrightarrow
K^{rel}_{\frak B}(M) \longrightarrow
K^{alg}_{\frak B}(M) \longrightarrow
K^{top}_{\frak B}(M).
\end{equation}

Choose partially flat connections $\nabla^{{\cal E}_1}$, $\nabla^{{\cal E}_2}$
which are compatible with the flat structures.  For $u \in [0,1]$, put
$\nabla^{\cal E}(u) = u \nabla^{{\cal E}_1} + (1-u) \,
\alpha^* \nabla^{{\cal E}_2}$. Note that for $u \in (0,1)$, 
$\nabla^{\cal E}(u)$ may not be partially flat on ${\cal E}_1$.
\begin{definition}
The relative Chern-Simons class $CS\left(\nabla^{{\cal E}_1},
\nabla^{{\cal E}_2} \right) \in \overline{\Omega}^{\prime \prime,odd}
(M, {\frak B})$ is
\begin{equation} \label{CSdef}
CS\left(\nabla^{{\cal E}_1},
\nabla^{{\cal E}_2} \right) = 
- \int_0^1 \Tr \left( \left(
\partial_u \nabla^{\cal E}(u) \right)
e^{- \left( \nabla^{\cal E}(u)\right)^2} \right) du.
\end{equation}
\end{definition}

By construction,
\begin{equation} \label{dCS}
d CS\left(\nabla^{{\cal E}_1},
\nabla^{{\cal E}_2} \right) = \ch(\nabla^{{\cal E}_1}) -  
\ch(\nabla^{{\cal E}_2})
\end{equation}
vanishes in $\overline{\Omega}^{\prime \prime,even}(M, {\frak B})$.
Thus there is a class
$\left[ CS\left(\nabla^{{\cal E}_1, flat},
\nabla^{{\cal E}_2, flat} \right) \right] \in 
H^{\prime \prime,odd}_{\frak B}(M)$
which turns out to only depend on $[{\cal E}_1, {\cal E}_2, \alpha]
\in K^{rel}_{\frak B}(M)$. In particular,
$\left[ CS\left(\nabla^{{\cal E}_1, flat},
\nabla^{{\cal E}_2, flat} \right) \right]$ is independent of the choice
of the partially flat connections $\nabla^{{\cal E}_1}$, $\nabla^{{\cal E}_2}$
and only depends on $\alpha$ through its isotopy class; this will also
follow from Proposition \ref{CSvarprop}.

From (\ref{isos}), 
\begin{equation}
\left[ CS\left(\nabla^{{\cal E}_1, flat},
\nabla^{{\cal E}_2, flat} \right) \right] \in
\left( \bigoplus_p \HH^{p+1}\left(M; 
\frac{\overline{\Omega}_p({\frak B})}{\overline{B}_p({\frak B})}\right) \right)
\oplus \left( \bigoplus
\begin{Sb}
p > q + 1\\
p + q \: odd
\end{Sb}
\HH^p(M; \overline{H}_q({\frak B})) \right).
\end{equation}
The next proposition is implicitly contained 
in \cite[p. 444-448]{Connes-Karoubi (1988)}. We give a simpler proof.
\begin{proposition} \label{simpler}
$\left[ CS\left(\nabla^{{\cal E}_1, flat},
\nabla^{{\cal E}_2, flat} \right) \right]$ actually lies in
$\bigoplus 
\begin{Sb}
p > q \\
p + q \: odd
\end{Sb}
\HH^p(M; \overline{H}_q({\frak B}))$.
\end{proposition}
\begin{pf}
Let $CS^{p+1,p} \in \overline{\Omega}^{p+1,p}(M, {\frak B})$ denote
the $(p+1,p)$-component of the explicit differential form in (\ref{CSdef}).
From (\ref{dCS}),  $d^{1,0} \: CS^{p+1,p} = 0$ and so
$CS^{p+1, p}$ defines an element of 
$\HH^{p+1} \left (M;
\frac{\overline{\Omega}_p({\frak B})}{\overline{B}_p({\frak B})} \right)$.
We show that this element lies in 
$\HH^{p+1}(M, \overline{H}_p({\frak B}))$.
For $i \in \{1,2\}$, define
\begin{equation} 
\partial^{{\cal E}_i} : 
{\cal E}_i \rightarrow \Omega_1({\frak B}) \otimes_{\frak B}
{\cal E}_i
\end{equation}
as in (\ref{partialdef}).
Put
\begin{equation}
g(x) = \frac{e^{-{x^2}}-1}{x}.
\end{equation} 
If $z \in Z_{p+1}(M; \C)$ then in $\overline{\Omega}_{p+1}({\frak B})$,
\begin{align}
d^{\frak B} 
\int_z  CS^{p+1,p} & = \int_z d^{0,1} \: CS^{p+1,p} = 
\int_z \left( d^{1,0} \: CS^{p,p+1} + d^{0,1} \: CS^{p+1,p} \right) \\ 
& = \int_z \left[  \Tr \left( e^{- 
\left( \nabla^{{\cal E}_1,1,0} \partial^{{\cal E}_1} \right)^2} \right)
- \Tr \left( e^{- 
\left( \nabla^{{\cal E}_2,1,0} \partial^{{\cal E}_2}
\right)^2} \right) \right] \notag \\ & =
\int_z \left[  \Tr \left( \left( \nabla^{{\cal E}_1,1,0}  
\partial^{{\cal E}_1} \right) 
g(\nabla^{{\cal E}_1,1,0} \partial^{{\cal E}_1}) \right) \right. \notag \\ &
\left. \hspace{1in}
- \Tr \left( \left( \nabla^{{\cal E}_2,1,0}  \partial^{{\cal E}_2} \right)
  g(\nabla^{{\cal E}_2,1,0} \partial^{{\cal E}_2}) \right) 
 \right] \notag \\ & =
\int_z d^{1,0} \left[  \Tr \left( \partial^{{\cal E}_1} 
g(\nabla^{{\cal E}_1,1,0} \partial^{{\cal E}_1})
\right) 
-  \Tr \left( \partial^{{\cal E}_2} g(\nabla^{{\cal E}_2,1,0} 
\partial^{{\cal E}_2}) \right) \right]
= 0. \notag
\end{align}
The proposition follows.
\end{pf}

If ${\cal E}^\prime$ is a ${\frak B}^\prime$-vector 
bundle with a flat structure and a partially flat connection
$\nabla^{{\cal E}^\prime}$ then
\begin{equation} \label{csmul}
CS \left(\nabla^{{\cal E}_1 \otimes_\C {\cal E}^\prime},
\nabla^{{\cal E}_2 \otimes_\C {\cal E}^\prime} \right) =
CS \left(\nabla^{{\cal E}_1},
\nabla^{{\cal E}_2} \right) \cdot
\ch \left( \nabla^{{\cal E}^\prime} \right)
\in
\overline{\Omega}^{\prime \prime,odd}
(M, {\frak B} \otimes_\C {\frak B}^\prime).
\end{equation} 

Finally, for future reference we define a trace on an algebra of
integral operators on ${\cal E}$. Suppose that $M$ is compact and Riemannian.
Let ${\frak F}$ be a Fr\'echet algebra
containing ${\frak B}$. Let
$\Hom^{\infty}_{\frak B}({\cal E}, {\frak F} \otimes_{\frak B} {\cal E})$
be the algebra of integral operators 
$$
T : C^\infty(M; {\cal E}) \rightarrow C^\infty(M;
{\frak F} \otimes_{\frak B} {\cal E})$$
with smooth kernels $T(m_1, m_2) \in 
\Hom_{\frak B}({\cal E}_{m_2}, {\frak F} \otimes_{\frak B} {\cal E}_{m_1})$.
That is, for $s \in C^\infty(M; {\cal E})$,
\begin{equation}
(Ts)(m_1) = \int_M T(m_1, m_2) s(m_2) \: d\vol(m_2) \in 
{\frak F} \otimes_{\frak B} {\cal E}_{m_1}.
\end{equation}
Put 
\begin{equation} \label{trace}
\TR(T) = \int_M \Tr(T(m,m)) \: d\vol(m) \in {\frak F}/
\overline{[{\frak F},{\frak F}]}.
\end{equation}
Then $\TR$ is a trace on 
$\Hom^{\infty}_{\frak B}({\cal E}, {\frak F} \otimes_{\frak B} {\cal E})$.
If ${\cal E}$ is $\Z_2$-graded then there is a supertrace $\TR_s$ on 
$\Hom^{\infty}_{\frak B}({\cal E}, {\frak F} \otimes_{\frak B} {\cal E})$.

\section{Groups and Covering Spaces} \label{Groups and Covering Spaces}
In this section we review the calculation of the cyclic cohomology of
a group algebra. We then describe the relationship between analysis on a 
normal covering space $M^\prime \rightarrow M$ and on a certain
${\frak B}$-vector
bundle ${\cal E}$ over $M$. We put an explicit partially flat connection
on ${\cal E}$ and compute the pairing of its Chern character with the
cohomology of the covering group.
\subsection{Cyclic Cohomology of Group Algebras}
\label{Cyclic Cohomology of Group Algebras}
Let $\Gamma$ be a discrete group.  
Let $\C \Gamma$ be the group algebra of $\Gamma$.  Let $\langle \Gamma
\rangle$ denote the conjugacy classes of $\Gamma$, and $\langle \Gamma
\rangle^\prime \left( {\mbox{\rm resp. }}\langle \Gamma \rangle^{\prime 
\prime} \right)$
those represented by elements of finite (resp. infinite) order.  For $x \in
\Gamma$, let $Z_x$ denote its centralizer in $\Gamma$ and put
$N_x = Z_x/\{x\}$, the quotient of $Z_x$ by the cyclic group generated by
$x$.  If $x$ and $x^\prime$ are conjugate then $N_x$ and $N_{x^\prime}$ are
isomorphic groups, and we will write $N_{\langle x \rangle}$ for their
isomorphism class.  Let $\C[z]$ be a polynomial ring in a variable $z$ of
degree $2$.  Then the cyclic cohomology of $\C \Gamma$ is given
\cite{Burghelea (1985)} by
\begin{equation}
HC^*(\C \Gamma) = \left( \bigoplus_{\langle x \rangle \in \langle \Gamma
\rangle^\prime} \HH^*(N_{\langle x \rangle}; \C) \otimes \C[z] \right)
\oplus  \bigoplus_{\langle x \rangle \in \langle \Gamma
\rangle^{\prime \prime}} \HH^*(N_{\langle x \rangle}; \C).
\end{equation}

We will need explicit cocycles for $HC^*(\C \Gamma)$.  Fix a representative
$x \in \langle x \rangle$.  Put
\begin{align} \label{cocycle}
C^k_x = \{ \tau : & \Gamma^{k+1} \rightarrow \C : \tau \: {\mbox{\rm is skew
and for all }} (\gamma_0, \ldots, \gamma_k) \in \Gamma^{k+1} \: 
{\mbox{\rm and }} z \in Z_x, \\
& \tau (\gamma_0 z, \gamma_1 z, \ldots, \gamma_k z) = \tau (\gamma_0, \gamma_1,
\ldots, \gamma_k) \:
{\mbox{\rm and }} \notag \\
& \tau (\gamma_0 x , \gamma_1, \ldots, \gamma_k) = 
\tau (\gamma_0, \gamma_1, \ldots, \gamma_k)\}. \notag
\end{align}
Let $\delta$ be the usual coboundary operator :
\begin{equation}
(\delta \tau)(\gamma_0, \ldots, \gamma_{k+1}) = \sum_{j=0}^{k+1}
(-1)^j \: \tau(\gamma_0, \ldots, \widehat{\gamma_j}, \ldots, \gamma_{k+1}).
\end{equation}
Denote the resulting cohomology groups by $H^k_x$. Then $H^k_x$ is isomorphic
to $\HH^k(N_{\langle x \rangle}; \C)$ and for each cocycle $\tau \in
Z^k_x$, there is a cyclic cocycle $Z_\tau \in ZC^k(\C \Gamma)$ given by
\begin{equation} \label{cycliccocycle}
Z_\tau (\gamma_0, \gamma_1, \ldots, \gamma_k) =
\begin{cases}
0 \text{ if } \gamma_k \ldots \gamma_0 \notin \langle
x \rangle \\
\tau(\gamma_0 g, \gamma_1 \gamma_0 g, \ldots, \gamma_k \ldots
\gamma_0 g) \text{ if } \gamma_k \ldots \gamma_0 = g x g^{-1}.
\end{cases}
\end{equation}
For $k > 0$, these are in fact reduced cocycles. In particular, from
(\ref{cycliciso}), they pair with $\overline{H}_k(\C \Gamma)$.

\subsection{Noncommutative Geometry of Covering Spaces} \label{Covering Spaces}
The material in this subsection is essentially taken from 
\cite{Lott (1992)}, with a change from right modules to left modules.

Let $\Gamma$ be a finitely generated discrete group.
Let $\parallel
\circ \parallel$ be a right-invariant word-length metric on $\Gamma$.
Put
\begin{equation}
{\frak B}^{\omega} = \{ b : \Gamma \rightarrow \C : {\mbox{\rm for all}}
\: q \in \Z, \sup_{g \in \Gamma} 
\left( e^{q \parallel g \parallel} |b(g)| \right) < \infty \}.
\end{equation}
Then ${\frak B}^\omega$ is independent of the choice of $\parallel
\circ \parallel$ and is a Fr\'echet locally $m$-convex algebra with unit.
Note that ${\frak B}^\omega$ is generally not stable under the 
holomorphic functional calculus in the reduced group $C^*$-algebra
$C^*_r \Gamma$. For this reason, we will eventually replace it with
a larger algebra.  But let us continue with
${\frak B}^\omega$ for the moment.

Let $M$ be a smooth connected compact Riemannian manifold. 
Let $\rho : \pi_1(M) \rightarrow \Gamma$ be a surjective
homomorphism. There is an
induced connected
normal $\Gamma$-covering $M^\prime$ of $M$, on which $g \in \Gamma$
acts on the left by $L_g \in \Diff(M^\prime)$.
Let $\pi : M^\prime \rightarrow M$ be the projection map. Put
\begin{equation}
{\cal D}^\omega = {\frak B}^\omega \times_\Gamma M^\prime.
\end{equation}
Then ${\cal D}^\omega$ is a ${\frak B}^\omega$-vector bundle on $M$ with
a flat structure. 

Let $E$ be a complex vector bundle on $M$ with connection $\nabla^E$ and let
$E^\prime$ be the pulled-back vector bundle $\pi^* E$ on $M^\prime$ with
connection $\nabla^{E^\prime} = \pi^* \nabla^E$. Define
\begin{equation}
{\cal E}^\omega = {\cal D}^\omega \otimes E,
\end{equation}
a ${\frak B}^\omega$-vector bundle on $M$.
Fix a basepoint $x_0 \in M^\prime$. There is an isomorphism
between $C^\infty\left( M; {\cal E}^\omega \right)$ and
\begin{align} \label{srep}
\{ s \in C^\infty(M^\prime; E^\prime) : & \: {\mbox{\rm for all }}
q \in \Z \: {\mbox{\rm and multi-indices }} \alpha, \\ 
& \sup_{x \in M^\prime} 
\left( e^{q d(x_0, x)} |\nabla^\alpha s(x)| \right) < \infty \}. \notag
\end{align}
The action of ${\frak B}^\omega$ on 
$C^\infty\left( M; {\cal E}^\omega \right)$ is given explicitly by
saying that
$b = \sum_{g \in \Gamma} b_g \: g \in {\frak B}^\omega$ sends
$s \in C^\infty(M^\prime; E^\prime)$ to
\begin{equation}
b\cdot s = \sum_{g \in \Gamma} b_g \: L_{g^{-1}}^* s.
\end{equation}

We now construct an explicit partially flat connection 
$\nabla^{{\cal D}^\omega}$ on ${\cal D}^\omega$.
The $(1,0)$-part $\nabla^{{\cal D}^\omega,1,0}$ is determined by the
flat structure on ${\cal D}^\omega$. It remains to construct
\begin{equation}
\nabla^{{\cal D}^\omega,0,1} : C^\infty(M; {\cal D}^\omega)
\rightarrow C^\infty(M; \Omega_1({\frak B}^\omega) \otimes_{{\frak B}^\omega}
{\cal D}^\omega).
\end{equation}
Let $h \in C^\infty_0(M^\prime)$ be a real-valued function satisfying
\begin{equation} \label{heqn}
\sum_{g \in \Gamma} L_g^* h = 1.
\end{equation}
Given $s \in C^\infty(M; {\cal D}^\omega)$, considering it to be an element
of $C^\infty(M^\prime)$ by (\ref{srep}), define its covariant derivative
to be
\begin{equation}
\nabla_g s = h \cdot L_g^* s \in C^\infty(M^\prime).
\end{equation}
\begin{proposition} \label{hprop} \cite[Prop. 9]{Lott (1992)}
\begin{equation}
\nabla^{{\cal D}^\omega,0,1} s = \sum_{g \in \Gamma}
dg \otimes_{C^\infty(M; {\frak B}^\omega)} \nabla_g s
\end{equation}
defines the $(0,1)$-part of a partially flat connection on
${\cal D}^\omega$.
\end{proposition}

We will use the inclusion
\begin{equation}
\Omega^* \left(M,{\frak B}^\omega; {\cal D}^\omega \right) \rightarrow
\Omega_*({\frak B}^\omega) \otimes_{{\frak B}^\omega} \Omega^*(M^\prime).
\end{equation}
The curvature of $\nabla^{{\cal D}^\omega}$, acting on $s \in
C^\infty(M^\prime)$, is computed to be
\begin{equation} \label{curv}
\left( \nabla^{{\cal D}^\omega} \right)^2 s = -
\sum_{g \in \Gamma} dg \left(  d^{M^\prime}h \right) 
L_g^* s + \sum_{g, g^\prime \in 
\Gamma} dg \: dg^\prime \: h \left( L_{g}^* h 
\right) L_{g^\prime g}^* s.
\end{equation}

Let $\tau \in Z^k_x$ be as in (\ref{cocycle}) and let $Z_\tau
\in ZC^k(\C \Gamma)$ be the corresponding cyclic cocycle.
Suppose that there are constants $C, D > 0$ such that for all
$(\gamma_0, \ldots, \gamma_k) \in \Gamma^{k+1}$,
\begin{equation}
|Z_\tau (\gamma_0, \ldots, \gamma_k)| \le C 
e^{D (\parallel \gamma_0 \parallel + \ldots \parallel \gamma_k \parallel)}.
\end{equation}
Then $Z_\tau$ extends to an element of $ZC^k({\frak B}^\omega)$.

The cover $M^\prime$ of $M$ is
classified by a map $\nu : M \rightarrow B\Gamma$, defined up to homotopy.
If $x = e$, we can think of $[\tau]$ as an element of 
$\HH^k(\Gamma; \C) \cong \HH^k(B\Gamma; \C)$. Recall that
$\left[ \ch \left( \nabla^{{\cal D}^\omega} \right) \right] \in
H^{even}_{{\frak B}^\omega}(M)$.

\begin{proposition} \label{injective}
The pairing $\left\langle Z_\tau,
\ch \left( \nabla^{{\cal D}^\omega} \right)
 \right\rangle \in \HH^*(M; \C)$ is given by
\begin{equation}
\left\langle Z_\tau, \ch \left( \nabla^{{\cal D}^\omega} \right) 
\right\rangle = 
\begin{cases}
0 & \text{ if } x \ne e \\
c_k \: \nu^* [\tau] & \text{ if } x=e, \notag
\end{cases}
\end{equation}
where $c_k$ is a nonzero constant which only depends on $k$.
\end{proposition}
\begin{pf}
Let $c_k$ denote a generic nonzero $k$-dependent constant.
We use equation (\ref{curv}) for the curvature of $\nabla^{{\cal D}^\omega}$.
Consider first the term in $\ch \left( \nabla^{{\cal D}^\omega} \right)$
coming from $\left(- \sum_{g \in \Gamma} dg \left(  d^{M^\prime}h \right) 
L_g^* \right)^k$. For $s \in C^\infty(M^\prime)$, we have
\begin{align}
& \left( \sum_{g \in \Gamma} dg \left(  d^{M^\prime}h \right) L_g^* 
\right)^k s = \\
& \sum_{g_1 \ldots g_k} \left[dg_1 \left(  d^{M^\prime}h \right) L_{g_1}^*
\right] \ldots \left[
dg_k \left( d^{M^\prime}h \right) L_{g_k}^* \right] s = \notag \\
& c_k 
\sum_{g_1 \ldots g_k} dg_1 \ldots dg_k  \left(  d^{M^\prime}h \right) 
 \left(  L_{g_1}^* d^{M^\prime} h \right) \ldots
 \left(  L_{g_{k-1} \ldots g_1}^* d^{M^\prime} h \right) 
L_{g_k \ldots g_1}^* s = \notag \\
& c_k 
\sum_{g_1 \ldots g_k} dg_1 \ldots dg_k  \left(  d^{M^\prime}h \right) 
 \left(  L_{g_1}^* d^{M^\prime} h \right) \ldots
 \left(  L_{g_{k-1} \ldots g_1}^* d^{M^\prime} h \right) 
(g_k \ldots g_1)^{-1} \cdot s = \notag \\
& c_k 
\sum_{g_1 \ldots g_k} dg_1 \ldots dg_k  (g_k \ldots g_1)^{-1}
L_{(g_k \ldots g_1)^{-1}}^* \left[
\left( d^{M^\prime}h \right)
\left( L_{g_1}^* d^{M^\prime} h \right) \ldots
 L_{g_{k-1} \ldots g_1}^* d^{M^\prime} h \right] s = \notag \\
& c_k 
\sum_{g_1 \ldots g_k} dg_1 \ldots dg_k  (g_k \ldots g_1)^{-1}
\left( L_{(g_k \ldots g_1)^{-1}}^* d^{M^\prime} h \right) \ldots
\left( L_{g_{k}^{-1}}^* d^{M^\prime} h \right)s.
\end{align}
The contribution of this term to 
$\left\langle Z_\tau, \ch \left( \nabla^{{\cal D}^\omega} \right) 
\right\rangle$, or more precisely the pullback of the contribution to 
$M^\prime$, is
\begin{align}
& c_k 
\sum_{g_1 \ldots g_k} Z_\tau \left(
dg_1 \ldots dg_k  (g_k \ldots g_1)^{-1} \right)
\left( L_{(g_k \ldots g_1)^{-1}}^* d^{M^\prime} h \right) \ldots
\left( L_{g_{k}^{-1}}^* d^{M^\prime} h \right) = \notag \\
& c_k 
\sum_{g_1 \ldots g_k} Z_\tau \left(
(g_k \ldots g_1)^{-1} dg_1 \ldots dg_k  \right)
\left( L_{(g_k \ldots g_1)^{-1}}^* d^{M^\prime} h \right) \ldots
\left( L_{g_{k}^{-1}}^* d^{M^\prime} h \right). \notag
\end{align}
It is clear at this point that a nonzero contribution only arises when
$x = e$, in which case we get
\begin{align}
& c_k 
\sum_{g_1 \ldots g_k} \tau \left(
(g_k \ldots g_1)^{-1}, \ldots, g_k^{-1}, e \right)
\left( L_{(g_k \ldots g_1)^{-1}}^* d^{M^\prime} h \right) \ldots
\left( L_{g_{k}^{-1}}^* d^{M^\prime} h \right) = \\
& c_k 
\sum_{\gamma_1 \ldots \gamma_k} \tau \left(
\gamma_1, \ldots, \gamma_k, e \right)
\left( L_{\gamma_1}^* d^{M^\prime} h \right) \ldots
\left( L_{\gamma_k}^* d^{M^\prime} h \right).
\end{align}
One can show
\cite[Lemma 3]{Lott (1992)} 
that there is a closed form $\omega \in \Omega^k(M)$ such that
\begin{equation}
\sum_{\gamma_1 \ldots \gamma_k} \tau \left(
\gamma_1, \ldots, \gamma_k, e \right)
\left( L_{\gamma_1}^* d^{M^\prime} h \right) \ldots
\left( L_{\gamma_k}^* d^{M^\prime} h \right) = \pi^* \omega.
\end{equation}
Furthermore, the de Rham cohomology class $[\omega] \in \HH^k(M; \C)$ of 
$\omega$ satisfies \cite[Prop. 14]{Lott (1992)}
\begin{equation}
[\omega] = c_k \: \nu^* [\tau].
\end{equation}

We now argue that this is in fact the only nonzero contribution to
$\left\langle Z_\tau, \ch \left( \nabla^{{\cal D}^\omega} \right) 
\right\rangle$. First, looking at the group element factors in
$\left( \nabla^{{\cal D}^\omega} \right)^2$ and the structure of
$Z_\tau$, it is clear that $\left\langle Z_\tau, \ch 
\left( \nabla^{{\cal D}^\omega} \right) \right\rangle$
vanishes if $x \ne e$. Next, consider the possible contributions of the term
$\sum_{g, g^\prime \in 
\Gamma} dg \: dg^\prime \: h \left( L_{g}^* h 
\right) L_{g^\prime g}^*$ to 
$\left\langle Z_\tau, \Tr \left( \nabla^{{\cal D}^\omega} \right)^{2j} 
\right\rangle$.
For example, consider the case $j=1$. Then for $s \in C^\infty(M^\prime)$,
\begin{align}
& \sum_{g, g^\prime \in 
\Gamma} dg \: dg^\prime \: h \left( L_{g}^* h 
\right) L_{g^\prime g}^* s = \\
& \sum_{g, g^\prime \in 
\Gamma} dg \: dg^\prime \: h \left( L_{g}^* h 
\right) (g^\prime g)^{-1} \cdot s = \notag\\
& \sum_{g, g^\prime \in 
\Gamma} dg \: dg^\prime (g^\prime g)^{-1} L_{(g^\prime g)^{-1}}^* \left[ 
h \left( L_{g}^* h 
\right) \right] s = \notag \\
& \sum_{g, g^\prime \in 
\Gamma} dg \: dg^\prime (g^\prime g)^{-1} \left( L_{(g^\prime g)^{-1}}^*  
h \right) \left( L_{{g^\prime}^{-1}}^* h 
\right) s. \notag
\end{align}
If $\tau \in Z^2_e$ then
\begin{align}
& \left\langle Z_\tau, \Tr \left( \sum_{g, g^\prime \in 
\Gamma} dg \: dg^\prime \: h \left( L_{g}^* h 
\right) L_{g^\prime g}^* \right) 
\right\rangle = \\
&\sum_{g, g^\prime \in 
\Gamma} Z_\tau \left( dg \: dg^\prime (g^\prime g)^{-1} \right)
\left( L_{(g^\prime g)^{-1}}^*  
h \right) \left( L_{{g^\prime}^{-1}}^* h 
\right) \notag \\
&\sum_{g, g^\prime \in 
\Gamma} Z_\tau \left( (g^\prime g)^{-1} \: dg \: dg^\prime \right)
\left( L_{(g^\prime g)^{-1}}^*  
h \right) \left( L_{{g^\prime}^{-1}}^* h 
\right) \notag \\
&\sum_{g, g^\prime \in 
\Gamma} \tau \left( (g^\prime g)^{-1}, {g^\prime}^{-1},e \right)
\left( L_{(g^\prime g)^{-1}}^*  
h \right) \left( L_{{g^\prime}^{-1}}^* h 
\right) \notag \\
&\sum_{\gamma, \gamma^\prime \in 
\Gamma} \tau \left( \gamma, \gamma^\prime,e \right)
\left( L_{\gamma}^*  
h \right) \left( L_{\gamma^\prime}^* h 
\right).
\end{align}
Because of the antisymmetry of $\tau$, this vanishes.  A similar argument
using antisymmetry
applies to all terms in $\left\langle Z_\tau, 
\ch \left( \nabla^{{\cal D}^\omega} \right) 
\right\rangle$ involving $\sum_{g, g^\prime \in 
\Gamma} dg \: dg^\prime \: h \left( L_{g}^* h 
\right) L_{g^\prime g}^*$.
\end{pf}
{\bf Remark 1 : } There is a  
universal $\C \Gamma$-vector bundle ${\cal D}^0$ on $B\Gamma$.
Working simplicially \cite[Chapitre V]{Karoubi (1987)},
one can define a natural partially flat connection $\nabla^{{\cal D}^0}$ on 
${\cal D}^0$. Provided that one relaxes the regularity condition on $h$ to
being Lipschitz, one can realize $\nabla^{{\cal D}^\omega}$ as
$\nu^* \nabla^{{\cal D}^0}$, extended from $\C \Gamma$ to 
${\frak B}^\omega$.\\ \\
{\bf Remark 2 : } Let $C^*_r \Gamma$ denote the reduced group $C^*$-algebra
of $\Gamma$. Suppose that there is a Fr\'echet locally $m$-convex algebra 
${\frak B}^\infty$ such that\\
1. ${\frak B}^\omega \subset {\frak B}^\infty \subset C^*_r \Gamma$.\\
2. ${\frak B}^\infty$ is dense in $C^*_r \Gamma$.\\
3. ${\frak B}^\infty$ is
stable under the holomorphic functional calculus in $C^*_r \Gamma$.\\
We can complete ${\cal D}^\omega$ to a ${\frak B}^\infty$-vector bundle
${\cal D}^\infty$ on $M$ and to a 
$C^*_r \Gamma$-vector bundle
${\cal D}$ on $M$.  The latter represents an element $[{\cal D}]  \in 
K^{top}_{C^*_r \Gamma}(M) \cong
KK(\C, C(M) \otimes C^*_r \Gamma)$ and so gives a map
$\alpha : K_*(M) \rightarrow K_*(C^*_r \Gamma) \cong
K_*({\frak B}^\infty)$. Composing with the Chern character
gives $(\ch \circ \alpha)_\C : K_*(M) \otimes \C \rightarrow 
\overline{H}_*({\frak B}^\infty)$, which
coincides with the map coming from $\left[ 
\ch \left( \nabla^{{\cal D}^\infty} \right) \right]
\in \bigoplus_{p+q \: even} \HH^p \left(M; 
\overline{H}_q({\frak B}^\infty)\right)$. 
Suppose that for each $[\tau] \in
\HH^*(\Gamma; \C)$, there is a representative $\tau \in Z^*(\Gamma; \C)$
such that the cyclic cocycle $Z_\tau \in ZC^*(\C \Gamma)$ extends to a
continuous cyclic cocycle on
${\frak B}^\infty$. Proposition 
\ref{injective} shows that if $\nu^* [ \tau] \in \HH^*(M ;\C)$ is nontrivial
then $Z_\tau$ pairs nontrivially with
$\Image (\ch \circ \alpha)_\C$. 
Taking $M$ to be a sufficiently good approximation
to $B\Gamma$, we conclude that the
Strong Novikov Conjecture (SNC) holds for $\Gamma$,
meaning that the assembly map
$K_*(B \Gamma) \otimes \C \rightarrow K_*(C^*_r \Gamma) \otimes \C$ is
injective. The fact that the existence of ${\frak B}^\infty$ implies
SNC is well-known \cite[III.5]{Connes (1994)}, but we wish to emphasize how it
comes from the computation of $\ch \left( \nabla^{{\cal D}^\infty} \right)$.

If $\Gamma$ acts properly and cocompactly on a smooth manifold $X$ then
one can form the ${\frak B}^\omega$-vector bundle ${\frak B}^\omega
\times_\Gamma X$ on the orbifold $\Gamma/X$ and carry out a similar
analysis. 
The upshot is that if a finitely-generated discrete group $\Gamma$ satisfies
Hypothesis \ref{hypo} below 
then the Baum-Connes map \cite[II.10.$\epsilon$]{Connes (1994)} 
is rationally injective.\\ \\
{\bf Remark 3 : } In \cite{Lott (1992)} we gave a heat kernel proof of the
higher index theorem.  This proof can be reinterpreted using partially flat
connections.  For example, let $M$ be a even-dimensional
closed connected spin Riemannian
manifold and let $E$ be a Hermitian vector bundle on $M$ with Hermitian
connection $\nabla^E$. Then \cite[Prop. 12]{Lott (1992)} can be interpreted
as saying
\begin{equation} \label{higherindex}
\lim_{s \rightarrow 0} \left\langle Z_\tau,
 \TR_s \left( e^{- D_s^2} \right) \right\rangle = 
\int_M \widehat{A} \left(\nabla^{TM} \right) \wedge \ch \left( \nabla^E \right)
\wedge \left\langle Z_\tau,
\ch \left( \nabla^{{\cal D}^\omega} \right) \right\rangle.
\end{equation}
Here $\TR_s$ is the supertrace,
$s > 0$ is a factor which rescales the metric on $M$ and
$D_s$ denotes the (rescaled) Dirac operator on $M$, coupled to
${\cal E}^\omega$ using the connection $\nabla^{{\cal E}^\omega}$.
One can also prove (\ref{higherindex}) using the methods of
Section \ref{Fiber Bundles} of the present paper.

\section{${\frak B}$-Hermitian Metrics and Characteristic Classes}
\label{Hermitian Metrics and Characteristic Classes}
In this section we discuss the basic properties of a ${\frak B}$-valued
Hermitian metric on a ${\frak B}$-vector bundle.
We use such a Hermitian metric to define a
characteristic class of a ${\frak B}$-vector bundle with a flat structure.
(Related ideas occur in \cite[6.31-6.32]{Karoubi (1987)}.)
We show that the explicit partially flat connection described in the
previous section, in the context of covering spaces, is self-adjoint. 
We give an application to the
question of the rational injectivity of the algebraic K-theory assembly
map.
\subsection{${\frak B}$-Hermitian Metrics} \label{Hermitian Metrics}
Let $M$, ${\frak B}$ and ${\frak E}$
be as in Section \ref{Noncommutative Bundle Theory}. We assume that
$M$ is compact, possibly with boundary.
Suppose that ${\frak B}$ has an anti-involution, meaning a $\C$-antilinear
map $* : {\frak B} \rightarrow {\frak B}$ 
such that $(b_1 b_2)^* = b_2^* b_1^*$ and
$(b^*)^* = b$. We extend $*$ to $\Omega_*({\frak B})$ by requiring that
$(db)^* = - d(b^*)$. Let
$\overline{\frak E}^*$ be the vector space of $\C$-antilinear maps
$t : {\frak E} \rightarrow {\frak B}$ such that
$t(be) = t(e) b^*$ for all $b \in {\frak B}$ and $e \in {\frak E}$.
It is a left ${\frak B}$-module. 
If ${\cal E}$ is a
${\frak B}$-vector bundle on $M$ then there is an associated ${\frak B}$-vector
bundle $\overline{\cal E}^*$ such that $(\overline{\cal E}^*)_m =
\overline{{\cal E}_m}^*$. If ${\cal E}$ has a flat structure then so does
$\overline{\cal E}^*$. An element $t \in C^\infty(M; \overline{\cal E}^*)$
extends to a $\C$-antilinear map $t : \Omega(M, {\frak B}; {\cal E})
\rightarrow \Omega(M, {\frak B})$ such that 
\begin{equation}
t(\omega \otimes_{\frak B} e)
= t(e) \omega^*
\end{equation}
 for all $\omega \in \Omega_*({\frak B})$ and $e \in
C^\infty(M; {\cal E})$.
If $\nabla^{\cal E}$ is a connection on ${\cal E}$ then
there is an induced connection $\nabla^{\overline{\cal E}^*}$ on 
$\overline{\cal E}^*$ given by
\begin{equation} \label{funny}
d \left( t(e) \right) = \left( \nabla^{\overline{\cal E}^*} t \right)(e) -
t \left(\nabla^{\cal E} e \right) \in \Omega_1({\frak B})
\end{equation}
for all $t \in  C^\infty(M; \overline{\cal E}^*)$ and $e \in
C^\infty(M; {\cal E})$. (The funny sign in (\ref{funny}) comes from the
definition of the the involution on $\Omega_1({\frak B})$.)
If $\nabla^{\cal E}$ is partially flat then so is
$\nabla^{\overline{\cal E}^*}$.
\begin{definition}
1. A Hermitian form on 
${\frak E}$ is a map $\langle \cdot, \cdot \rangle : {\frak E} \times
{\frak E} \rightarrow {\frak B}$ which is $\C$-linear in the first variable,
$\C$-antilinear in the second variable and satisfies
$\langle b_1 e_1, b_2 e_2 \rangle = b_1 \langle e_1, e_2 \rangle
b_2^*$
for all $b_1, b_2 \in {\frak B}$ and $e_1, e_2 \in {\frak E}$.\\
2. A Hermitian form $\langle \cdot, \cdot \rangle$ is nondegenerate if
it induces an isomorphism $h^{\frak E} : {\frak E} \rightarrow 
\overline{\frak E}^*$
by $\left( h^{\frak E}(e_1) \right)(e_2) = \langle e_1, e_2 \rangle$.
\end{definition} 

We can extend $\langle \cdot, \cdot \rangle$ to a Hermitian form on
$\Omega_*({\frak B}) \otimes_{\frak B} {\frak E}$ by requiring that
\begin{equation} \label{formform}
\langle \omega_1 \otimes_{\frak B} e_1, \omega_2 \otimes_{\frak B} e_2 
\rangle = \omega_1 \langle e_1, e_2 \rangle \omega_2^* \in \Omega_*({\frak B})
\end{equation}
for all $\omega_1, \omega_2 \in \Omega_*({\frak B})$ and
$e_1, e_2 \in {\frak E}$.

There is a canonical Hermitian form $\langle \cdot, \cdot \rangle^0$
on ${\frak B}^n$ given by
$\langle \{x_i\}_{i=1}^n, \{y_i\}_{i=1}^n \rangle^0 =
\sum_{i=1}^n x_i y_i^*$.
 
\begin{definition}
A Hermitian metric on ${\frak E}$ is a Hermitian form 
$\langle \cdot, \cdot \rangle$ on ${\frak E}$ which is positive-definite, 
meaning that there is an embedding $i : {\frak E} \rightarrow
{\frak B}^n$ for some $n$ such that $\langle \cdot, \cdot \rangle =
i^* \langle \cdot, \cdot \rangle^0$.
\end{definition}
The method of proof of  
\cite[Lemme 2.7]{Karoubi (1980)} shows that a Hermitian metric is 
nondegenerate.

Since ${\frak E}$ is a finitely-generated projective ${\frak B}$-module, it
is clear that it admits some Hermitian metric.
The method of proof of  \cite[Lemme 2.9]{Karoubi (1980)} gives the following
proposition.
\begin{proposition} \label{K}
If $\langle \cdot, \cdot \rangle_0$ and $\langle \cdot, \cdot \rangle_1$
are Hermitian metrics on ${\frak E}$ then there is a
smooth $1$-parameter family $\{\alpha_t\}_{t \in [0,1]}$ in 
$\Aut_{{\frak B}}({\frak E})$ such that $\alpha_0 = \Id_{\frak E}$ and 
$\langle \cdot, \cdot \rangle_0 = \alpha_1^* \langle \cdot, \cdot \rangle_1$.
\end{proposition}
\begin{definition} 
A Hermitian metric on a ${\frak B}$-vector bundle ${\cal E}$
is given by a smooth family of Hermitian metrics on the fibers 
$\{{\cal E}_m\}_{m \in M}$. 
\end{definition}
\begin{proposition} \label{isotopy}
There is a Hermitian metric on a ${\frak B}$-vector bundle 
${\cal E}$.
Any two such Hermitian metrics are related by an automorphism which is
isotopic to the identity.
\end{proposition}
\begin{pf}
The algebra $C^\infty(M; {\frak B})$ is a Fr\'echet locally $m$-convex algebra
in a natural way.  Furthermore, $C^\infty(M; {\cal E})$ is a
finitely-generated projective $C^\infty(M; {\frak B})$-module. (The proof
is similar to that of the usual case when ${\frak B} = \C$, the essential
tool being that $\Inv(M_N({\frak B}))$ is open in $M_N({\frak B})$.)
A Hermitian
metric on the ${\frak B}$-vector bundle ${\cal E}$ is the same as a 
Hermitian metric on the $C^\infty(M; {\frak B})$-module
$C^\infty(M; {\cal E})$. The result now follows from Proposition
\ref{K}. 
\end{pf}

A Hermitian metric on ${\cal E}$ gives a $C^\infty(M; 
{\frak B})$-linear isomorphism
$h^{\cal E} : {\cal E} \rightarrow \overline{\cal E}^{*}$.

\begin{definition}
Given a connection $\nabla^{{\cal E}}$ on ${\cal E}$, its adjoint connection is
\begin{equation}
\left( \nabla^{{\cal E}} \right)^* = \left( h^{{\cal E}} 
\right)^{-1} \circ
\nabla^{\overline{\cal E}^*} \circ h^{{\cal E}},
\end{equation}
another connection on ${\cal E}$.
\end{definition}
Explicitly,
\begin{equation}
d \left\langle e_1, e_2 \right\rangle = \left\langle \nabla^{{\cal E}} 
e_1, e_2 \right\rangle -
\left\langle e_1, \left( \nabla^{{\cal E}} \right)^* e_2 \right\rangle
\in \Omega_1(M, {\frak B}).
\end{equation}
We say that $\nabla^{{\cal E}}$ is self-adjoint if
$\left( \nabla^{{\cal E}} \right)^* = \nabla^{{\cal E}}$.

Suppose that ${\cal E}$ has a flat structure. 
The triple $\left({\cal E}, \overline{\cal E}^*, h^{\cal E} \right)$
defines an element of $K^{rel}_{\frak B}(M)$. 
In this case, we write
\begin{equation}
CS \left( \nabla^{\cal E}, h^{\cal E} \right) = 
CS \left( \nabla^{\cal E}, \nabla^{\overline{\cal E}^*} \right)
\in \overline{\Omega}^{\prime \prime,odd}
(M, {\frak B}).
\end{equation}
Proposition
\ref{isotopy} implies that
$$\left[ CS\left(\nabla^{{\cal E}}, h^{\cal E} \right) \right] \in 
\bigoplus 
\begin{Sb}
p > q \\
p + q \: odd
\end{Sb}
\HH^p(M; \overline{H}_q({\frak B}))$$
only depends on the flat structure on ${\cal E}$. To put it another
way, the assignment of 
$\left({\cal E}, \overline{\cal E}^*, h^{\cal E} \right)$ to
${\cal E}$ gives an explicit map $K^{alg}_{\frak B}(M) \rightarrow
K^{rel}_{\frak B}(M)$. We can then apply $CS$ to obtain an invariant of
$K^{alg}_{\frak B}(M)$. In total, we have defined a map
\begin{align} \label{chcs}
\ch \oplus CS : K^{alg}_{\frak B}(M) \rightarrow
& \left(\bigoplus_p \HH^p(M; \overline{Z}_p({\frak B})) \right)
\oplus
\left(\bigoplus
\begin{Sb}
p < q \\
p + q \: even
\end{Sb}
\HH^p(M; \overline{H}_q({\frak B}))\right) \oplus \\
& \left( \bigoplus 
\begin{Sb}
p > q \\
p + q \: odd
\end{Sb}
\HH^p(M; \overline{H}_q({\frak B})) \right). \notag
\end{align}

Let $(M, *)$ and $(M^\prime, *^\prime)$ be smooth connected manifolds with
basepoints.  Let ${\cal E}$ be a ${\frak B}$-vector
bundle on $M$ and similarly for ${\cal E}^\prime$. 
Let ${\cal T} = M \times {\cal E} \big|_{*}$ 
denote the trivial ${\frak B}$-vector bundle on $M$ with the same fiber
at $*$ as ${\cal E}$, and similarly for ${\cal T}^\prime$. 
Then $[{\cal E}] - [{\cal T}]$ is an element of
the reduced group $\widetilde{K}^{top}_{{\frak B}}(M)$, 
and similarly for
$[{\cal E}^\prime] - [{\cal T}^\prime]$. The
virtual ${\frak B} \otimes_\C {\frak B}^{\prime}$-vector bundle
${\cal E} \otimes_\C {\cal E}^\prime - {\cal T} \otimes_\C {\cal E}^\prime
- {\cal E} \otimes_\C {\cal T}^\prime + {\cal T} \otimes_\C {\cal T}^\prime$ on
$M \times M^\prime$ is trivial on $(M \times \{*^\prime\}) \cup 
(\{*\} \times M^\prime)$ and so
passes to an element of 
$\widetilde{K}^{top}_{{\frak B} \otimes_\C {\frak B}^{\prime}}
(M \wedge M^\prime)$ which represents the product
$\left( [{\cal E}] - [{\cal T}] \right) \cdot
\left( [{\cal E}^\prime] - [{\cal T}^\prime] \right)$. 
If ${\cal E}$ and ${\cal E}^\prime$ have flat structures then we get the
product in $\widetilde{K}^{alg}$ \cite[Chapitre II]{Loday (1976)}.

Let $\nabla^{\cal E}$ and $\nabla^{{\cal E}^\prime}$ be connections on
${\cal E}$ and ${\cal E}^\prime$, respectively.  There are induced
connections on ${\cal T}$ and ${\cal T}^\prime$.
Let $z \in Z_*(M,*;\C)$ and $z^\prime \in 
Z_*(M^\prime,*^\prime;\C)$ be relative cycles. Let $z z^\prime \in Z_*(M\wedge
M^\prime,*;\C)$ be the product. Then with an obvious notation,
(\ref{chmull}) implies that
\begin{equation} \label{chmul}
\int_{zz^\prime}\ch \left( \nabla^{\left( {\cal E} - {\cal T} \right) \cdot
\left( {\cal E}^\prime - {\cal T}^\prime \right)} \right) =
\int_z \ch \left( \nabla^{{\cal E} - {\cal T} } \right) \cdot
\int_{z^\prime} \ch \left( \nabla^{
{\cal E}^\prime - {\cal T}^\prime } \right).
\end{equation}
Suppose that ${\cal E}$ and ${\cal E}^\prime$ have  partially flat connections
and Hermitian metrics. Suppose that 
$\nabla^{{\cal E}^\prime}$ is self-adjoint. Then from (\ref{csmul}),
\begin{equation} \label{chcsmul}
\int_{z z^\prime} CS \left( \nabla^{\left( {\cal E} - {\cal T} \right) \cdot
\left( {\cal E}^\prime - {\cal T}^\prime \right)}, 
h^{\left( {\cal E} - {\cal T} \right) \cdot
\left( {\cal E}^\prime - {\cal T}^\prime \right)}
\right) =
\int_z CS \left( \nabla^{{\cal E} - {\cal T}}, 
h^{{\cal E} - {\cal T}} \right) \cdot
\int_{z^\prime} \ch \left( \nabla^{
{\cal E}^\prime - {\cal T}^\prime } \right).
\end{equation}
\subsection{${\frak B}$-Hermitian Metrics,  Group Algebras and Assembly Maps}
\label{Hermitian Metrics and Group Algebras}
We use the notation of Subsection \ref{Covering Spaces}. Define an
involution on ${\frak B}^\omega$ by $* \left( \sum_{g \in \Gamma} c_g g \right)
= \sum_{g \in \Gamma} \overline{c_g} g^{-1}$. Considering ${\frak B}^\omega$
as a left-module over itself, it has a Hermitian form given by
$\langle b_1, b_2 \rangle = b_1 b_2^*$. We can transfer this Hermitian
form fiberwise to ${\cal D}^\omega$. In what follows, we will freely
identify differential forms on $M$ and $\Gamma$-invariant differential
forms on $M^\prime$.
\begin{definition}
Given $s_1, s_2 \in C^\infty(M; {\cal D}^\omega)$, consider them to be
elements of $C^\infty(M^\prime)$ by (\ref{srep}). Then the Hermitian form
$\langle s_1, s_2 \rangle^\omega \in C^\infty(M; {\frak B}^\omega)$ is given by
\begin{align}
\langle s_1, s_2 \rangle^\omega(x) & = \sum_{g,g^\prime \in \Gamma} g \:
(L_{g g^\prime}^* s_1)(x) \: (L_{g^\prime}^*\overline{s_2})(x) \\
& = \sum_{g,g^\prime \in \Gamma} g \:
s_1(g g^\prime x) \: \overline{s_2}(g^\prime x). \notag
\end{align}
\end{definition}
We do not claim that $\langle \cdot, \cdot \rangle^\omega$ is a
Hermitian metric, in that it may be degenerate.
\begin{proposition} \label{sa}
$\nabla^{{\cal D}^\omega}$ is self-adjoint with respect to $\langle \cdot,
\cdot \rangle^\omega$, meaning that for all $s_1, s_2 \in C^\infty(M; 
{\cal D}^\omega)$,
\begin{equation} \label{eq1}
d \left\langle s_1, s_2 \right\rangle^\omega = \left\langle 
\nabla^{{\cal D}^\omega} 
s_1, s_2 \right\rangle^\omega -
\left\langle s_1, \nabla^{{\cal D}^\omega} s_2 \right\rangle^\omega
\in \Omega_1(M, {\frak B}^\omega).
\end{equation}
\end{proposition}
\begin{pf}
As $\Gamma$-invariant differential forms on $M^\prime$, we have
\begin{equation} \label{eq2}
d \left\langle s_1, s_2 \right\rangle^\omega 
= \sum_{g,g^\prime} \left[ dg \left(
L_{g g^\prime}^* s_1 \right) L_{g^\prime}^*\overline{s_2} + g \left(
d^{M^\prime} L_{g g^\prime}^* s_1 \right) 
L_{g^\prime}^*\overline{s_2} + g \left(
L_{g g^\prime}^* s_1 \right) d^{M^\prime} L_{g^\prime}^*\overline{s_2} \right]
\end{equation}
and
\begin{align} \label{eq3}
\left\langle \nabla^{{\cal D}^\omega,1,0} s_1, s_2 \right\rangle^\omega & =
\left\langle d^{M^\prime} s_1, s_2 \right\rangle^\omega \\
& = \sum_{g, g^\prime} g \left( L_{g g^\prime}^* d^{M^\prime} s_1 
\right)
L_{g^\prime}^*\overline{s_2}, \notag
\end{align}
\begin{align} \label{eq3.5}
\left\langle \nabla^{{\cal D}^\omega,0,1} s_1, s_2 \right\rangle^\omega & = 
\left\langle \sum_\gamma d\gamma \left( h L_\gamma^* 
s_1 \right), s_2 \right\rangle^\omega \\
& = \sum_{g, g^\prime, \gamma} 
d\gamma \: g \: L_{g g^\prime}^*
\left( h L_\gamma^* s_1 \right) L_{g^\prime}^*\overline{s_2}  \notag \\
& = \sum_{g, g^\prime, \gamma} 
\left[ d(\gamma g) - \gamma \: dg \right] \: L_{g g^\prime}^*
h \left( L_{\gamma g g^\prime}^* s_1 \right) L_{g^\prime}^*\overline{s_2} 
\notag \\
& = \sum_{g, g^\prime, \gamma} 
dg \: L_{\gamma^{-1} g g^\prime}^*
h \left( L_{g g^\prime}^* s_1 \right) L_{g^\prime}^*\overline{s_2}  -
\sum_{g, g^\prime, \gamma} 
\gamma \: dg \: L_{g g^\prime}^*
h \left( L_{\gamma g g^\prime}^* s_1 \right) L_{g^\prime}^*\overline{s_2}. 
\notag \\
& = \sum_{g, g^\prime} 
dg \:
\left( L_{g g^\prime}^* s_1 \right) L_{g^\prime}^*\overline{s_2} -
\sum_{g, g^\prime, \gamma} 
\gamma \: dg \: L_{g g^\prime}^*
h \left( L_{\gamma g g^\prime}^* s_1 \right) L_{g^\prime}^*\overline{s_2}. 
\notag
\end{align}
Switching $s_1$ and $s_2$ gives
\begin{equation}
\left\langle \nabla^{{\cal D}^\omega,1,0} s_2, s_1 \right\rangle^\omega
 = \sum_{g, g^\prime} g \left( L_{g g^\prime}^* d^{M^\prime} s_2 
\right) L_{g^\prime}^*\overline{s_1}
\end{equation}
and
\begin{equation}
\left\langle \nabla^{{\cal D}^\omega,0,1} s_2, s_1 \right\rangle^\omega =
\sum_{g, g^\prime, \gamma} 
d\gamma \: g \: L_{g g^\prime}^*
\left( h L_\gamma^* s_2 \right) L_{g^\prime}^*\overline{s_1}.
\end{equation}
Then
\begin{align} \label{eq4}
\left\langle s_1, \nabla^{{\cal D}^\omega,1,0} s_2 \right\rangle^\omega 
& = - \sum_{g, g^\prime} g^{-1} \left( L_{g^\prime}^*{s_1} \right)
L_{g g^\prime}^* d^{M^\prime} \overline{s_2} \\
& = - \sum_{g, g^\prime} g  \left( L_{g^\prime}^*{s_1} \right)
L_{g^{-1} g^\prime}^* d^{M^\prime} \overline{s_2} \notag \\
& = - \sum_{g, g^\prime}  g \left( L_{g g^\prime}^*{s_1} \right)
L_{g^\prime}^* d^{M^\prime} \overline{s_2} \notag
\end{align}
and
\begin{align} \label{eq5}
\left\langle s_1, \nabla^{{\cal D}^\omega,0,1} s_2 \right\rangle^\omega & =
 - \sum_{g, g^\prime, \gamma} g^{-1} \: 
d{\gamma}^{-1} \: L_{g^\prime}^*s_1 \left( L_{g g^\prime}^* h \right)
L_{\gamma g g^\prime}^* \overline{s_2} \\
& = - \sum_{g, g^\prime, \gamma} \gamma \: 
dg \: L_{g^\prime}^*s_1 \left( L_{\gamma^{-1} g^\prime}^* h \right)
L_{g^{-1} \gamma^{-1} g^\prime}^* \overline{s_2} \notag \\
& = - \sum_{g, g^\prime, \gamma} \gamma \: 
dg \: L_{\gamma g g^\prime}^*s_1 \left( L_{g g^\prime}^* h \right)
L_{g^\prime}^* \overline{s_2}. \notag
\end{align}
Combining (\ref{eq2}), (\ref{eq3}), (\ref{eq3.5}), (\ref{eq4}) and 
(\ref{eq5}) gives (\ref{eq1}). 
\end{pf} 

We now give examples in which the map (\ref{chcs}) is nontrivial.
Recall the notation of Subsection \ref{Cyclic Cohomology of Group Algebras}.
\begin{hypothesis} \label{hypo}
There is an involutive Fr\'echet locally $m$-convex algebra ${\frak B}^\infty$
such that\\
1. ${\frak B}^\omega \subset {\frak B}^\infty \subset C^*_r \Gamma$.\\
2. ${\frak B}^\infty$ is dense in $C^*_r \Gamma$ and stable under the
holomorphic functional calculus in $C^*_r \Gamma$.\\
3. For each $\langle x \rangle \in 
\langle \Gamma \rangle^{\prime}$ and 
$[\tau] \in \HH^*(N_{\langle x \rangle}; \C)$, 
there is a representative
$\tau \in Z^*_x$ such that the cyclic cocycle $Z_\tau
\in ZC^*(\C \Gamma)$ extends to a continuous cyclic cocycle on
${\frak B}^\infty$. 
\end{hypothesis}

Hypothesis \ref{hypo} is known to
be satisfied by virtually nilpotent groups \cite{Ji (1992)} and 
Gromov-hyperbolic groups \cite{Ogle (1992)}. 

We can extend ${\cal D}^\omega$ to a ${\frak B}^\infty$-vector bundle
${\cal D}^\infty$ on $M$, $\nabla^{{\cal D}^\omega}$ to a connection
$\nabla^{{\cal D}^\infty}$ on ${\cal D}^\infty$ and $\langle \cdot, \cdot
\rangle^\omega$ to a Hermitian metric $\langle \cdot, \cdot \rangle^\infty$ on
${\cal D}^\infty$.

By \cite[Chapitre III]{Karoubi (1987)}, 
we can consider an element $k$ of $K_n^{alg}(\Z \Gamma)$ to
be given by a formal difference of homology $n$-spheres 
$HS^n$ equipped with flat
bundles ${\cal E}^0$ of finitely-generated projective $\Z \Gamma$-modules.
For simplicity, we just consider a single $HS^n$.  
As in \cite[p. 98]{Karoubi (1987)}, we may approximate $HS^n$ by
compact manifolds (possibly with boundary), so for simplicity we
assume that $HS^n$ is a
compact manifold.  Let $[HS^n] \in 
\HH^n(HS^n; \C)$ denote its ``fundamental class''. Putting 
${\cal E}^\infty = {\frak B}^\infty \otimes_{\Z \Gamma} {\cal E}^0$, the
algebraic K-theory class of
$\left[{\cal E}^\infty \right]$ represents the image
of $k$ under the map $K_n^{alg}(\Z \Gamma) \rightarrow 
K_n^{alg} \left({\frak B}^\infty \right)$.
Now apply $(\ref{chcs})$ to $\left[{\cal E}^\infty \right]$, 
pair the result with $Z_\tau$ and
integrate over $\left[ HS^n \right]$ to get a number.
This procedure gives a map
\begin{equation} \label{procedure}
K_n^{alg}(\Z \Gamma) \otimes_\Z \C \rightarrow \bigoplus_{\langle x \rangle
\in \langle \Gamma^\prime \rangle} \left( \HH_n(N_{\langle x \rangle}; \C) 
\oplus \left( \bigoplus_{k=0}^{\left[ \frac{n-1}{2} \right]} 
\HH_{n-1-2k}(N_{\langle x \rangle}; \C) \right) \right)
\end{equation}
which is conjecturally injective. \\ \\
{\bf Example 2 :} Take $\Gamma$ to be the trivial group and 
${\cal B}^\infty = \C$. In this case, (\ref{chcs}) becomes
\begin{equation}
ch \oplus CS : K^{alg}_{\C}(M) \rightarrow \HH^0(M; \C) \oplus \left(
\bigoplus_{p \: odd} \HH^p(M; \C) \right).
\end{equation}
Applied to a flat complex vector bundle $E$ on $M$, this represents
the rank of $E$ along with its Borel classes 
\cite[Section Ig]{Bismut-Lott (1995)}. It is known that 
\begin{equation}
K^{alg}_n(\Z) \otimes \C = 
\begin{cases}
\C & \text{ if } n = 0 \\
\C & \text{ if } n \equiv 1 \text{ mod } 4, \: n > 1\\
0 & \text{ otherwise,}
\end{cases}
\end{equation}
with the higher terms being detected by the Borel classes
\cite{Borel (1974)}.
Thus for all $n \equiv 1 \mod 4, \: n > 1$, there 
is a homology $n$-sphere $HS^n$ and a flat bundle ${\cal E}^0$
of finitely-generated projective $\Z$-modules on $HS^n$ such that
if ${\cal E} = \C \otimes_\Z {\cal E}^0$ then
$\int_{[HS^n]} CS \left( {\cal E} \right) \neq 0$. \\ \\
{\bf Example 3 : } Take $\Gamma$ to be a finite group and ${\frak B}^\infty
= \C \Gamma$. We can write $\C \Gamma = \bigoplus_{\rho_i \in \widehat{\Gamma}}
M_{n_i}(\C)$, where $n_i = \dim(\rho_i)$. Then (\ref{chcs}) becomes
\begin{equation}
ch \oplus CS : K^{alg}_{\C \Gamma}(M) 
\rightarrow \bigoplus_{\rho \in \widehat{\Gamma}}
\left[ \HH^0(M; \C) \oplus \left(
\bigoplus_{p \: odd} \HH^p(M; \C) \right) \right].
\end{equation}

Consider the case $n = 1$ of (\ref{procedure}). Take
the homology sphere to be a circle $S^1$. Given $T \in 
GL_r(\Z \Gamma)$, form a flat $\Z \Gamma$-bundle  ${\cal E}^0$
on $S^1$ by gluing the ends of $[0,1] \times \left(\Z \Gamma\right)^r$ 
using $T$.  Then
${\cal E}^\infty = \C \Gamma \otimes_{\Z \Gamma} {\cal E}^0$ 
is a flat $\C \Gamma$-vector bundle on $S^1$ with holonomy
$T$. One computes that 
\begin{equation}
\int_{S^1} CS \left(\left[{\cal E}^\infty\right] \right) = 
\bigoplus_{\rho \in \widehat{\Gamma}} 
\: 2 \: \log |\det(\rho(T))|.
\end{equation}
Thus $CS$ detects all of $K^{alg}_1(\Z \Gamma) \otimes \C$ 
\cite{Oliver (1988)}.\\ \\
{\bf Example 4 : } 
Suppose that $\Gamma$ satisfies Hypothesis \ref{introhypo} of the
introduction.
For an algebra $A$, let ${\bf K}_A$ denote the algebraic K-theory spectrum
of $A$, with $({\bf K}_A)_0 = K_0(A) \times BGL(A)_\delta^+$.
For an abelian group $G$, let ${\bf H}_G$ denote the Eilenberg-Maclane
spectrum of $G$.
We can think of the class $CS$ as arising from a map
\begin{equation}
{\bf K}_{{\frak B}^\infty} \stackrel{CS}{\longrightarrow}
\prod
\begin{Sb}
p > q \\
p + q \: odd
\end{Sb}
\Sigma^p({\bf H}_{\overline{H}_q({\frak B}^\infty)}).
\end{equation}

We recall the assembly map of \cite[Chapitre IV]{Loday (1976)}, 
extended from $\Z \Gamma$
to ${\frak B}^\infty$.
The inclusion of $\Gamma$ into $GL(\Z \Gamma)$, as a matrix with one
nonzero entry in the upper left corner, induces maps 
\begin{equation} \label{ass1}
B\Gamma \rightarrow BGL(\Z \Gamma)_\delta \rightarrow 
BGL(\Z \Gamma)_\delta^+ \rightarrow BGL\left({\frak B}^\infty\right)_\delta^+,
\end{equation} 
which extend to
\begin{equation}
B\Gamma \rightarrow K_0(\Z \Gamma) \times BGL(\Z \Gamma)_\delta \rightarrow 
K_0(\Z \Gamma) \times BGL(\Z \Gamma)_\delta^+ \rightarrow 
K_0 \left({\frak B}^\infty\right) 
\times BGL\left({\frak B}^\infty\right)_\delta^+.
\end{equation}
Smashing with ${\bf K}_{\Z}$ gives the assembly map
\begin{equation} \label{ass2}
{\bf K}_\Z \wedge B\Gamma \rightarrow {\bf K}_{\Z \Gamma}
\rightarrow {\bf K}_{{\frak B}^\infty}.
\end{equation}

We can compose with $CS$ to get
\begin{equation} \label{ass3}
{\bf K}_\Z \wedge B\Gamma \rightarrow {\bf K}_{\Z \Gamma}
\rightarrow {\bf K}_{{\frak B}^\infty} \rightarrow
\prod
\begin{Sb}
p > q \\
p + q \: odd
\end{Sb}
\Sigma^p ( {\bf H}_{\overline{H}_q({\frak B}^\infty)} ).
\end{equation}
Taking homotopy groups gives
\begin{equation}
\widetilde{\HH}_p \left(B\Gamma; {\bf K}_\Z \right) \rightarrow 
K_{p}(\Z \Gamma)
\rightarrow K_p({\frak B}^\infty) \rightarrow
\bigoplus
\begin{Sb}
q<p \\
q + p \: odd
\end{Sb}
\overline{H}_q({\frak B}^\infty).
\end{equation}
Then tensoring with $\C$, we obtain 
\begin{equation} \label{tensoring}
\bigoplus_{l+m=p} (K_l(\Z) \otimes_{\Z} \C) 
\otimes_\C \widetilde{\HH}_m(\Gamma; \C) \rightarrow 
K_{p}(\Z \Gamma) \otimes_\Z \C
\rightarrow K_p({\frak B}^\infty) \otimes_\Z \C \rightarrow
\bigoplus
\begin{Sb}
q<p \\
q + p \: odd
\end{Sb}
\overline{H}_q({\frak B}^\infty).
\end{equation}

Let $HS^l$ and ${\cal E}$ be as in Example 2. 
Let $[{\cal E}]_{\C} \in K_l(\Z) \otimes_{\Z} \C$ be the corresponding 
K-theory class.
Then (\ref{tensoring}) gives a map
\begin{equation} \label{mapp}
[{\cal E}]_{\C} : \widetilde{\HH}_{p-l}(\Gamma; \C) \rightarrow
K_{p}(\Z \Gamma) \otimes_\Z \C
\rightarrow K_p({\frak B}^\infty) \otimes_\Z \C \rightarrow
\bigoplus
\begin{Sb}
q<p \\
q + p \: odd
\end{Sb}
\overline{H}_q({\frak B}^\infty). 
\end{equation}
Tracing through the definitions, (\ref{mapp}) 
can be interpreted concretely as follows. 
Assume that $B\Gamma$ is a manifold, possibly with boundary.
(Otherwise, approximate it by manifolds.) Let ${\cal E}^\prime$ be
the ${\frak B}^\infty$-vector bundle ${\cal D}^\infty$ on
$B \Gamma$. The map (\ref{ass1}) is realized
geometrically by $[{\cal E}^\prime] - [{\cal T}^\prime]
\in \widetilde{K}^{alg}_{{\frak B}^\infty}(B\Gamma)$, where  
${\cal T}^\prime$ is as in the discussion above equation (\ref{chmul}). 
Consider the
${\frak B}^\infty$-vector bundle $({\cal E} - {\cal T}) \otimes_\C
({\cal E}^\prime - {\cal T}^\prime)$ on $HS^l \wedge B\Gamma$.
Then
\begin{align}
CS \left(\left[({\cal E} - {\cal T}) \otimes_\C
({\cal E}^\prime - {\cal T}^\prime )\right]\right) \in 
& \bigoplus 
\begin{Sb}
q < p \\
q+p \: odd
\end{Sb}
\widetilde{\HH}^p(HS^l \wedge B \Gamma; \overline{H}_q({\frak B}^\infty)) \\
& =
\bigoplus 
\begin{Sb}
q < p \\
q+p \: odd
\end{Sb}
\widetilde{\HH}^{p-l}(\Gamma; \C) \otimes \overline{H}_q({\frak B}^\infty). 
\notag
\end{align}
The map (\ref{mapp}) comes from pairing $\widetilde{\HH}_{p-l}(\Gamma; \C)$ 
with
$CS \left(\left[({\cal E} - {\cal T}) \otimes_\C
({\cal E}^\prime - {\cal T}^\prime )\right]\right)$.

From Proposition \ref{sa}, $\nabla^{{\cal E}^\prime}$ is self-adjoint. 
Thus we can apply (\ref{chcsmul}), with $z \in Z_{l}(HS^l, *; \C)$ and 
$z^\prime \in Z_{p-l}(B\Gamma, *^\prime; \C)$. 
If $l \equiv 1 \mod 4$  and $l > 1$ then by Example 2,  
we can choose $HS^l$ and
${\cal E}$ so that $\int_z CS \left( \nabla^{{\cal E} - {\cal T}}, 
h^{{\cal E} - {\cal T}} \right)$ is nonzero. Applying Proposition
\ref{injective} to
$\int_{z^\prime} \ch \left( \nabla^{
{\cal E}^\prime - {\cal T}^\prime } \right)$, we conclude that
(\ref{tensoring}) gives an injection
\begin{align}
& \bigoplus_{k=1}^{\left[ \frac{p-1}{4}\right]} 
\widetilde{\HH}_{p-1-4k}(\Gamma; \C) 
\rightarrow K_p(\Z \Gamma) \otimes_\Z \C 
\rightarrow K_p({\frak B}^\infty) \otimes_\Z \C
\rightarrow \\
& \hspace{2in} \bigoplus_{k=0}^{\left[ \frac{p-1}{2}\right]} 
\overline{H}_{p-1-2k}({\frak B}^\infty) \rightarrow 
\bigoplus_{k=0}^{\left[ \frac{p-1}{2}\right]} 
\HH_{p-1-2k}(\Gamma; \C). \notag
\end{align}
Thus
\begin{equation}
\bigoplus_{k=1}^{\left[ \frac{p-1}{4}\right]} 
\widetilde{\HH}_{p-1-4k}(\Gamma; \C) 
\rightarrow K_p(\Z \Gamma) \otimes_\Z \C
\end{equation}
is injective.
Including the contribution of the characteristic class $\ch$ and taking 
more care with reduced vs. unreduced homology gives the injectivity of
\begin{equation} \label{ratinj}
\HH_p(\Gamma; \C) \oplus \left(
\bigoplus_{k=1}^{\left[ \frac{p-1}{4}\right]} \HH_{p-1-4k}(\Gamma; \C) \right)
\rightarrow
K_p(\Z \Gamma) \otimes_\Z \C.
\end{equation}
This is not an optimal result, as (\ref{ratinj}) is known to be
injective for all groups $\Gamma$ such that $\dim_\C \HH_k(\Gamma; \C) <
\infty$ for all $k \in \N$, regardless of whether or not they satisfy
Hypothesis \ref{introhypo} \cite{Boekstedt-Hsiang-Madsen (1993)}.
The proof of \cite{Boekstedt-Hsiang-Madsen (1993)} uses more
complicated methods. 

There is a conjecture that (\ref{ratinj}) is an isomorphism if $\Gamma$ is
torsion-free. This is known to be true when $\Gamma$ is a discrete
cocompact subgroup of a Lie group with a finite number of connected
components
\cite{Farrell-Jones (1993)}.
\section{Noncommutative Superconnections}
In this section we first extend the results of the preceding sections from
connections to superconnections. For basic information about 
the superconnection formalism,
we refer to the book \cite{Berline-Getzler-Vergne (1992)}.
We then use superconnections to prove
a finite-dimensional analog of our fiber-bundle index theorem.  We also
construct the associated finite-dimensional analytic torsion form and relate
it to various versions of the Reidemeister torsion. The main technical
problems of this section involve the large-time behaviour of
heat kernels in Fr\'echet spaces.
\label{Noncommutative Superconnections}
\subsection{Partially Flat Superconnections}
\label{Partially Flat Superconnections}
Let $M$, ${\frak B}$ and ${\cal E}$ be as in Subsection 
\ref{Hermitian Metrics}.
Suppose that ${\cal E}$ is $\Z$-graded
as a direct sum
\begin{equation}
{\cal E} = \bigoplus_{i=1}^n {\cal E}^i
\end{equation}
of ${\frak B}$-vector bundles on $M$. We use the induced $\Z_2$-grading
on ${\cal E}$ when defining supertraces. The algebra $\Omega \left(
M; \Hom_{\frak B} ({\cal E}, \Omega_*({\frak B}) \otimes_{\frak B} 
{\cal E} ) \right)$ has a trigrading as 
\begin{equation}
\Omega \left(
M; \Hom_{\frak B} ({\cal E}, \Omega_*({\frak B}) \otimes_{\frak B} 
{\cal E} )\right) =
\bigoplus
\begin{Sb}
p,q,r \in \Z\\
p,q \ge 0
\end{Sb}
\Omega^{p,q,r}(M, {\frak B}, \End({\cal E})),
\end{equation}
where by definition
\begin{equation}
\Omega^{p,q,r}(M, {\frak B}, \End({\cal E})) = 
\Omega^p \left(
M; \Hom_{\frak B} ({\cal E}^\bullet, \Omega_q({\frak B}) \otimes_{\frak B}
{\cal E}^{\bullet + r} ) \right).
\end{equation}
Define a subalgebra of $\Omega \left(
M; \Hom_{\frak B} ({\cal E}, \Omega_*({\frak B}) \otimes_{\frak B} 
{\cal E} )\right)$ by
\begin{equation} \label{omegaprime}
\Omega^\prime \left(
M; \Hom_{\frak B} ({\cal E}, \Omega_*({\frak B}) \otimes_{\frak B} 
{\cal E} )\right) = \bigoplus_{p + r \le q} \Omega^{p,q,r}
(M, {\frak B}, \End({\cal E})).
\end{equation} 
\begin{definition}
A degree-$1$ superconnection $A^\prime$ on ${\cal E}$ is a sum
\begin{equation} \label{supercomp}
A^\prime = \sum_{p \ge 0} A_p^\prime
\end{equation}
where\\
1. $A_1^\prime$ is a connection $\nabla^{\cal E}$
on ${\cal E}$ which preserves the
$\Z$-grading.\\
2. For $p \neq 1$, $A_p^\prime \in \bigoplus_{p + q + r = 1}  
\Omega^{p,q,r}(M, {\frak B}, \End({\cal E}))$.
\end{definition}

We will sometimes omit the phrase ``degree-$1$''.
For $p \neq 1$, let $A^\prime_{p,q,r}$ denote the component of
$A_p^\prime$ in $\Omega^{p,q,r}(M, {\frak B}, {\cal E})$. As in 
(\ref{conndecomp}), we write 
$\nabla^{\cal E} = 
\nabla^{{\cal E},1,0} \oplus 
\nabla^{{\cal E},0,1}$.

The superconnection $A^\prime$ gives a $\C$-linear map
\begin{equation}
A^\prime : C^\infty(M; {\cal E}) \rightarrow \Omega^*(M, {\frak B};
 {\cal E})
\end{equation}
which satisfies the Leibniz rule.  We extend $A^\prime$ to a $\C$-linear
map on $\Omega^*(M, {\frak B}; {\cal E})$
by requiring that for all $\omega \in \Omega^k(M, {\frak B})$ and $s \in
\Omega^{l}(M, {\frak B}; {\cal E})$, 
\begin{equation}
\nabla^{\cal E}(\omega s) = (-1)^k \: \omega \wedge 
\nabla^{\cal E} s + d\omega 
\otimes_{C^\infty(M; {\frak B})} s. 
\end{equation}
The curvature of $A^\prime$ is
\begin{equation}
\left( A^\prime \right)^2 \in \bigoplus_{p+q+r=2} 
\Omega^{p,q,r}(M, {\frak B}, \End({\cal E})).
\end{equation}
Let $\left( A^\prime \right)^2_{p,q,r}$ denote the component of
$\left( A^\prime \right)^2$ in $\Omega^{p,q,r}
(M, {\frak B}, \End({\cal E}))$.

The Chern character of $A^\prime$ is
\begin{equation}
\ch \left( A^\prime \right) = \Tr_s \left( e^{-
\left(A^\prime\right)^2}\right) \in
\overline{\Omega}^{even}(M, {\frak B}).
\end{equation}
It is a closed form whose cohomology class $\left[\ch \left( A^\prime \right)
\right] \in H^{even}_{\frak B}(M)$ is independent of the choice of $A^\prime$.
\begin{definition} \label{superparflat}
The superconnection $A^\prime$ is partially flat if
$\left(A^\prime \right)^2_{p,0,2-p} = 0$ for all $p \ge 0$.
\end{definition}
\begin{definition}
A superflat structure on ${\cal E}$ is given by a degree-$1$ superconnection
\begin{equation} \label{superflateqn}
A^{\prime,flat} : C^\infty(M; {\cal E}) \rightarrow
\Omega^*(M; {\cal E})
\end{equation}
which is ${\frak B}$-linear and whose extension to $\Omega^*(M; {\cal E})$
satisfies $\left(A^{\prime,flat} \right)^2 = 0$.
\end{definition}

Note that the map in (\ref{superflateqn}) does not involve any
${\frak B}$-differentiation.
A partially flat 
superconnection determines a superflat structure on ${\cal E}$
by 
\begin{equation}
A^{\prime, flat} = A^\prime_{0,0,1} + \nabla^{{\cal E},1,0} + 
\sum_{p=2}^\infty A^\prime_{p,0,1-p}.
\end{equation}
Conversely, given a superflat structure on ${\cal E}$,
there is a partially flat superconnection on ${\cal E}$ 
which is compatible with the superflat structure, 
although generally not a unique one.\\ \\
{\bf Example 5 :} If ${\frak B} = \C$ then a partially flat degree-$1$
superconnection on ${\cal E}$ is the same as a flat degree-$1$ 
superconnection on ${\cal E}$ in the
sense of \cite[Section IIa]{Bismut-Lott (1995)}.\\ \\
{\bf Example 6 :} If ${\cal E}$ is concentrated in degree $0$ then
a partially flat degree-$1$ superconnection on ${\cal E}$ is the same
as a partially flat connection on ${\cal E}$ in the sense of
Subsection \ref{Connections and Chern Character}.\\

Hereafter, we assume that $A^\prime$ is a partially flat degree-$1$
superconnection.
\begin{proposition}
We have $\ch( A^\prime ) \in \overline{\Omega}^{\prime, even}(M, {\frak B})$.
\end{proposition}
\begin{pf}
As $\left(A^\prime \right)^2 \in \Omega^\prime \left(
M; \Hom_{\frak B} ({\cal E}, \Omega_*({\frak B}) \otimes_{\frak B} 
{\cal E} )\right)$, the same is true for $e^{-
\left(A^\prime\right)^2}$.
As $\Tr_s$ vanishes outside of
$\bigoplus_{p,q \ge 0} \Omega^{p,q,0}(M, {\frak B}, \End({\cal E}))$, 
the proposition follows.
\end{pf}

Thus
$[\ch(A^\prime)] \in H^{\prime, even}_{\frak B}(M)$.

Let ${\cal E}_1$ and ${\cal E}_2$ be smooth ${\frak B}$-vector bundles on $M$
with superflat structures.  Suppose that there is a smooth isomorphism
$\alpha : {\cal E}_1 \rightarrow {\cal E}_2$ of 
${\cal E}_1$ and ${\cal E}_2$ as topological ${\frak B}$-vector bundles.
Choose partially flat superconnections 
$A^\prime_1$, $A^\prime_2$ on ${\cal E}_1$ and ${\cal E}_2$, 
respectively,
which are compatible with the superflat structures.
For $u \in [0,1]$, put 
$A(u) = u A^\prime_1 + (1-u) \, \alpha^* A^\prime_2$.
Note that for $u \in (0,1)$, $A(u)$ may not be partially flat on
${\cal E}_1$.
\begin{definition}
The relative Chern-Simons class $CS\left(A^\prime_1, A^\prime_2
\right) \in \overline{\Omega}^{\prime \prime,odd}(M, {\frak B})$ is
\begin{equation} \label{superCSdef}
CS\left(A^\prime_1,
A^\prime_2 \right) = 
- \int_0^1 \Tr_s \left( \left(
\partial_u A(u) \right)
e^{- A^2(u)} \right) du.
\end{equation}
\end{definition}

By construction,
\begin{equation} \label{dsuperCS}
d CS\left(A^\prime_1,
A^\prime_2 \right) = \ch(A^\prime_1) -  
\ch(A^\prime_2)
\end{equation}
vanishes in $\overline{\Omega}^{\prime \prime,even}(M, {\frak B})$. 
Thus there is a class
$\left[ CS\left(A^{\prime}_1,
A^{\prime}_2 \right) \right] \in H^{\prime \prime, odd}_{\frak B}(M)$.
\begin{proposition}
$\left[ CS\left(A^{\prime}_1,
A^{\prime}_2 \right) \right]$ actually lies in
$\bigoplus 
\begin{Sb}
p > q \\
p + q \: odd
\end{Sb}
\HH^p(M; \overline{H}_q({\frak B}))$.
\end{proposition}
\begin{pf} The proof is like that of Proposition \ref{simpler}.
We omit the details.
\end{pf}
\begin{proposition} \label{CSvarprop}
The class $\left[ CS\left(A^{\prime}_1,
A^{\prime}_2 \right) \right]$ is independent of the choice of
partially flat connections $A^\prime_1$, $A^\prime_2$; hence we denote
it by $\left[ CS\left(A^{\prime,flat}_1,
A^{\prime,flat}_2 \right) \right]$.  It only 
depends on $\alpha$ through its isotopy class. More precisely,
let $\{ \alpha(\epsilon)\}_{\epsilon \in \R}$ be a smooth $1$-parameter
family of $\alpha$'s.  Then the variation of
$CS\left(A^\prime_1, A^\prime_2
\right) \in \overline{\Omega}^{\prime \prime,odd}(M, {\frak B})$ is given by
\begin{align} \label{CSvar}
\partial_\epsilon CS\left(A^\prime_1, A^\prime_2
\right) & =
d \left( \int_0^1 \Tr_s \left( \alpha^{-1}(\epsilon) \:
\frac{d\alpha(\epsilon)}{d \epsilon} \:
e^{- A^2(u)}  \right) du
+ \int_0^1 \int_0^1  u(1-u) \right. \\
& \left. \hspace{.5in} \Tr_s \left(
\alpha^{-1}(\epsilon) 
\frac{d\alpha(\epsilon)}{d \epsilon} 
\left[ A^\prime_1 - \alpha^* A^\prime_2 ,
e^{-r A^2(u)} 
\left( A^\prime_1 - \alpha^* A^\prime_2 \right) 
e^{-(1-r) A^2(u)} \right] \right) dr du
\right) \notag
\end{align}
\end{proposition}
\begin{pf}
We first prove (\ref{CSvar}). Put $\widetilde{M} = \R \times M$. Let
$p : \R \times M \rightarrow M$ be the projection onto the second factor.  
For $i \in \{
1,2 \}$, put $\widetilde{\cal E}_i = p^* {\cal E}_i$,
$\widetilde{A}^\prime_i = p^* A^\prime_i$. Then
$\widetilde{A}^\prime_i$ is partially flat on $\widetilde{\cal E}_i$. 
Define $\widetilde{\alpha} :
\widetilde{\cal E}_1 \rightarrow \widetilde{\cal E}_2$ by saying that
$\widetilde{\alpha}\big|_{\{\epsilon\} \times M} = \alpha(\epsilon)$.
Let $\widetilde{d} = d\epsilon \: 
\partial_\epsilon + d$ denote the differential
on $\overline{\Omega}^{\prime \prime, *}(\widetilde{M}, {\frak B})$. Write
$CS\left(\widetilde{A}^\prime_1, \widetilde{A}^\prime_2 \right) \in
\overline{\Omega}^{\prime \prime, odd}(\widetilde{M}, {\frak B})$ as
\begin{equation}
CS\left(\widetilde{A}^\prime_1, \widetilde{A}^\prime_2 \right) =
CS\left(A^\prime_1, A^\prime_2 \right) (\epsilon) + d\epsilon \wedge
T(\epsilon)
\end{equation}
with $T(\epsilon) \in \overline{\Omega}^{\prime \prime, even}(M, {\frak B})$.
Then the equation $\widetilde{d} \: CS\left(\widetilde{A}^\prime_1, 
\widetilde{A}^\prime_2 \right) = 0$ implies that
\begin{equation} \label{vareqn}
\partial_\epsilon CS\left(A^\prime_1, A^\prime_2 \right) = d T(\epsilon).
\end{equation}
It remains to work out $T(\epsilon)$ explicitly.

With an obvious notation, we have
\begin{align}
\widetilde{A}^\prime_1 & = d\epsilon \wedge \partial_\epsilon + A^\prime_1, \\
\widetilde{\alpha}^* \widetilde{A}^\prime_2 & = 
d\epsilon \wedge \left[ \partial_\epsilon + \alpha^{-1}(\epsilon) 
\frac{d\alpha(\epsilon)}{d \epsilon} \right]
 + \alpha(\epsilon)^* A^\prime_2. \notag
\end{align}
Then
\begin{align}
\widetilde{A}(u) & = d\epsilon \wedge
\left[ \partial_\epsilon + (1-u) \: \alpha^{-1}(\epsilon) \: 
\frac{d\alpha(\epsilon)}{d \epsilon} \right] + u A^\prime_1 + (1-u) \:
\alpha(\epsilon)^* A^\prime_2 \\
& = d\epsilon \wedge \left[ \partial_\epsilon + 
(1-u) \: \alpha^{-1}(\epsilon) \: 
\frac{d\alpha(\epsilon)}{d \epsilon} \right] + A(u),
\notag \\
\partial_u \widetilde{A}(u) & = 
- \: d\epsilon \wedge \alpha^{-1}(\epsilon) \:
\frac{d\alpha(\epsilon)}{d \epsilon} + A^\prime_1 - \alpha(\epsilon)^* 
A^\prime_2. \notag
\end{align}
As $\alpha(\epsilon)^* A^\prime_2 = \alpha^{-1}(\epsilon) \circ A^\prime_2
\circ \alpha(\epsilon)$, it follows that
\begin{equation}
\partial_\epsilon \left[\alpha(\epsilon)^* A^\prime_2 \right] =
\left[ \alpha(\epsilon)^* A^\prime_2, \alpha^{-1}(\epsilon) 
\frac{d\alpha(\epsilon)}{d \epsilon} \right]. 
\end{equation}
Then one finds
\begin{equation}
\widetilde{A}^2(u)
= u(1-u) \: d\epsilon \wedge \left[
\alpha^{-1}(\epsilon) 
\frac{d\alpha(\epsilon)}{d \epsilon},
A^\prime_1 - \alpha(\epsilon)^* 
A^\prime_2 \right]
+ A^2(u).
\end{equation}
Thus $d\epsilon \wedge T(\epsilon)$ is the $d\epsilon$-term of
\begin{align}
& - \int_0^1 \Tr_s \left( \left(
\partial_u \widetilde{A}(u) \right)
e^{-  \widetilde{A}^2(u)} \right) du = \\
& - \int_0^1 \Tr_s \left( \left( - \: 
d\epsilon \wedge \alpha^{-1}(\epsilon) \:
\frac{d\alpha(\epsilon)}{d \epsilon}
+ A^\prime_1 - \alpha^* A^\prime_2 \right) \right. \notag \\
& \left. \hspace{1in} e^{- \left( u(1-u) \: d\epsilon \wedge 
\left[\alpha^{-1}(\epsilon) 
\frac{d\alpha(\epsilon)}{d \epsilon},
A^\prime_1 - \alpha(\epsilon)^* 
A^\prime_2 \right] +  A^2(u) \right)} \right) du, \notag
\end{align}
giving
\begin{align} \label{finally}
d\epsilon \wedge
T(\epsilon)  = \: & d\epsilon \wedge \int_0^1 \Tr_s \left( \alpha^{-1}
(\epsilon) \:
\frac{d\alpha(\epsilon)}{d \epsilon} \:
e^{- A^2(u)}  \right) du \\
& + \int_0^1 \int_0^1 \Tr_s \left( \left(
A^\prime_1 - \alpha^* A^\prime_2 \right) 
e^{-(1-r) A^2(u)} \right. \notag \\
& \left. \hspace{.5in} u(1-u) \: d\epsilon \wedge 
\left[ \alpha^{-1}(\epsilon) 
\frac{d\alpha(\epsilon)}{d \epsilon},
A^\prime_1 - \alpha^* A^\prime_2 \right]
e^{-r A^2(u)} \right) dr du \notag \\
= \:& d\epsilon \wedge \int_0^1 \Tr_s \left( \alpha^{-1}
(\epsilon) \:
\frac{d\alpha(\epsilon)}{d \epsilon} \:
e^{- A^2(u)}  \right) du \notag \\
& + d\epsilon \wedge \int_0^1 \int_0^1  u(1-u) \Tr_s \left(
\left[ \alpha^{-1}(\epsilon) 
\frac{d\alpha(\epsilon)}{d \epsilon},
A^\prime_1 - \alpha^* A^\prime_2 \right]
e^{-r A^2(u)}  \right. \notag \\
& \left. \hspace{1in} \left( A^\prime_1 - \alpha^* A^\prime_2 \right) 
e^{-(1-r) A^2(u)} \right) dr du \notag \\
= \:& d\epsilon \wedge \int_0^1 \Tr_s \left( \alpha^{-1}
(\epsilon) \:
\frac{d\alpha(\epsilon)}{d \epsilon} \:
e^{- A^2(u)}  \right) du \notag \\
& + d\epsilon \wedge \int_0^1 \int_0^1  u(1-u) \Tr_s \left(
\alpha^{-1}(\epsilon) 
\frac{d\alpha(\epsilon)}{d \epsilon}  \right. \notag \\
& \left. \hspace{1in} \left[ A^\prime_1 - \alpha^* A^\prime_2 ,
e^{-r A^2(u)} 
\left( A^\prime_1 - \alpha^* A^\prime_2 \right) 
e^{-(1-r) A^2(u)} \right] \right) dr du \notag
\end{align}
Equation (\ref{CSvar}) follows from combining (\ref{vareqn}) and 
(\ref{finally}).

Thus $\left[ CS\left(A^{\prime}_1,
A^{\prime}_2 \right) \right]$ only depends on $\alpha$ through its
isotopy class. 
A similar argument, working on $\R \times M$, shows that
$\left[ CS\left(A^{\prime}_1,
A^{\prime}_2 \right) \right]$ is independent of the choice of
partially flat connections $A^\prime_1$, $A^\prime_2$. 
\end{pf}

Put 
\begin{equation}
v = A^\prime_{0,0,1} \in C^\infty \left(M; \Hom_{\frak B}
\left( {\cal E}^\bullet, {\cal E}^{\bullet + 1} \right) \right).
\end{equation}
Then the partial flatness condition implies that
\begin{align} \label{pfimplications}
v^2 & = 0, \\
\nabla^{{\cal E},1,0} \: v & = 0, \notag \\
\left( \nabla^{{\cal E},1,0} \right)^2 + \left[ v, A^\prime_{2,0,-1} \right] 
& = 0. \notag
\end{align}
Thus there is a cochain complex of ${\frak B}$-vector bundles
\begin{equation}
({\cal E}, v) : 0 \rightarrow {\cal E}^0 \stackrel{v}{\rightarrow}
{\cal E}^1 \stackrel{v}{\rightarrow}
\ldots \stackrel{v}{\rightarrow}
{\cal E}^n \rightarrow 0.
\end{equation}
\begin{definition}
For $m \in M$, let $H({\cal E},v)_m = \bigoplus_{i=0}^n \HH^i({\cal E},v)_m$ 
be the cohomology of the complex $({\cal E},v)_m$ over $m$.
\end{definition}

We cannot conclude immediately that $H({\cal E},v)_m$ is a projective module. 
Put $\overline{\cal E} = \Lambda \otimes_{\frak B} {\cal E}$ and
$\overline{v}_m = \Id_\Lambda \otimes_{\frak B} v_m$.
\begin{hypothesis}  \label{hypo2}
For all $m \in M$, the map $\overline{v}_m$ has closed image.
\end{hypothesis}
\noindent
{\bf Remark 4 : } If ${\frak B} =\Lambda = \C$ then, 
as the fiber of ${\cal E}$ is
finitely-generated, Hypothesis \ref{hypo2} is automatically satisfied.
\begin{proposition} \label{projmodule}
Under Hypothesis \ref{hypo2}, $H({\cal E},v)_m$ is a finitely-generated 
projective ${\frak B}$-module.
\end{proposition}
\begin{pf}
The claim is that if we have a cochain complex
\begin{equation}
({\frak E}, v) : 0 \rightarrow {\frak E}^0 \stackrel{v}{\rightarrow}
{\frak E}^1 \stackrel{v}{\rightarrow}
\ldots \stackrel{v}{\rightarrow}
{\frak E}^n \rightarrow 0.
\end{equation}
of finitely-generated projective ${\frak B}$-modules and if 
$\Image(\overline{v})$ is
closed then $\HH^*({\frak E}, v)$ is a finitely-generated projective
${\frak B}$-module. We first prove a small lemma.

\begin{lemma} \label{smalllemma}
Suppose that $E$ is a finitely-generated projective left ${\frak B}$-module.
Put $\overline{E} = \Lambda \otimes_{\frak B} E$. Given $T \in \End_{\frak B}
(E)$, put $\overline{T} = \Id_\Lambda \otimes_{\frak B} T \in
\End_{\Lambda} (\overline{E})$. If $\overline{T}$ is invertible then 
$T$ is invertible.
\end{lemma}
\begin{pf}
Write $E = {\frak B}^N e$ for some projection $e \in M_N({\frak B})$. Then
$\overline{E} = \Lambda^N \overline{e}$, with $\overline{e} \in M_N(\Lambda)$.
We can consider $T$ to be an element $T^\prime \in M_N({\frak B})$
satisfying $e T^\prime = T^\prime e = T^\prime$. Put
$S = T^\prime + 1 - e \in M_N({\frak B})$ and $\overline{S} = 
\overline{T^\prime} + 1 - \overline{e} \in M_N(\Lambda)$. Then
$\overline{S}$ is invertible. Hence $S$ is invertible 
\cite[Proposition A.2.2]{Bost (1990)}. The inverse of $T$ is given
by the restriction of $S^{-1}$ to $\Image(e)$.
\end{pf}

Now put ${\frak B}$-Hermitian metrics on $\{{\frak E}^i\}_{i=1}^n$. Let
$v^* \in \Hom_{\frak B}\left({\frak E}^\bullet, 
{\frak E}^{\bullet - 1}\right) $ be the adjoint to $v$, 
defined using these metrics. Put $\triangle = v v^* + v^* v$. 
Put $\overline{\frak E} = \Lambda \otimes_{\frak B}  {\frak E}$,
$\overline{v} =  \Id_\Lambda \otimes_{\frak B} v$ and $\overline{\triangle} =
\overline{v} \, \overline{v}^* + \overline{v}^* \, \overline{v}$. Then
$(\overline{\frak E}, \overline{v})$ is a cochain complex of finitely-generated
projective Hilbert $\Lambda$-modules with
$\Image(\overline{v})$ closed.
We use \cite[Theorem 15.3.8]{Wegge-Olsen (1993)}, about operators with
closed image, throughout. It implies that $\Ker(\overline{v})$ is
a finitely-generated projective Hilbert $\Lambda$-module with 
$\Image(\overline{v})$ as a Hilbert $\Lambda$-submodule. As usual,
\begin{equation}
\Ker(\overline{v}) \cap \Image(\overline{v})^\perp =
\Ker(\overline{v}) \cap \Ker(\overline{v}^*) =
\Ker(\overline{\triangle}).
\end{equation}
As $\overline{v}^*$ is conjugate to $\overline{v}$, 
$\Image({\overline{v}^*})$ is also closed. 
There is an inclusion map $r : \Image(\overline{\triangle}) \rightarrow
\Image(\overline{v}) \oplus 
\Image(\overline{v}^*)$. We claim
that $r$ is onto. We have that $\overline{v}$ is an isomorphism between
$\Ker(\overline{v})^\perp = \Image(\overline{v}^*)$ and $\Image(\overline{v})$.
Thus if $z \in \Image(\overline{v})$ then there is a $y$ such that 
$z = \overline{v} \, \overline{v}^* (y)$. Similarly, there is an $x$ such that
$\overline{v}^* (y) = \overline{v}^* \, \overline{v} (x)$, giving that
$z = \overline{\triangle}(\overline{v}(x))$. The same argument applies if
$z \in \Image(\overline{v}^*)$. Thus $r$ is an isomorphism.

In particular, 
$\Image(\overline{\triangle}) \cong
\Image(\overline{v}) \oplus \Image(\overline{v}^*)$ is closed, implying 
that $\overline{\triangle}$ restricts to an isomorphism between 
$\Ker(\overline{\triangle})^\perp = \Image(\overline{\triangle})$ and 
$\Image(\overline{\triangle})$. It follows that $0$ is isolated in the
spectrum $\sigma(\overline{\triangle})$ of $\overline{\triangle}$. By
Lemma \ref{smalllemma}, $\sigma(\triangle) = \sigma(\overline{\triangle})$.
Hence we can take a small loop $\gamma$ around $0$ and form the projection
operator
\begin{equation} \label{Pkernel}
P^{Ker(\triangle)} = \frac{1}{2 \pi i} \int_\gamma 
\frac{d\lambda}{\lambda - \triangle}.
\end{equation} 
It follows that $\Ker(\triangle)$ is a finitely-generated projective
${\frak B}$-module.

If $\gamma^\prime$ is a contour around $\sigma(\triangle) - \{0\}$ then
the Green's operator of $\triangle$ is given by
\begin{equation} 
G = \frac{1}{2 \pi i} \int_{\gamma^\prime} \frac{1}{\lambda} \: 
\frac{d\lambda}{\lambda - \triangle}.
\end{equation} 

For $x \in {\frak E}$, let $\overline{x} = 1 \otimes_{\frak B} x$ denote
its image in $\overline{\frak E}$. If $x \in \Ker(v) \cap \Ker(v^*)$ then
$x \in \Ker(\triangle)$. Conversely, if $x \in \Ker(\triangle)$ then
$\overline{x} \in \Ker(\overline{\triangle})$, implying that
$\overline{x} \in \Ker(\overline{v}) \cap \Ker(\overline{v})$. 
Hence $x \in \Ker(v) \cap \Ker(v^*)$, showing that 
$\Ker(v) \cap \Ker(v^*) = \Ker(\triangle)$.
 
Finally, consider the map $s : \Ker(v) \cap \Ker(v^*) \rightarrow 
\HH({\frak E},v)$. We claim that $s$ is an isomorphism.  If
$x \in \Ker(s)$ then $x = v(y)$ for some $y$.
Then $v^* v (y) = 0$, so $\overline{v}^* \overline{v}(\overline{y}) = 0$,
so $\overline{v} (\overline{y}) = 0$, so $x = 0$. Thus $s$ is injective.
If $h \in \HH({\frak E},v)$, find some $y \in \Ker(v)$ in the equivalence
class of $h$. Clearly $y - vGv^*(y) \in \Ker(v)$. By usual arguments,
$\overline{y - vGv^*(y)} \in \Ker(\overline{v}^*)$ and hence
$y - vGv^*(y) \in \Ker(v^*)$. As $s(y - vGv^*(y)) = h$,
$s$ is onto.

We have shown that $\HH({\frak E},v) \cong \Ker(v) \cap \Ker(v^*) =
\Ker(\triangle)$ is a finitely-generated projective ${\frak B}$-module. 
\end{pf}

Hereafter, we assume that Hypothesis \ref{hypo2} is satisfied.
\begin{proposition} \label{fittogether}
The $\{H({\cal E},v)_m\}_{m \in M}$ fit together to form a $\Z$-graded 
${\frak B}$-vector bundle $H({\cal E},v)$ on $M$ with a flat structure.
\end{proposition}
\begin{pf} By (\ref{pfimplications}), 
$v$ is covariantly-constant with respect to 
the connection $\nabla^{{\cal E},1,0}$. Given $m \in M$, we can use the
parallel transport of $\nabla^{{\cal E},1,0}$ to extend the result of
Proposition \ref{projmodule} uniformly to a neighborhood of $m$, giving
the ${\frak B}$-vector bundle structure on $H({\cal E},v)$.
The flat structure on $H({\cal E},v)$ comes from 
\cite[Prop. 2.5]{Bismut-Lott (1995)}.
\end{pf}

There is a Hermitian metric $h^{H({\cal E},v)}$ on 
$H({\cal E},v)$ coming from its
identification with $\Ker(\triangle) \subset {\cal E}$. Letting
$P^{Ker(\triangle)}$ be as in the proof of Proposition 
\ref{projmodule}, there
is an induced connection 
\begin{equation} \label{inducedconn}
\nabla^{H({\cal E},v)} = P^{Ker(\triangle)} \nabla^{\cal E}
\end{equation}
on $H({\cal E},v)$.
\begin{proposition} \label{dualconn}
The connection $\nabla^{H({\cal E},v)}$ is partially flat and compatible
with the flat structure on $H({\cal E},v)$. Furthermore,
$\left( \nabla^{H({\cal E},v)} \right)^* = 
P^{Ker(\triangle)} \left( \nabla^{\cal E} \right)^*$.
\end{proposition}
\begin{pf}
The proof is similar to that of \cite[Prop. 2.6]{Bismut-Lott (1995)}. We
omit the details.
\end{pf}
\subsection{A Finite-Dimensional Index Theorem} 
\label{Hermitian Metrics and Superconnections}
Let $\langle \cdot, \cdot \rangle$ be a Hermitian metric on
${\cal E}$ as in Subsection \ref{Hermitian Metrics}, which respects the
$\Z$-grading on ${\cal E}$. As in
that subsection, there is a partially flat degree-$1$ 
superconnection $\overline{A}^{\prime *}$ on $\overline{\cal E}^*$ and an
adjoint partially flat degree-$1$
superconnection $A^{\prime \prime} = \left( A^{\prime} \right)^*$ on 
${\cal E}$ given by
\begin{equation}
A^{\prime \prime} = \left( h^{{\cal E}} 
\right)^{-1} \circ
\overline{A}^{\prime *} \circ h^{{\cal E}}.
\end{equation}
Explicitly,
define an adjoint operation on $\Omega \left(
M; \Hom_{\frak B} ({\cal E}, \Omega_*({\frak B}) \otimes_{\frak B} 
{\cal E})  \right)$ by requiring that\\
1. For $\alpha, \alpha^\prime \in \Omega \left(
M; \Hom_{\frak B} ({\cal E}, \Omega_*({\frak B}) \otimes_{\frak B} 
{\cal E} ) \right)$, $( \alpha \alpha^\prime )^* = \alpha^{\prime *} 
\alpha^*$.\\
2. If $V \in C^\infty\left(
M; \Hom_{\frak B} ({\cal E}, \Omega_*({\frak B}) \otimes_{\frak B} 
{\cal E} ) \right)$ then $V^*$ is the adjoint defined using the Hermitian form
(\ref{formform}).\\
3. If $\omega \in \Omega^1(M)$ then its extension by the identity
to become an element of $\Omega^1 \left(
M; \End_{\frak B} ({\cal E}) \right)$ satisfies
$\omega^* = - \overline{\omega}$.

Then for $A^\prime$ as in (\ref{supercomp}),
\begin{equation}
A^{\prime \prime} = \sum_p A^{\prime \prime}_p,
\end{equation}
where $A^{\prime \prime}_1 = \left( \nabla^{\cal E} \right)^*$ and
for $p \neq 1$, $A^{\prime \prime}_p = \left( A^{\prime}_p \right)^*$. 
We write
\begin{equation}
CS \left(A^\prime, h^{\cal E} \right) = 
CS \left(A^\prime, \overline{A}^{\prime *} \right)
\in \overline{\Omega}^{\prime \prime,odd}
(M, {\frak B}).
\end{equation}

Let $N$ be the number operator on ${\cal E}$, meaning that $N$ acts
on $C^\infty \left(M; {\cal E}^j \right)$ as multiplication by $j$.
For $t > 0$, let $\langle \cdot, \cdot \rangle_t$ be the Hermitian metric
on ${\cal E}$ such that if $e_1, e_2 \in C^\infty \left(M; {\cal E}^j \right)$
then 
\begin{equation}
\langle e_1, e_2 \rangle_t = t^j \langle e_1, e_2 \rangle \in 
C^\infty(M; {\frak B}).  
\end{equation}
Letting $h^{\cal E}_t : {\cal E} \rightarrow \overline{\cal E}^*$ be the
isomorphism induced from $\langle \cdot, \cdot \rangle_t$, we have
$h^{\cal E}_t = h^{\cal E} t^N$. Letting
$A^{\prime \prime}_t$ denote the adjoint of $A^\prime$ with 
respect to $\langle \cdot, \cdot \rangle_t$, we have
$A^{\prime \prime}_t = t^{-N} \: A^{\prime \prime} \: t^N$.
\begin{proposition}
For $u \in [0,1]$, put $A(u) = u A^\prime + (1-u) A^{\prime \prime}_t$. Then
\begin{align} \label{pt}
\partial_t CS\left(A^\prime, h^{\cal E}_t
\right) = \: & \frac1t \: 
d \left( \int_0^1 \Tr_s \left( N
e^{- A^2(u)}  \right) du
+ \int_0^1 \int_0^1  u(1-u)  \right. \\
& \left. \hspace{.5in} \Tr_s \left(
N \left[ A^\prime - A^{\prime \prime}_t ,
e^{-r A^2(u)} 
\left( A^\prime - A^{\prime \prime}_t \right) 
e^{-(1-r) A^2(u)} \right] \right) dr du
\right). \notag
\end{align}
\end{proposition}
\begin{pf}
This follows from Proposition \ref{CSvarprop}.
\end{pf}

To make the equations more symmetric, put
$B^\prime_t = t^{N/2} \: A^\prime \: t^{-N/2}$ and
$B^{\prime \prime}_t = t^{-N/2} \: A^{\prime \prime} \: t^{N/2}$. Then
$B^{\prime \prime}_t$ is the adjoint of $B^\prime_t$ with respect to
$\langle \cdot, \cdot \rangle$. Explicitly,
\begin{align}
B^\prime_t = & \sum_{p \ge 0} t^{(1-p)/2} A^\prime_p, \\
B^{\prime \prime}_t = & \sum_{p \ge 0} t^{(1-p)/2} A^{\prime \prime}_p. \notag
\end{align}
\begin{proposition} \label{tt}
For $u \in [0,1]$, put $B_t(u) = u B^\prime_t + (1-u) B^{\prime \prime}_t$.
Define ${\cal T}(t) \in 
\overline{\Omega}^{\prime \prime, even}(M, {\frak B})$
by
\begin{align} \label{torsiondef}
{\cal T}(t) = 
\: &- \: \frac1t \: 
\left( \int_0^1 \Tr_s \left( N
e^{- B_t^2(u)}  \right) du
+ \int_0^1 \int_0^1  u(1-u)  \right. \\
& \left. \hspace{.5in} \Tr_s \left(
N \left[ B_t^\prime - B^{\prime \prime}_t ,
e^{-r B_t^2(u)} 
\left( B_t^\prime - B^{\prime \prime}_t \right) 
e^{-(1-r) B_t^2(u)} \right] \right) dr du
\right). \notag
\end{align}
Then
\begin{equation} \label{CSeqn}
CS\left(A^\prime,
h^{\cal E}_t \right) = CS\left(B^\prime_t,
h^{\cal E} \right) = 
- \int_0^1 \Tr_s \left( \left(
B^\prime_t - B^{\prime \prime}_t \right)
e^{- B_t^2(u)} \right) du
\end{equation}
and
\begin{equation} \label{partialtcs}
\partial_t CS\left(B^\prime_t, h^{\cal E} \right) =
- \: d {\cal T}(t).
\end{equation}
\end{proposition}
\begin{pf}
This follows from (\ref{superCSdef}) and (\ref{pt}) by conjugating within the
supertrace by $t^{N/2}$. 
\end{pf}

We now discuss the large-$t$ asymptotics of $\ch \left(B_t(u)\right)$,
$CS\left(B^\prime_t,
h^{\cal E}\right)$ and ${\cal T}(t)$. We must first specify the
notion of convergence.
Define $i_j$ and $\parallel \cdot \parallel_j$ as after
equation (\ref{projlimit}).

Let ${\frak E}$ be a finitely-generated projective left ${\frak B}$-module
with a ${\frak B}$-Hermitian metric.
Write ${\frak E} = {\frak B}^N e$ for some fixed projection
$e \in M_N({\frak B})$. Put $e_j = i_j(e) \in M_N(B_j)$ and
$E_j = B_j^N e_j$. Then $E_j$ inherits a Banach space structure as a
closed subspace of $B_j^N$. Furthermore,
$\End_{B_j} (E_j)$ inherits a Banach algebra structure as a closed subalgebra
of $\End(E_j)$.  Note that $\End_{B_0}(E_0)$ is the same underlying
algebra as the $C^*$-algebra $\End_\Lambda(\overline{\frak E})$, 
but may have a different norm.

We can identify $\End_{\frak B}({\frak E})$ with the projective limit of
Banach algebras
\begin{equation}
\ldots \longrightarrow 
\End_{B_{j+1}}(E_{j+1}) \longrightarrow \End_{B_{j}}(E_{j})
\longrightarrow \ldots \longrightarrow \End_{B_{0}}(E_{0}). 
\end{equation}
We again write $\parallel \cdot \parallel_j$ for the induced 
submultiplicative seminorm
on $\End_{\frak B}({\frak E})$. Given $T \in \End_{\frak B}({\frak E})$,
let $T_j$ be its image in $\End_{B_{j}}(E_{j})$. Then
\begin{equation}
\ldots \supseteq \sigma(T_{j+1}) \supseteq \sigma(T_j) \supseteq \ldots
\supseteq \sigma(T_0)
\end{equation} 
and $\sigma(T) = \bigcup_{j=0}^\infty \sigma(T_j)$. As $\End_{B_{0}}(E_{0})$
is the same underlying algebra as the $C^*$-algebra 
$\End_\Lambda(\overline{\frak E})$, $\sigma(T_0) = \sigma(\overline{T})$. By
Lemma \ref{smalllemma}, $\sigma(T) = \sigma(\overline{T})$. Thus
$\sigma(T) = \sigma(T_j) = \sigma(\overline{T})$ for all $j$.

Using the description of $\Omega_*({\frak B})$ in 
\cite[Section II]{Lott (1992)}, there is a sequence of seminorms 
$\{ \parallel \cdot \parallel_j\}_{j=0}^\infty$ on each $\Omega_k({\frak B})$
coming from the norms on $B_j$. We obtain seminorms 
$\parallel \cdot \parallel_j$ on $\Hom_{\frak B} \left( {\frak E}, 
\Omega_k({\frak B}) \otimes_{\frak B} {\frak E} \right)$ and 
$\overline{\Omega}_k({\frak B})$, with respect to which (\ref{tracemap})
is continuous.  Convergence of $\ch \left(B_t(u)\right)$
or $CS\left(B^\prime_t,
h^{\cal E}\right)$ will mean convergence in all seminorms
$\{\parallel \cdot \parallel_j\}_{j=0}^\infty$.

\begin{proposition} \label{csconv}
For all $u \in (0,1)$,
as $t \rightarrow \infty$,
\begin{equation} \label{chernasymp}
\ch \left(B_t(u) \right) = 
\ch \left( \nabla^{H({\cal E},v)}(u)
\right) + O(t^{-1/2})
\end{equation}
uniformly on $M$. Also,
\begin{equation} \label{csasymp}
CS\left(B^\prime_t, h^{\cal E}\right) = 
CS\left( \nabla^{H({\cal E},v)}, h^{H({\cal E},v)} 
\right) + O(t^{-1/2})
\end{equation}
uniformly on $M$.
\end{proposition}
\begin{pf}
We will only prove (\ref{csasymp}), as the proof of 
(\ref{chernasymp}) is similar but easier.
We begin with some generalities.  Suppose that $E$ is a finitely-generated
projective left ${\frak B}$-module, $T_1$ is an element of
$\End_{\frak B} (E)$ and $T_2$ is an element of
$\End_{\frak B} (E) \otimes {\frak G}$ for some Grassmann algebra
${\frak G}$. We assume that $T_2$ has positive Grassmann degree.
Suppose that $\sigma(T_1) \subset \R$
and that $0$ is isolated in $\sigma(T_1)$. Let $\gamma_1$ be a small
loop around $0$ and let $\gamma_2$ be a contour around $\sigma(T_1) - \{ 0 \}$.
We orient $\gamma_1$ and $\gamma_2$ counterclockwise. Then we
can write
\begin{equation} \label{starteqn}
P^{Ker(T_1)} = \int_{\gamma_1} \frac{1}{z-T_1} \: \frac{dz}{2 \pi i}
\end{equation}
and
\begin{equation}
P^{Im(T_1)} = \int_{\gamma_2} \frac{1}{z-T_1} \: \frac{dz}{2 \pi i}.
\end{equation}
As $\sigma(T_1 + T_2) = \sigma(T_1)$, we can write
\begin{equation}
e^{- t (T_1 +  T_2)^2} = 
\int_{\gamma_1} \frac{e^{- t z^2}}{z - (T_1 +  T_2)} \: 
\frac{dz}{2 \pi i} +
\int_{\gamma_2} \frac{e^{- t z^2}}{z - (T_1 +  T_2)} \: 
\frac{dz}{2 \pi i}.
\end{equation}
Using the series
\begin{equation}
\frac{1}{z - (T_1 +  T_2)} = \frac{1}{z - T_1} + 
\frac{1}{z - T_1} \:  T_2 \: \frac{1}{z - T_1} + \ldots,
\end{equation}
the first contour integral becomes
\begin{align} \label{cont1}
\int_{\gamma_1} \frac{e^{- t z^2}}{z - (T_1 +  T_2)} 
\: \frac{dz}{2 \pi i} = \: & 
P^{Ker(T_1)} -  P^{Ker(T_1)} \:  T_2 \: G - 
G \:  T_2 \: P^{Ker(T_1)} \\
& + P^{Ker(T_1)} \:  T_2 \: G \:  T_2 \: G +
G \:  T_2 \: P^{Ker(T_1)} \:  T_2 \: G \notag \\
& + G \:  T_2 \: G \:  T_2 \: P^{Ker(T_1)} -
t \: P^{Ker(T_1)} \:  T_2 \: P^{Ker(T_1)} \: 
 T_2 \: P^{Ker(T_1)} \notag \\
& + \ldots \notag
\end{align}
where $G$ is the Green's operator of $T_1$.

Writing out $B_t^\prime$ and $B_t^{\prime \prime}$ explicitly, we have
\begin{align} \label{B_t}
B_t^\prime - B_t^{\prime \prime}  \:& = \sqrt{t} \: (v - v^*) + 
\nabla^{\cal E} - \left( \nabla^{\cal E} \right)^*  +
O( t^{-1/2}), \\
B_t(u) & = \sqrt{t} \: (u v + (1-u) v^*) + 
u \nabla^{\cal E} + (1-u) \left( \nabla^{\cal E} \right)^*  +
O( t^{-1/2}). \notag
\end{align}
We apply equations (\ref{starteqn}) - (\ref{cont1}) with
\begin{align} \label{t1t2}
T_1 \: & = u v + (1-u) v^*, \\
 T_2 \: & = t^{-1/2} B_t(u) - T_1 = 
t^{-1/2} \left( u \nabla^{\cal E} + (1-u) \left( \nabla^{\cal E} \right)^* 
\right) +
O( t^{-1}). \notag
\end{align}
For $u \in (0,1)$, $\Ker(T_1) = \Ker(u(1-u)\triangle) =
\Ker(\triangle)$. 
Let us write
\begin{equation}
CS \left( B^\prime_t, h^{\cal E} \right) = CS_1 + CS_2,
\end{equation}
with
\begin{align}
CS_1 \: & = - \int_0^1 \Tr_s \left( \left(
B^\prime_t - B^{\prime \prime}_t \right)
\int_{\gamma_1} \frac{e^{- t z^2}}{z - (T_1 +  T_2)} \: 
\frac{dz}{2 \pi i} \right) du, \\
CS_2 \: & = - \int_0^1 \Tr_s \left( \left(
B^\prime_t - B^{\prime \prime}_t \right)
\int_{\gamma_2} \frac{e^{- t z^2}}{z - (T_1 +  T_2)} \: 
\frac{dz}{2 \pi i} \right) du. \notag
\end{align}
Substituting the series (\ref{cont1}) into $CS_1$, we see that 
the leading
terms in $t$ come from the terms in (\ref{cont1}) without any 
factors of $G$. Using (\ref{inducedconn}) and Proposition \ref{dualconn},
one finds
\begin{align} \label{CS1}
CS_1 \: & = - \int_0^1 \Tr_s \left( \left(
\sqrt{t} \: (v - v^*) + 
\nabla^{\cal E} - \left( \nabla^{\cal E} \right)^* \right) P^{Ker(\triangle)}
\right. \\
& \left. \hspace{.5in} 
\int_{\gamma_1} \frac{e^{- t z^2}}{z - t^{-1/2} P^{Ker(\triangle)}
\left( u \nabla^{\cal E} + 
(1-u) \left( \nabla^{\cal E} \right)^* \right)
P^{Ker(\triangle)}} \: P^{Ker(\triangle)} \: 
\frac{dz}{2 \pi i} \right) du + O(t^{-1/2}) \notag \\
& = - \int_0^1 \Tr_s \left( \left(
\nabla^{H({\cal E},v)} - \left( \nabla^{H({\cal E},v)} \right)^* \right)
e^{- {\left( u \nabla^{H({\cal E},v)} + 
(1-u) \left( \nabla^{H({\cal E},v)} 
\right)^* \right)}^2} \right) du + O(t^{-1/2}) 
\notag \\
& = CS\left( \nabla^{H({\cal E},v)}, h^{H({\cal E},v)} 
\right) + O(t^{-1/2}). \notag
\end{align}
Note that only a finite number of terms of the series (\ref{cont1})
contribute to the component of $CS_1$ of a given degree in 
$\overline{\Omega}^{\prime \prime,odd}
(M, {\frak B})$. Thus the derivation of (\ref{CS1}) is purely algebraic.

It remains to estimate $CS_2$. Put $X = T_1 T_2 + T_2 T_1 + T_2^2$ and
$e^{-r T_1^{\prime 2}} =
P^{Im(T_1)} e^{-r T_1^2} P^{Im(T_1)}$.
We use the heat kernel expansion
\begin{align} \label{Duhamel}
\int_{\gamma_2} \frac{e^{- t z^2}}{z - (T_1 +  T_2)} \: 
\frac{dz}{2 \pi i} = 
& \int_0^\infty e^{-r_0 T_1^{\prime 2}} \delta(r_0 - t) \: dr_0 \\
& - \int_0^\infty \int_0^\infty  e^{-r_0 T_1^{\prime 2}}
X e^{-r_1 T_1^{\prime 2}} \delta(r_0 + r_1 - t)
\: dr_0 dr_1 \notag \\
& + \int_0^\infty P^{Ker(T_1)} X e^{-r_0 T_1^{\prime 2}}
\theta(r_0 - t) \: dr_0 \notag \\
& + \int_0^\infty e^{-r_0 T_1^{\prime 2}} X P^{Ker(T_1)}
\theta(r_0 - t) \: dr_0 + \ldots \notag
\end{align}
Here $\theta$ is the Heaviside function : $\theta(r) = 
\frac{1 + sign(r)}{2}$.
The series (\ref{Duhamel}) is similar to the Duhamel expansion of
$e^{-t(T_1^2 + X)}$, with each intermediate
factor of $e^{-rT_1^2}$ in the Duhamel expansion being replaced by
either $e^{-r T_1^{\prime 2}}$ or $P^{Ker(T_1)}$. 
A term on the right-hand-side of (\ref{Duhamel}) with 
$k$ $X$'s and $l$ $P^{Ker(T_1)}$'s, $k \ge l$, will have a factor of
\begin{equation*}
\begin{array}{lll}
(-1)^{k} \: \delta \left(\sum_{i=0}^{k} r_i - t \right) & 
\text{ if } & l = 0, \\
\frac{(-1)^{k+l}}{(l-1)!} \: \theta \left(\sum_{i=0}^{k-l} r_i - t \right)
\left(\sum_{i=0}^{k-l} r_i - t \right)^{l-1} & \text{ if } &
l \ge 1. \notag
\end{array}
\end{equation*}
In our case, from (\ref{t1t2}),
\begin{equation}
X = t^{-1/2} \left( u^2 \nabla^{\cal E} v + 
u(1-u) \left( \left(\nabla^{\cal E} \right)^* v +
\nabla^{\cal E} v^* \right) + (1-u)^2 \left(\nabla^{\cal E} \right)^* 
v^* \right) + O(t^{-1}).
\end{equation}
Put 
\begin{equation}
e^{-r \triangle^\prime} = P^{Im(\triangle)} e^{-r\triangle} 
P^{Im(\triangle)}.
\end{equation}   
Using (\ref{Duhamel}) gives
\begin{align} \label{cs2}
CS_2 =
& - \int_0^1 \Tr_s \left[
\left( \nabla^{\cal E} - \left( \nabla^{\cal E} \right)^* 
\right) e^{- t u (1-u) \triangle^\prime}
\right] du \\
& +  \int_0^1 \int_0^\infty \int_0^\infty \Tr_s \left[ (v - v^*) 
e^{- r_0 u (1-u) \triangle^\prime} \left( u^2 \nabla^{\cal E} v + 
u(1-u) \left( \left(\nabla^{\cal E} \right)^* v +
\nabla^{\cal E} v^* \right) \right. \right. \notag \\
& \left. \left. \hspace{.5in} + (1-u)^2 \left(\nabla^{\cal E} \right)^* v^* 
\right)  e^{- r_1 u (1-u) \triangle^\prime} 
\right] \delta(r_0 + r_1 - t) \: dr_0 dr_1 du + \ldots \notag
\end{align} 

We now use the crucial fact \cite[Theorem 1.22]{Davies (1980)}
that if $\{\alpha_r\}_{r > 0}$ is a 1-parameter
semigroup in a Banach algebra then the number
\begin{equation} \label{start1}
a = \lim_{r \rightarrow \infty} r^{-1} \log \parallel \alpha_r \parallel
\end{equation}
exists. Furthermore, 
for all $r > 0$, the spectral radius of $\alpha_r$ is given by
\begin{equation} \label{start2}
\SpRad(\alpha_r) = e^{ar}.
\end{equation}
Let $\lambda_0 > 0$ be the infimum of the nonzero spectrum of $\triangle$. 
Then by the spectral mapping theorem, $\SpRad(e^{-r\triangle^\prime}) =
e^{-r \lambda_0}$. Thus for any $j \ge 0$, there is a constant $C_j > 0$ such
that for all $r > 1$,
\begin{equation} \label{end1}
\parallel  e^{-r\triangle^\prime} \parallel_j \: \le \: C_j 
\: e^{-r \lambda_0/2}.
\end{equation}

Consider the $\overline{\Omega}^{\prime \prime,1}
(M, {\frak B})$-component of $CS_2$, which is given explicitly in
(\ref{cs2}).  First, we have
\begin{equation}
\parallel \int_0^1 \Tr_s \left[
\left( \nabla^{\cal E} - \left( \nabla^{\cal E} \right)^* 
\right) e^{- t u (1-u) \triangle^\prime}
\right] du \parallel_j \: \le \const \int_0^1 e^{-tu(1-u) \lambda_0/2} du
= O(t^{-1}). 
\end{equation}
Next, consider the second term in (\ref{cs2}). 
Recall that 
$CS \left( B_t^\prime, h^{\cal E} \right)$ lives in the quotient space 
$\overline{\Omega}^{\prime \prime,odd}
(M, {\frak B})$ defined in 
(\ref{formspaces}). It follows from (\ref{pfimplications}) that
$\nabla^{\cal E} v$ will not contribute to the 
$\overline{\Omega}^{\prime \prime,1}
(M, {\frak B})$-component of $CS_2$, and similarly for
$\left( \nabla^{\cal E} \right)^* v^*$. Hence
the second term in (\ref{cs2}) is bounded above in the $\parallel \cdot
\parallel_j$-seminorm by
\begin{equation}
\const \int_0^1 \int_0^\infty \int_0^\infty 
u(1-u) e^{-t u(1-u) \lambda_0/2}
\delta(r_0 + r_1 - t) \: dr_0 dr_1 du = O(t^{-1}). 
\end{equation}

Thus we have shown that for each $j$, the 
$\parallel \cdot \parallel_j$-seminorm of the
$\overline{\Omega}^{\prime \prime,1}
(M, {\frak B})$-component of $CS_2$ is $O(t^{-1})$.
One can carry out a similar analysis for all of the terms in $CS_2$.
The point
is that for large $t$, the $e^{-t u(1-u) \lambda_0/2}$ factor ensures that
in the $u$-integral, only the behavior near $u = 0$ and $u = 1$ is important.
Consider, for example, what happens when $u$ is close to $1$. 
When $u = 1$, $B_t(1)$ is partially flat and so $X$ lies in
$\Omega^\prime \left(
M; \Hom_{\frak B} ({\cal E}, \Omega_*({\frak B}) \otimes_{\frak B} 
{\cal E} )\right)$, as defined in (\ref{omegaprime}). 
Consider the value at $u = 1$ of a given term of
$CS_2$. The contribution from $B_t^\prime - B_t^{\prime \prime}$ lies
in $\bigoplus_{p + q \pm r = 1} \Omega^{p,q,r}(M, {\frak B}, \End({\cal E}))$
and the contribution from the $X$'s lies in
$\bigoplus_{p^\prime + r^\prime \le q^\prime} 
\Omega^{p^\prime,q^\prime,r^\prime}(M, {\frak B}, \End({\cal E}))$.
For the supertrace to be nonzero, we must have $r+r^\prime = 0$.
The explicit factor of $t$ appearing is
\begin{equation}
t^{\frac12} t^{- \: \frac{p+p^\prime + q+q^\prime}{2}} t^{\frac{p^\prime
+ q^\prime + r^\prime}{2}} = t^{\frac{1 - p - q - r}{2}}.
\end{equation}
Suppose first that $p + q + r = 1$. The integral near $u = 1$ of 
$e^{-t u(1-u) \lambda_0/2}$ yields an additional factor of $t^{-1}$, giving
a total estimate of $O(t^{-1})$. Suppose now that $p + q - r = 1$.
The explicit factor of 
$t$ is $t^{-r}$. Along with the integral of $e^{-t u(1-u) \lambda_0/2}$,
we obtain a total estimate of $O(t^{-1-r})$. 
For this to be more significant than
$O(t^{-1})$, we must have $r < 0$. As $r = p + q - 1$ and $p, q \ge 0$,
the only possibility is $p = q = 0$ and $r = -1$. Then $r^\prime = 1$ and
so $p^\prime + 1 \le q^\prime$. The supertrace of such a term lies in
$\overline{\Omega}^{p^\prime, q^\prime}(M, {\frak B}) \subset 
\overline{\Omega}^{\prime, odd}(M, {\frak B})$ and so vanishes
in $\overline{\Omega}^{\prime \prime,odd}
(M, {\frak B})$.

One can carry out a similar analysis near $u = 0$.
Finally, the convergence is clearly uniform on $M$.
\end{pf}
\begin{corollary} \label{corr1}
We have
\begin{equation}
\left[ \ch \left( A^\prime \right) \right] =
\left[ \ch \left( \nabla^{H({\cal E}, v)} \right) \right]
\text{ in } H^{\prime, even}_{\frak B}(M)
\end{equation}
and
\begin{equation} \label{equiv}
\left[ CS \left( A^\prime, h^{\cal E} \right) \right] =
\left[ CS \left( \nabla^{H({\cal E}, v)}, h^{H({\cal E}, v)} \right) \right]
\text{ in }
\bigoplus 
\begin{Sb}
p > q \\
p + q \: odd
\end{Sb}
\HH^p(M; \overline{H}_q({\frak B})). 
\end{equation}
\end{corollary}
\begin{pf}
For all $u \in (0,1)$,  $\left[ \ch \left( A(u) \right) \right] =
\left[ \ch \left( A^\prime \right) \right]$ and
$\left[ \ch \left( \nabla^{H({\cal E}, v)}(u) \right) \right] =
\left[ \ch \left( \nabla^{H({\cal E}, v)} \right) \right]$. 
As $A^\prime = B_1^\prime$, 
the corollary follows from 
(\ref{partialtcs}) and Proposition \ref{csconv}.
\end{pf}
{\bf Remark 5 : } If ${\frak B} = \C$ then (\ref{equiv})
is equivalent to \cite[Theorem 2.14]{Bismut-Lott (1995)}. 
\begin{proposition} \label{torsionconv}
As $t \rightarrow \infty$,
\begin{equation} \label{torsionasymp}
{\cal T}(t) = O(t^{-3/2})
\end{equation}
uniformly on $M$.
\end{proposition}
\begin{pf}
Let $\epsilon$ be a new variable which commutes
with the other variables and satisfies $\epsilon^2 = 0$. Then
\begin{align}
{\cal T}(t) = 
\: &- \: \frac1t \: 
\left( \int_0^1  \Tr_s \left( N
e^{- B_t^2(u)}  \right) du
+ \partial_\epsilon \int_0^1  u(1-u)  \right. \\
& \left. \hspace{.5in} \Tr_s \left(
N \left[ B_t^\prime - B^{\prime \prime}_t ,
e^{- B_t^2(u) +\epsilon \left( B_t^\prime - B^{\prime \prime}_t \right)} 
\right] \right) du
\right). \notag
\end{align}
The implication is that we can again write ${\cal T}(t)$  as 
${\cal T}_1(t) + {\cal T}_2(t)$, where ${\cal T}_i(t)$ is a
contour integral around $\gamma_i$.
Using the method of proof of Proposition \ref{csconv}, one finds
\begin{align} 
{\cal T}_1(t) = 
\: &- \: \frac1t \: 
\left( \int_0^1 \Tr_s \left( N
e^{- \nabla^{H({\cal E},v)^2}(u)}  \right) du
+ \int_0^1 \int_0^1  u(1-u) 
\Tr_s \left(
N \left[ \nabla^{H({\cal E},v)} -
\left( \nabla^{H({\cal E},v)} \right)^* , \right. \right. \right.\\
& \hspace{.5in} \left. \left. \left. e^{-r \nabla^{H({\cal E},v)^2}(u)} 
\left( \nabla^{H({\cal E},v)} -
\left( \nabla^{H({\cal E},v)} \right)^* \right) 
e^{-(1-r) \nabla^{H({\cal E},v)^2}(u)} \right] \right) dr du
\right) + O(t^{-3/2})  \notag \\
= \: &- \: \frac1t \: \int_0^1
\Tr_s \left( N
e^{- \nabla^{H({\cal E},v)^2}(u)}  \right) du
+ O(t^{-3/2}).  \notag
\end{align}
As $\Tr_s \left( N
e^{- \nabla^{H({\cal E},v)^2}(u)}  \right) \in \overline{\Omega}^{\prime,even}
(M, {\frak B})$,
it follows that ${\cal T}_1(t)$ is $O(t^{-3/2})$ in 
$\overline{\Omega}^{\prime \prime,even}(M, {\frak B})$.

Next, consider ${\cal T}_2(t)$. Counting powers of $t$ as in the proof of
Proposition \ref{csconv}, one finds that ${\cal T}_2(t)$ is $O(t^{-3/2})$
except for possible terms which decay like $t^{-1}$ and lie in
$\bigoplus_p \overline{\Omega}^{p, p}(M, {\frak B})  \mod
\overline{\Omega}^{\prime,even}(M, {\frak B})$. Let us
write such a term as $t^{-1} {\cal T}_2^{p,p}$, with
${\cal T}_2^{p,p} \in \overline{\Omega}^{p, p}(M, {\frak B})$.
From (\ref{partialtcs}), we see that $t^{-1} d^{1,0}
{\cal T}_2^{p,p}$ comes from the $t$-derivative of the 
$\overline{\Omega}^{p+1, p}(M, {\frak B})$-component of $CS \left(
B_t^\prime, h^{\cal E} \right)$. However, as $CS \left(
B_t^\prime, h^{\cal E} \right)$ has no $\log(t)$-term
in its asymptotics, it follows that
$d^{1,0} {\cal T}_2^{p,p} = 0$. Thus ${\cal T}_2^{p,p}$ lies in
$Z^p \left( M; \overline{\Omega}_p({\frak B}) \right)$. As we quotient
by this subspace in defining 
$\overline{\Omega}^{\prime \prime, even}(M, {\frak B})$, equation
(\ref{torsionasymp}) follows.     

Again, the convergence is clearly uniform on $M$.   
\end{pf}
{\bf Remark 6 : } There may seem to be a contradiction between
Proposition \ref{torsionconv} and
\cite[Theorem 2.13]{Bismut-Lott (1995)}, in which a nonzero $O(t^{-1})$
term for ${\cal T}(t)$ was found.  However, in the present paper
we quotient by $Z^k \left( M; \overline{\Omega}_k({\frak B}) \right)$ in
defining 
$\overline{\Omega}^{2k} \left( M; \overline{\Omega}_p({\frak B}) \right)$.
When ${\frak B} = \C$, as in \cite{Bismut-Lott (1995)}, this quotienting
removes the $O(t^{-1})$ term of \cite[Theorem 2.13]{Bismut-Lott (1995)}
and so there is no contradiction.

\subsection{The Analytic Torsion Form} \label{The Analytic Torsion Form}

In this subsection we consider the special case when ${\cal E}$ has
not only a partially flat degree-$1$ superconnection, but has a partially
flat connection in the sense of Definition
\ref{partiallyflat}. Let
\begin{equation}
({\cal E}, v) : 0 \rightarrow {\cal E}^0 \stackrel{v}{\rightarrow}
{\cal E}^1 \stackrel{v}{\rightarrow}
\ldots \stackrel{v}{\rightarrow}
{\cal E}^n \rightarrow 0
\end{equation}
be a cochain complex of ${\frak B}$-vector bundles on $M$. 
Let 
\begin{equation}
\nabla^{\cal E} = \bigoplus_{i=0}^n \nabla^{{\cal E}^i}
\end{equation}
be a partially flat connection on 
${\cal E} = \bigoplus_{i=0}^n {{\cal E}^i}$.
Suppose that $\nabla^{{\cal E},1,0} v = 0$. Put
\begin{equation}
A^\prime = v + \nabla^{\cal E}.
\end{equation}
Then $A^\prime$ is a partially flat degree-$1$ superconnection.
In the notation of Subsection \ref{Hermitian Metrics and Superconnections},
\begin{equation}
\begin{align}
B^\prime_t = & \sqrt{t} \: v + \nabla^{\cal E}, \\
B^{\prime \prime}_t = & \sqrt{t} \: v^* + \left( \nabla^{\cal E}
\right)^*. \notag
\end{align}\end{equation}
\begin{proposition} \label{smallt}
As $t \rightarrow 0$,
\begin{equation} \label{smalltch}
\ch\left(B_t(u) \right) = 
\ch \left( \nabla^{\cal E}(u)
\right) + O(t),
\end{equation}
\begin{equation} \label{smalltcs}
CS\left(B^\prime_t, h^{\cal E}\right) = 
CS\left( \nabla^{\cal E}, h^{\cal E}
\right) + O(t)
\end{equation}
and
\begin{equation} \label{smallttorsion}
{\cal T}(t) = O(1).
\end{equation}
\end{proposition}
\begin{pf}
Equations (\ref{smalltch}) and (\ref{smalltcs}) are 
evident. From (\ref{torsiondef}),
the $t^{-1}$-term of ${\cal T}$ is
\begin{align} 
&- \: \frac1t \: 
\left( \int_0^1 \Tr_s \left( N
e^{- \nabla^{{\cal E}^2}(u)}  \right) du
+ \int_0^1 \int_0^1  u(1-u) 
\Tr_s \left(
N \left[ \nabla^{{\cal E}} -
\left( \nabla^{{\cal E}} \right)^* , \right. \right. \right.\\
& \hspace{.5in} \left. \left. \left. e^{-r \nabla^{{\cal E}^2}(u)} 
\left( \nabla^{{\cal E}} -
\left( \nabla^{{\cal E}} \right)^* \right) 
e^{-(1-r) \nabla^{{\cal E}^2}(u)} \right] \right) dr du
\right) \notag \\
= \: &- \: \frac1t \: \int_0^1
\Tr_s \left( N
e^{- \nabla^{{\cal E}^2}(u)}  \right) du. \notag
\end{align}
As $\Tr_s \left( N
e^{- \nabla^{{\cal E}^2}(u)}  \right) \in \overline{\Omega}^{\prime,even}
(M, {\frak B})$,
this vanishes in $\overline{\Omega}^{\prime \prime, even}(M, {\frak B})$.
It is easy to check that there is no $O(t^{-1/2})$-term.
\end{pf}
\begin{corollary} \label{pushforward}
We have
\begin{equation}
\left[ \ch \left( \nabla^{\cal E} \right) \right] =
\left[ \ch \left( \nabla^{H({\cal E}, v)} \right) \right]
\text{ in } H^{\prime,even}_{\frak B}(M)
\end{equation}
and
\begin{equation}
\left[ CS \left( \nabla^{\cal E}, h^{\cal E} \right) \right] =
\left[ CS \left( \nabla^{H({\cal E}, v)}, h^{H({\cal E}, v)} \right) \right]
\text{ in }
\bigoplus 
\begin{Sb}
p > q \\
p + q \: odd
\end{Sb}
\HH^p(M; \overline{H}_q({\frak B})). 
\end{equation}
\end{corollary}
\begin{pf}
For all $u \in (0,1)$, $\left[ \ch \left( \nabla^{\cal E}(u) \right) \right]
= \left[ \ch \left( \nabla^{\cal E} \right) \right]$ and
$\left[ \ch \left( \nabla^{H({\cal E}, v)}(u) \right) \right] =
\left[ \ch \left( \nabla^{H({\cal E}, v)} \right) \right]$. The corollary
now follows from Corollary \ref{corr1} and Proposition \ref{smallt}.
\end{pf} 
{\bf Remark 7 : } If ${\frak B} = \C$ then Corollary \ref{pushforward}
is equivalent to \cite[Theorem 2.19]{Bismut-Lott (1995)}. 
\begin{definition}
The analytic torsion form
${\cal T} \in \overline{\Omega}^{\prime \prime, even}(M, {\frak B})$
is given by
\begin{equation} \label{torsioneqn}
{\cal T} = \int_0^\infty {\cal T}(t) \: dt.
\end{equation}
\end{definition}

By Propositions \ref{torsionconv} and \ref{smallt}, the integral in
(\ref{torsioneqn}) makes sense.
\begin{proposition}
We have
\begin{equation}
d {\cal T} = CS \left( \nabla^{\cal E}, h^{\cal E} \right) -
CS \left( \nabla^{H({\cal E}, v)}, h^{H({\cal E}, v)} \right)
\text{ in } \overline{\Omega}^{\prime \prime, odd}(M, {\frak B}).
\end{equation}
\end{proposition}
\begin{pf}
This follows from (\ref{partialtcs}) and Propositions \ref{csconv} and 
\ref{smallt}.
\end{pf}

Let us look more closely at ${\cal T}_{[0]}$, the component of ${\cal T}$ in 
$\overline{\Omega}^{\prime \prime, 0}(M, {\frak B})$.
Assume for simplicity that $M$ is connected. Then
\begin{equation}
\overline{\Omega}^{\prime \prime, 0}(M, {\frak B}) =
\left( C^\infty(M)/\C\right) \otimes 
\left({\frak B}/\overline{[{\frak B},{\frak B}]}\right).
\end{equation}
From (\ref{torsiondef}) and (\ref{torsioneqn}),
\begin{align} \label{torequiv}
{\cal T}_{[0]} & \equiv  
- \int_0^\infty 
\left( \int_0^1 \Tr_s \left( N
e^{- tu(1-u)\triangle}  \right) du 
+ t \int_0^1 \int_0^1  u(1-u) \right. \\
& \left. \hspace{.5in} \Tr_s \left(
N \left[ v - v^*,
e^{-ru(1-u)\triangle} 
\left( v - v^* \right) 
e^{-(1-r)u(1-u)\triangle} \right] \right) dr du 
\right)\frac{dt}{t} \notag \\
& = - \int_0^\infty  \int_0^1 \Tr_s \left( N
\left( 1 - 2 t u(1-u) \triangle \right)
e^{- tu(1-u)\triangle}  \right) du 
\frac{dt}{t} \notag
\end{align}

To give a specific lifting of  ${\cal T}_{[0]}$ to 
$C^\infty(M) \otimes \left({\frak B}/\overline{[{\frak B},{\frak B}]}\right)$,
we use the fact that $\Tr_s \left( N \big|_{\cal E} \right)$ and 
$\Tr_s \left( N \big|_{H(E, v)} \right)$ are constant on $M$. 
Put 
\begin{equation} \label{g(t)}
g(t) = - \int_0^1 (1 - 2u(1-u) t)e^{-tu(1-u)} du.
\end{equation}
The asymptotics of $g$ are given by
\begin{equation} \label{gasymp}
g(0) = -1 \text{ and } \lim_{t \rightarrow \infty}
g(t) = 0.
\end{equation}
Then we can
define the lifting to be 
\begin{align}
{\cal T}_{[0]} & =  \int_0^\infty \left[ \Tr_s \left(
N g \left(t \triangle\right) \right) 
- \left( \Tr_s \left( N \big|_{\cal E} \right) -
\Tr_s \left( N \big|_{H(E, v)} \right) \right) g(t) \right. \\
& \left. \hspace{1in} +
\Tr_s \left( N \big|_{H(E, v)} \right) 
\right] \frac{dt}{t}. \notag
\end{align}
Let $\triangle^\prime$ be the restriction of $\triangle$ to 
$\Image(\triangle)$.
Then
\begin{equation}
{\cal T}_{[0]} =  \int_0^\infty \left[ \Tr_s \left(
N g \left(t \triangle^\prime \right) \right) 
- \Tr_s \left( N \big|_{Im(\triangle)} \right)
g(t) \right] \frac{dt}{t}.
\end{equation}
It follows from (\ref{gasymp}) that for $\lambda > 0$,
\begin{equation}
\int_0^\infty \left[ g(\lambda t) - g(t) \right] \frac{dt}{t} = \log
(\lambda).
\end{equation}
Thus by the holomorphic functional calculus,
\begin{equation}
{\cal T}_{[0]} = \Tr_s \left( N \log(\triangle^\prime) \right) \in
C^\infty(M) \otimes \left({\frak B}/\overline{[{\frak B},{\frak B}]}\right).
\end{equation}
{\bf Example 7 : } If ${\frak B} = \C$ then ${\cal T}_{[0]}$ is
the usual Reidemeister torsion of the cochain complex $({\cal E}, v)$
\cite{Ray-Singer (1971)}, considered to be a function on $M$.\\ \\
{\bf Example 8 : } Suppose that $\Gamma$ is a finite group and ${\frak B} =
\C \Gamma$.
Then ${\cal T}_{[0]}$ is equivalent to the equivariant Reidemeister
torsion of \cite{Lott-Rothenberg (1991)}. \\ \\
{\bf Example 9 : } Suppose that $\Gamma$ is a discrete group and
${\frak B} = C^*_r \Gamma$.
Let $\tau$ be the trace on ${\frak B}$ given
by $\tau(\sum_{\gamma \in \Gamma} c_\gamma \gamma) = c_e$. 
Then $\tau \left( {\cal T}_{[0]} \right)$
is the $L^2$-Reidemeister torsion of 
\cite{Carey-Mathai (1990)}, in the case when there is a gap in the
spectrum.  (One can define the
$L^2$-Reidemeister torsion for a cochain complex of modules over the
group von Neumann algebra, not just the group $C^*$-algebra.)

\section{Fiber Bundles} \label{Fiber Bundles}
In this section we extend the results of Section \ref{Noncommutative
Superconnections} to the fiber bundle setting. The translation is that
the algebra of endomorphisms of a finitely-generated projective
${\frak B}$-module gets replaced by an algebra of 
${\frak B}$-pseudodifferential operators.
In the case ${\frak B} = \C$, we recover the 
fiber-bundle results of \cite{Bismut-Lott (1995)}.
We emphasize the necessary modifications to
\cite{Bismut-Lott (1995)} and refer to \cite{Bismut-Lott (1995)} for
some computations.
\subsection{${\frak B}$-Pseudodifferential Calculus}
Let $Z^n$ be a smooth closed manifold.  Let ${\cal E}^1$ and ${\cal E}^2$ 
be smooth 
${\frak B}$-vector bundles on $Z$, with fibers isomorphic to ${\frak E}^1$ and
${\frak E}^2$, respectively.
In the case when ${\frak B}$ is a $C^*$-algebra, an algebra 
$\Psi^{\infty}_{\frak B}(Z; {\cal E}^1, {\cal E}^2)$ of classical 
${\frak B}$-pseudodifferential operators was
defined in \cite[\S 3]{Mischenko-Fomenko (1979)}. We extend this notion to 
the Fr\'echet locally $m$-convex algebra ${\frak B}$ as follows.
First,  define seminorms $\{ \parallel \cdot \parallel
\}_{j=0}^\infty$ on $\Hom_{\frak B}({\frak E}^1,{\frak E}^2)$ similarly to 
the discussion after Proposition
\ref{tt}. Let $U \cong \R^n$ be a coordinate patch of $Z$
equipped with isomorphisms ${\cal E}^1 \big|_U \cong U \times
{\frak E}^1$ and ${\cal E}^2 \big|_U \cong U \times
{\frak E}^2$. We define an algebra $\Psi^{\infty}_{\frak B}
(U; {\cal E}^1, {\cal E}^2)$ of 
classical ${\frak B}$-pseudodifferential operators on $U$ by
requiring that the symbol $\sigma(z, \xi) \in 
\Hom_{\frak B}({\frak E}^1, {\frak E}^2)$
of an order-$m$ operator $T \in \Psi^m_{\frak B}(Z; {\cal E}^1, {\cal E}^2)$
be compactly supported in $z$ and satisfy
\begin{equation}
\parallel \partial_{z^\alpha} \partial_{\xi^\beta} \sigma(z, \xi)
\parallel_j \: \le \: C_{\alpha, \beta, j} (1+|\xi|)^{m - |\beta|}
\end{equation}
for all multi-indices $\alpha$ and $\beta$.
Then we define $\Psi^{\infty}_{\frak B}(Z; {\cal E}^1,{\cal E}^2)$ 
using a partition
of unity as in \cite[\S 3]{Mischenko-Fomenko (1979)}.

Using the representation of ${\frak B}$ as the projective limit of the 
sequence (\ref{projlimit}) of Banach algebras $\{B_j\}_{j=0}^\infty$, with
$B_0$ a $C^*$-algebra,
we can say that $\Psi^{\infty}_{\frak B}(Z; {\cal E}^1, {\cal E}^2)$ 
is the projective limit
of the sequence of pseudodifferential operator algebras
\begin{equation} \label{projlimit2}
\ldots \longrightarrow \Psi^{\infty}_{B_{j+1}} (Z; E^1_{j+1}, E^2_{j+1}) 
\longrightarrow 
\Psi^{\infty}_{B_j} (Z; E^1_j, E^2_j) \longrightarrow \ldots
\longrightarrow \Psi^{\infty}_{B_0} (Z; E^1_0, E^2_0).
\end{equation}   

Let ${\cal E}$ be a ${\frak B}$-vector bundle on $Z$.
Given $T \in \Psi^{\infty}_{\frak B}(Z; {\cal E}, {\cal E})$, let
$i_j(T)$ be its image in $\Psi^{\infty}_{B_j} (Z; E_j, E_j)$.
\begin{proposition} \label{samespectrum}
If $i_0(T)$ is invertible in $\Psi^{\infty}_{B_0} (Z; E_0, E_0)$ then
$T$ is invertible in $\Psi^{\infty}_{\frak B}(Z; {\cal E}, {\cal E})$.
\end{proposition}
\begin{pf}
It is enough to show that each $i_j(T)$ is invertible in 
$\Psi^{\infty}_{B_j} (Z; E_j, E_j)$, as then $T^{-1}$ will be the
inverse limit of $\{ \left( i_j(T) \right)^{-1} \}_{j=0}^\infty$.
So suppose that $B$ is a Banach algebra which is dense in a
$C^*$-algebra $\overline{B}$ and 
stable under the holomorphic functional calculus
in $\overline{B}$. Let $E$ be a $B$-vector bundle on $Z$.
Let $\overline{E} = \overline{B} \otimes_B E$ be the corresponding 
$\overline{B}$-vector bundle on $Z$. Given $T \in \Psi^{m}_{B} (Z; E, E)$, let
$\overline{T}$ be its image in $\Psi^{m}_{\overline{B}} 
(Z; \overline{E}, \overline{E})$. 
We will show that if $\overline{T}$ is invertible in 
$\Psi^{\infty}_{\overline{B}} (Z; \overline{E},\overline{E})$ then
$T$ is invertible in $\Psi^{\infty}_{B} (Z; E,E)$.

Write $E$ as the image under a projection $e \in 
C^\infty(Z; M_N(B))$ of a trivial $B$-vector bundle $Z \times B^N$. Let
$E^\prime = \Image(1-e)$ be the complementary $B$-vector bundle.
Choose $T^\prime \in \Psi^{m}_{B} (Z; E, E)$ such that 
$\overline{T^\prime}$ is invertible in 
$\Psi^{\infty}_{\overline{B}} (Z; \overline{E^\prime},\overline{E^\prime})$. 
If we can
show that $T \oplus T^\prime$ is invertible in $\Psi^{\infty}_{B} (Z; B^N,
B^N)$
then the inverse of $T$ will be given by the restriction of
$\left( T \oplus T^\prime \right)^{-1}$ to $\Image(e)$. 
So we may as well assume that $E$ is a trivial $B$-vector bundle with fiber
$B^N$.

Note that as $\overline{T}$ is invertible, $T$ is elliptic.
By the usual parametrix construction,
we can find $U \in \Psi^{-m}_B(Z; E,E)$ such that $TU = I-K$ with 
$K \in \Psi^{-\infty}_B(Z; E,E)$, and similarly for $UT$.   
Perturbing $U$ a bit if
necessary, we can assume that $\overline{U}$ is invertible. Then
$I - \overline{K}$ is invertible, with inverse $\overline{U}^{-1}
\overline{T}^{-1}$. If we can show that $I - K$ is invertible then
$T^{-1} = U (I - K)^{-1}$. 

Thus we are reduced to showing that if $K \in \Psi^{-\infty}_B(Z; E,E)$
and $I - \overline{K}$ is invertible then $I - K$ is invertible.
Fix a Riemannian metric on $Z$.
Let $\{ e_i \}_{i=1}^\infty$ be the orthonormal basis of smooth functions on 
$Z$ given by the eigenfunctions of the Laplacian.  For $M \ge 1$, let 
$p_M \in \Psi^{-\infty}_{\C}(Z; \C, \C)$ 
be the obvious projection operator from 
$L^2(Z)$ to $\left( \bigoplus_{i = 1}^M  e_i \right)$ and let
$P_M \in \Psi^{- \infty}_B(Z; E,E)$ be its extension by the identity on $B^N$. 
Consider the operator $D_M = I - (I - P_M)  K (I - P_M)$. We claim that if 
$M$ is large enough then $D_M$ is invertible.  To see this, write the
Schwartz kernel of $(I - P_M)  K (I - P_M)$ as
\begin{align}
\left[ (I - P_M)  K (I - P_M) \right](z, z^\prime)  = & 
K(z, z^\prime) - \int_Z P_M(z, w) K(w, z^\prime) dw
- \int_Z K(z, w^\prime) P_M(w^\prime, z^\prime) dw  \\
& +
\int_Z \int_Z P_M(z, w) K(w, w^\prime) P_M(w^\prime, z^\prime) dw dw^\prime.
\notag
\end{align}
The sequence of Schwartz kernels $\{P_M(w, w^\prime)\}_{M=1}^\infty$ forms
an approximate identity. By assumption, $K(z, z^\prime)$ is a smooth
function from $Z \times Z$ to $M_N(B)$.
It follows that for any $\epsilon > 0$, 
there is an $M \ge 1$ such that for all $z, z^\prime \in Z$, in the
Banach norm,
\begin{equation}
| \left[ (I - P_M)  K (I - P_M) \right](z, z^\prime) |
\: \le \epsilon.
\end{equation}
Taking $\epsilon$ small enough, the sum of convolutions
\begin{equation}
D_M^{-1} = \sum_{k=0}^\infty \left( (I - P_M)  K (I - P_M) \right)^k
\end{equation}
converges in the algebra $I + \Psi^{- \infty}_B(Z; E,E)$.

With respect to the decomposition $I = P_M + (I - P_M)$, write
\begin{equation}
I - K = 
\begin{pmatrix}
\alpha & \beta \\
\gamma & \delta
\end{pmatrix}.
\end{equation}
We have shown that $\delta$ is invertible. Then
\begin{equation}
\begin{pmatrix}
\alpha & \beta \\
\gamma & \delta
\end{pmatrix} =
\begin{pmatrix}
1 & \beta \delta^{-1} \\
0 & 1
\end{pmatrix}
\begin{pmatrix}
\alpha - \beta \delta^{-1} \gamma & 0 \\
0 & \delta
\end{pmatrix}
\begin{pmatrix}
1 & 0 \\
\delta^{-1} \gamma & 1
\end{pmatrix}.
\end{equation}
As $I - \overline{K}$ is invertible,
it follows that $\overline{\alpha - \beta \delta^{-1} \gamma}$ 
is invertible in $M_{MN}(\overline{B})$. 
Then $\alpha - \beta \delta^{-1} \gamma$ is invertible in $M_{MN}(B)$ 
\cite[Proposition A.2.2]{Bost (1990)}. Hence
\begin{equation}
(I - K)^{-1} = 
\begin{pmatrix}
1 & 0 \\
- \delta^{-1} \gamma & 1
\end{pmatrix}
\begin{pmatrix}
\left( \alpha - \beta \delta^{-1} \gamma \right)^{-1} & 0 \\
0 & \delta^{-1}
\end{pmatrix}
\begin{pmatrix}
1 & - \beta \delta^{-1} \\
0 & 1
\end{pmatrix}
\end{equation}
is well-defined in 
$I + \Psi^{- \infty}_B(Z; E,E)$.
\end{pf}

Note that $\Psi^{- \infty}_{\frak B}(Z; {\cal E},{\cal E})$ is an algebra
in its own right (without unit). Given $T \in \Psi^{- \infty}_{\frak B}(Z; 
{\cal E},{\cal E})$, let $\sigma_{\Psi^{- \infty}}(T)$ denote its spectrum in
$\Psi^{- \infty}_{\frak B}(Z; {\cal E},{\cal E})$ and let
$\sigma_{\Psi^{\infty}}(T)$ denote its spectrum in
$\Psi^{\infty}_{\frak B}(Z; {\cal E},{\cal E})$.
\begin{lemma} \label{spec1}
$\sigma_{\Psi^{-\infty}}(T) = \sigma_{\Psi^{\infty}}(T)$.
\end{lemma}
\begin{pf}
As $\Psi^{- \infty}_{\frak B}(Z; {\cal E},{\cal E})$ has no unit,
$0 \in \sigma_{\Psi^{- \infty}}(T)$. As $T$ is not invertible in 
$\Psi^{\infty}_{\frak B}(Z; {\cal E},{\cal E})$, 
$0 \in \sigma_{\Psi^{ \infty}}(T)$. If $\lambda \neq 0$ then by definition,
$\lambda \notin \sigma_{\Psi^{- \infty}}(T)$ if and only if we can solve
the equation $TU - \lambda U - \lambda^{-1} T = 
UT - \lambda U - \lambda^{-1} T = 0$ for some $U \in 
\Psi^{- \infty}_{\frak B}(Z; {\cal E},{\cal E})$. Thus if
$\lambda \notin \sigma_{\Psi^{- \infty}}(T)$ then in 
$\Psi^{\infty}_{\frak B}(Z; {\cal E},{\cal E})$, we have
$(T - \lambda)(U - \lambda^{-1}) = (U - \lambda^{-1})(T - \lambda) = I$ 
and hence 
$\lambda \notin \sigma_{\Psi^{\infty}}(T)$. Conversely, if
$\lambda \notin \sigma_{\Psi^{\infty}}(T)$ then we can solve
$(T - \lambda)(U - \lambda^{-1}) = (U - \lambda^{-1})(T - \lambda) = I$ 
for some $U \in
\Psi^{\infty}_{\frak B}(Z; {\cal E},{\cal E})$. By the pseudodifferential
operator calculus, $U = \lambda^{-1} T (T - \lambda)^{-1}
= \lambda^{-1} (T - \lambda)^{-1} T
\in \Psi^{-\infty}_{\frak B}(Z; {\cal E},{\cal E})$
and so $\lambda \notin \sigma_{\Psi^{- \infty}}(T)$.
\end{pf}

Fix a Riemannian metric on $Z$. 
Given a ${\frak B}$-vector bundle ${\cal E}$ on $Z$, write
${\cal E} = {\frak B}^N e$ for some $N$ and some projection $e \in 
C^\infty(Z; M_N({\frak B}))$. Then 
\begin{equation}
\Hom_{\frak B} \left({\cal E}_{z_2},{\cal E}_{z_1} \right) \cong
\{ k \in M_N({\frak B}) : k = e(z_1) k e(z_2)\}.
\end{equation}
Consider the algebra ${\frak A}$ of
integral operators whose kernels $K(z_1, z_2) \in \Hom_{\frak B}
\left({\cal E}_{z_2},{\cal E}_{z_1} \right)$ are continuous in
$z_1$ and $z_2$, with multiplication
\begin{equation}
(K K^\prime)(z_1, z_2) = \int_Z K(z_1,z) K^\prime(z,z_2) \: d\vol(z).
\end{equation}
Let $A_j$ be the analogous algebra with continuous kernels
$K(z_1, z_2) \in \Hom_{B_j}
\left((E_j)_{z_2},(E_j)_{z_1} \right)$. Give $\Hom_{B_j}
\left((E_j)_{z_2},(E_j)_{z_1} \right)$ the Banach space norm $|\cdot|_j$
induced from $\Hom \left( B_j^N, B_j^N \right)$.
Define a norm $| \cdot |_j$ on
$A_j$ by
\begin{equation}
|K|_j = (\vol(Z))^{-1} \max_{z_1, z_2 \in Z} |K(z_1, z_2)|_j.
\end{equation}
Then one can check that $A_j$ is a Banach algebra (without unit).
Furthermore, 
${\frak A}$ is the
projective limit of $\{A_j\}_{j \ge 0}$ and so is
a Fr\'echet locally $m$-convex algebra with seminorms $\{ \parallel \cdot 
\parallel_j \}_{j \ge 0}$ coming from $\{ | \cdot |_j \}_{j \ge 0}$. 
Any $T \in \Psi^{-\infty}_{\frak B}(Z; {\cal E},{\cal E})$ gives an element
of ${\frak A}$ through its Schwartz kernel. Let $\sigma_{\frak A}(T)$ be its
spectrum in ${\frak A}$.
\begin{lemma} \label{spec2}
$\sigma_{\frak A}(T) = \sigma_{\Psi^{- \infty}}(T)$.
\end{lemma}
\begin{pf}
As ${\frak A}$ and $\Psi^{-\infty}_{\frak B}(Z; {\cal E},{\cal E})$ have
no unit, $0 \in \sigma_{\frak A}(T)$ and $0 \in \sigma_{\Psi^{- \infty}}(T)$.
For $\lambda \neq 0$, suppose that $\lambda \notin 
\sigma_{\Psi^{- \infty}}(T)$. Then we
can solve $TU - \lambda U - \lambda^{-1} T = UT - \lambda U - \lambda^{-1} T 
= 0$ for some $U \in 
\Psi^{- \infty}_{\frak B}(Z; {\cal E},{\cal E})$. As $U$ defines an element
of ${\frak A}$, we have $\lambda \notin \sigma_{\frak A}(T)$. Now suppose that
$\lambda \notin \sigma_{\frak A}(T)$. Then we
can solve $TU - \lambda U - \lambda^{-1} T = UT - \lambda U - \lambda^{-1} T 
= 0$ for some $U \in 
{\frak A}$. As $U = \lambda^{-1} T (U - \lambda^{-1}) =
\lambda^{-1} (U - \lambda^{-1}) T$, it follows that
$U$ has a smooth kernel and so defines an element of
$\Psi^{-\infty}_{\frak B}(Z; {\cal E},{\cal E})$. Thus $\lambda \notin
\sigma_{\Psi^{- \infty}}(T)$.
\end{pf}

Define $\TR : {\frak A} \rightarrow {\frak B}/
\overline{[{\frak B},{\frak B}]}$ as in (\ref{trace}).

\begin{corollary} \label{decay}
Suppose that $\{\alpha_r\}_{r > 0}$ is a $1$-parameter semigroup in
$\Psi^{-\infty}_{\frak B}(Z; {\cal E},{\cal E})$ whose spectral radius
in $\Psi^{\infty}_{\frak B}(Z; {\cal E},{\cal E})$ is
given by $\SpRad(\alpha_r) = e^{ar}$ for some $a < 0$. Then  for all
$j \ge 0$, as $r \rightarrow \infty$, in
${\frak A}$ we have $\parallel \alpha_r \parallel_j = o \left( 
e^{ar/2} \right)$. In particular, $\TR \left( \alpha_r \right) = o \left( 
e^{ar/2} \right)$.
\end{corollary}
\begin{pf}
By Lemmas \ref{spec1} and \ref{spec2}, the spectral radius of
$\alpha_r$ in ${\frak A}$ is $e^{ar}$. Then its spectral radius in
the Banach algebra
$A_j$ is less than or equal to $e^{ar}$. By (\ref{start1}) and 
(\ref{start2}), $\parallel \alpha_r \parallel_j = o \left( 
e^{ar/2} \right)$. As $\TR$ is continuous on ${\frak A}$, the
corollary follows. 
\end{pf}
 
\subsection{Induced Superconnections} \label{Induced Superconnections}
Let $Z \rightarrow M \stackrel{\pi}{\rightarrow} B$ be a smooth
fiber bundle with total space $M$, compact base $B$ and 
connected closed fibers $\{Z_b\}_{b \in B}$ of dimension
$n$. We use the notation of
\cite[Section 3]{Bismut-Lott (1995)} when discussing the topology or geometry
of such a fiber bundle. Let $TZ$ be the vertical tangent bundle of the
fiber bundle, an $\R^n$-bundle on $M$, and let $o(TZ)$ be its orientation
bundle, a flat $\R$-bundle on $M$.
Let $T^H M$ be a horizontal
distribution for the fiber bundle. Let ${\cal E}$ be a ${\frak B}$-vector
bundle on $M$. There is an induced $\Z$-graded
${\frak B}$-vector bundle ${\cal W}$ on $B$ whose fiber over $b \in B$
consists of the smooth ${\cal E}$-valued differential forms on $Z_b$, i.e.
${\cal W}_b = \Omega^* \left(Z_b; {\cal E} \big|_{Z_b} \right)$. If
$n > 0$ then ${\cal W}_b$ is infinitely-generated.
Using the horizontal distribution, there is an isomorphism
\begin{equation} \label{iso}
\Omega(B, {\frak B}; {\cal W}) \cong \Omega(M, {\frak B}; {\cal E}).
\end{equation}

We equip ${\cal E}$ with a partially flat connection $\nabla^{\cal E}$ as in 
(\ref{connection}). Using (\ref{iso}), this induces
a partially flat degree-$1$ superconnection $A^\prime$ on ${\cal W}$.
The connection component $\nabla^{\cal W}$ of $A^\prime$ has two pieces :
\begin{align}
&\nabla^{{\cal W},1,0} : \Omega^{p,q}(M, {\frak B}; {\cal E}) \rightarrow
\Omega^{p+1,q}(M, {\frak B}; {\cal E}), \\
&\nabla^{{\cal W},0,1} : 
\Omega^{p,q}(M, {\frak B}; {\cal E}) \rightarrow 
\Omega^{p,q+1}(M, {\frak B}; {\cal E}). \notag
\end{align} 
As in
\cite[Proposition 3.4]{Bismut-Lott (1995)}, $\nabla^{{\cal W},1,0}$ is
given by Lie differentiation with respect to a horizontal vector field on
$M$. On the other hand, $\nabla^{{\cal W},0,1}$ comes from the action of
$\partial^{\cal E}$ as in (\ref{partialdef}). 
The other nonzero components of $A^\prime$ are 
$A^\prime_{0,0,1} = d^Z$ and $A^\prime_{2,0,-1} = i_T$. The degree-$1$
superconnection $A^{\prime,flat}$ defining the superflat
structure on ${\cal W}$ is essentially the same as the flat
degree-$1$ superconnection of \cite[Section 3b]{Bismut-Lott (1995)}.
The main difference between \cite[Section 3]{Bismut-Lott (1995)} 
and the present paper is that we take into account 
$\nabla^{{\cal W},0,1}$, so that $A^\prime$ is not completely flat.

Let $g^{TZ}$ be a family of vertical Riemannian metrics on the fiber
bundle.  Let $*$ be the 
corresponding fiberwise Hodge duality operator, extended
linearly from $C^\infty(M; \Lambda(T^* Z))$ to
$C^\infty(M; \Lambda(T^* Z) \otimes {\cal E}) \cong C^\infty(B; {\cal W})$. 
Let $h^{\cal E}$ be a Hermitian metric on ${\cal E}$.
Let $\left( \nabla^{\cal E} \right)^*$ be the adjoint connection to
$\nabla^{\cal E}$, with respect to $h^{\cal E}$. 
There is a self-adjoint connection on ${\cal E}$ given by
\begin{equation}
\nabla^{{\cal E},sa} = \frac12 \left( \nabla^{{\cal E}} + 
\left( \nabla^{{\cal E}} \right)^* \right).
\end{equation}
Put 
\begin{equation}
\psi = \left( \nabla^{{\cal E},1,0} \right)^* - \nabla^{{\cal E},1,0}
\in \Omega^1(M; \End_{\frak B}({\cal E})). 
\end{equation}
Let us assume that $\left( \partial^{\cal E} \right)^* = \partial^{\cal E}$;
this can always be achieved by replacing $\partial^{\cal E}$ by
$\frac12 \left( \partial^{\cal E} + \left( \partial^{\cal E} \right)^*
\right)$ if necessary.
Then
\begin{equation}
\nabla^{{\cal E},sa} = \nabla^{{\cal E}} + \frac{\psi}{2}.
\end{equation}

For notational convenience, put 
\begin{equation}
\nabla^{{\cal E}, flat, sa} = \frac12 \left( \nabla^{{\cal E},flat} + 
\left( \nabla^{{\cal E},flat} \right)^* \right).
\end{equation}
Then
\begin{equation}
\nabla^{{\cal E},flat, sa} = \nabla^{{\cal E}, flat} + \:
\frac{\psi}{2},
\end{equation}
\begin{equation}
\nabla^{{\cal E}, sa} = \nabla^{{\cal E}, flat,sa} + \partial^{\cal E}
\end{equation}
and
\begin{equation}
\left( \nabla^{{\cal E}, flat, sa} \right)^2 =
- \: \frac{\psi^2}{4}.
\end{equation}
Furthermore, one can check that 
$\nabla^{TM \otimes {\cal E}, flat, sa} \psi$ vanishes in 
$\Omega^2 \left( M; \End_{\frak B}({\cal E})\right)$. 

There is a Hermitian metric $h^{\cal W}$ on ${\cal W}$ such that for
$s, s^\prime \in {\cal W}_b$, 
\begin{equation}
\langle s, s^\prime \rangle_{h^{\cal W}} = \int_{Z_b}
\langle s \wedge * s^\prime \rangle_{h^{\cal E}} \in {\frak B}.
\end{equation}
Let $(A^\prime)^*$ be the adjoint superconnection to $A^\prime$, with
respect to $h^{\cal W}$. 
There is a self-adjoint connection $\nabla^{{\cal W},sa}$ on ${\cal W}$ given 
by
\begin{equation}
\nabla^{{\cal W},sa} = \frac12 \left( \nabla^{{\cal W}} + 
\left( \nabla^{{\cal W}} \right)^* \right).
\end{equation}

Let $\{e_j\}_{j=1}^n$ be a local orthonormal basis for $TZ$, with dual basis  
$\{\tau^j\}_{j=1}^n$. Let $E^j$ denote exterior multiplication by $\tau^j$ and
let $I^j$ denote interior multiplication by $e_j$. Put
\begin{align} \label{cc}
c^j & = E^j - I^j, \\
\widehat{c}^j & = E^j + I^j.
\notag
\end{align}
Then 
\begin{align}
c^j c^k + c^k c^j & = - 2 \delta^{jk},\\
\widehat{c}^j \widehat{c}^k + \widehat{c}^k \widehat{c}^j & =  
2 \delta^{jk}, \notag \\
c^j \widehat{c}^k + \widehat{c}^k c^j & = 0. \notag
\end{align}
Thus $c$ and $\widehat{c}$ generate two graded-commuting Clifford algebras.

Let $\nabla^{TZ}$ be Bismut's connection on $TZ$ 
\cite[p. 322]{Berline-Getzler-Vergne (1992)}, with curvature
$R^{TZ}$. Let $e \left( TZ, \nabla^{TZ} \right) \in \Omega^n(M; o(TZ))$ be the
corresponding Euler form. Define ${\cal R} \in 
\Omega^2 \left( M; \End (\Lambda(T^*Z) \otimes {\cal E})\right)$ by
\begin{equation}
{\cal R} = \frac14 \left( \langle e_j, R^{TZ} e_k \rangle_{g^{TZ}} \: 
\widehat{c}^j \widehat{c}^k \otimes I_{\cal E} \right) - \frac14
\left( I_{\Lambda(T^*Z)} \otimes \psi^2 \right). 
\end{equation}
Let $R \in C^\infty(M)$ be the scalar curvature of the fibers. 
Let ${\tau^\alpha}$ be a local basis of $T^*B$ and let $E^\alpha$
denote exterior multiplication by $\tau^\alpha$.
Define the superconnection $B_t(u)$ on ${\cal W}$ as in Proposition \ref{tt}.
\begin{proposition} \label{B_t(u)}
We have
\begin{align}
B_t(u) \: = \: & \frac{\sqrt{t}}{2} \: c^j 
\nabla^{TZ \otimes {\cal E},sa}_{e_j} -
\left( \frac12 - u \right) \sqrt{t} \:
\widehat{c}^j \nabla^{TZ \otimes {\cal E},sa}_{e_j}
- \frac{\sqrt{t}}{4} \: \widehat{c}^j \psi_j + 
 \left( \frac12 - u \right) \frac{\sqrt{t}}{2} \:
c^j \psi_j \\
& + \nabla^{{\cal W},sa} + \left( \frac12 - u \right) E^\alpha \psi_\alpha
+ \left( \frac12 - u \right) \omega_{\alpha j k} E^\alpha c^j \widehat{c}^k 
\notag \\
& + \left( \frac12 - u \right) \: \frac{1}{\sqrt{t}} \: 
\omega_{\alpha \beta j} 
E^\alpha E^\beta \widehat{c}^j + \frac{1}{2\sqrt{t}} \: 
\omega_{\alpha \beta j} E^\alpha E^\beta c^j, \notag
\end{align}
where
\begin{align}
\nabla^{{\cal W},sa,1,0} & = E^\alpha \left( 
\nabla^{TZ \otimes {\cal E},sa}_{e_\alpha} + \: \frac12 k_\alpha \right), \\ 
\nabla^{{\cal W},sa,0,1} & = \partial^{\cal E}. \notag
\end{align}
\end{proposition}
\begin{pf} This follows from a computation using
\cite[Prop. 3.5, 3.7]{Bismut-Lott (1995)}. We omit the details.
\end{pf}

Let $z$ be an odd Grassmann variable which anticommutes with all
of the Grassmann variables previously introduced.
Put 
\begin{align}
{\cal D}_j & = \nabla^{TZ \otimes {\cal E},sa}_{e_j} -
\frac{1}{2 \sqrt{t}} \: \omega_{\alpha j k} E^\alpha c^k -
\frac{1}{4t} \: \omega_{\alpha \beta j} E^\alpha E^\beta
- \frac{1}{2\sqrt{t}} \: z \widehat{c}^j, \\
{\cal D}^2 & = 
{\cal D}_j {\cal D}_j - {\cal D}_{\nabla^{TZ}_{e_j} e_j}. \notag 
\end{align}
\begin{proposition} \label{Lichno}
The following Lichnerowicz-type formula holds :
\begin{align} \label{Lichnomess}
B_{4t}^2(u) + 2 u (1-u) z \left( B_{4t}^\prime - B_{4t}^{\prime \prime}
\right) \: = \: & 
4 u(1-u) \left[ t \left( - {\cal D}^2 + \frac{R}{4} \right) \right. \\
& +
\frac{t}{2} c^i c^j {\cal R}(e_i, e_j) + \sqrt{t} c^i E^\alpha
{\cal R}(e_i, e_\alpha) + \frac12 E^\alpha E^\beta {\cal R}
(e_\alpha, e_\beta) \notag \\
& + t \left( \frac14 \psi_j^2 + \frac18 \widehat{c}^j \widehat{c}^k
[\psi_j, \psi_k] - \frac12 c^j \widehat{c}^k
\left( \nabla^{TZ \otimes {\cal E},sa}_{e_j} \psi_k \right) \right) \notag \\
& \left. - \frac{\sqrt{t}}{2} E^\alpha \widehat{c}^j 
\left( \nabla^{TZ \otimes {\cal E},sa}_{e_\alpha} \psi_j \right) 
- \frac{z \sqrt{t}}{2} c^j \psi_j - \frac{z}{2} E^\alpha \psi_\alpha \right]
\notag \\
& + \sqrt{t} \: c^j 
\left( \nabla^{{\cal E},sa}_{e_j} \partial^{\cal E} \right) - 2
\left( \frac12 - u
\right) \sqrt{t} \: 
\widehat{c}^j \left( \nabla^{{\cal E},sa}_{e_j} \partial^{\cal E} \right)
\notag \\
& + E^\alpha \left(\nabla^{{\cal E},sa}_{e_\alpha} \partial^{\cal E} \right)
+ \left( \frac12 - u \right)
E^\alpha [\psi_\alpha, \partial^{\cal E}]
 - \frac{\sqrt{t}}{2} \: \widehat{c}^j [\psi_j, \partial^{\cal E}] \notag \\
& + \sqrt{t}
\left( \frac12 - u \right) c^j [\psi_j, \partial^{\cal E}]
+ \left( \partial^{\cal E} \right)^2. \notag
\end{align}
\end{proposition}
\begin{pf} Let us write $B_{4t}(u) = u B^{\prime,flat}_{4t} + (1-u)
B^{\prime \prime,flat}_{4t} + \partial^{\cal E}$. Then
\begin{align}
B_{4t}^2(u) = & u(1-u) \left( B^{\prime,flat}_{4t} B^{\prime \prime,flat}_{4t}
 + B^{\prime \prime,flat}_{4t} B^{\prime,flat}_{4t} \right) \\
& + u 
\left[ B^{\prime, flat}_{4t}, \partial^{\cal E} \right] + (1-u) 
\left[ B^{\prime \prime, flat}_{4t}, \partial^{\cal E} \right]
+ \left( \partial^{\cal E} \right)^2. \notag
\end{align} 
A formula for $\frac14 \left( 
B^{\prime,flat}_{4t} B^{\prime \prime,flat}_{4t} + 
B^{\prime \prime,flat}_{4t} B^{\prime,flat}_{4t} \right)  + \frac12 \: z
\left(B^{\prime,flat}_{4t} - B^{\prime \prime,flat}_{4t} \right)$ was given in
\cite[Theorem 3.11]{Bismut-Lott (1995)}.
The rest of (\ref{Lichnomess}) can be derived using Proposition
\ref{B_t(u)}. 
\end{pf}
\subsection{Small Time Limits}

For $t > 0$ and $u \in (0,1)$, the restriction of $B_t^2(u)$ to
a fiber $Z_b$ is an element of $\Psi^{2}_{\frak B} \left(Z_b;
\Lambda(T^*Z_b) \otimes {\cal E}\big|_{Z_b}, 
\Lambda(T_b^*B) \otimes \Omega_*({\frak B}) \otimes_{\frak B}
(\Lambda(T^*Z_b) \otimes {\cal E}\big|_{Z_b}) \right)$ with principal
symbol $\sigma(z, \xi) = u(1-u)t |\xi|^2$.  It follows that on $Z_b$,
\begin{equation} 
e^{-B_t^2(u)} \in \Psi^{-\infty}_{\frak B} \left( Z_b;
\Lambda(T^*Z_b) \otimes {\cal E}\big|_{Z_b}, 
\Lambda(T_b^*B) \otimes \Omega_*({\frak B}) \otimes_{\frak B}
(\Lambda(T^*Z_b) \otimes {\cal E}\big|_{Z_b}) \right).
\end{equation}
Hence $e^{-B_t^2(u)}$ has a smooth kernel $e^{-B_t^2(u)}(z, z^\prime)$ and
using the notion of $\TR_s$ from Subsection \ref{Chern-Simons Classes}, 
we can define
$\TR_s \left( e^{-B_t^2(u)} \right) \in 
\overline{\Omega}^{even}(B, {\frak B})$.

Put $\nabla^{\cal E}(u) = u \nabla^{\cal E} + (1-u) \left(
\nabla^{\cal E} \right)^*$.
\begin{proposition} \label{Chernlimit}
For all $u \in (0,1)$, as $t \rightarrow 0$,
\begin{equation} \label{ttt}
\TR_s \left( e^{-B_t^2(u)} \right) =
\begin{cases}
\int_Z e \left(TZ, \nabla^{TZ} \right) \wedge \ch \left(
\nabla^{\cal E}(u) \right) + O(t) & \text{if $n$ is even,}\\
O(\sqrt{t}) & \text{if $n$ is odd}
\end{cases}
\end{equation}
uniformly on $B$.
\end{proposition}
\begin{pf}
Consider a rescaling in which
$\partial_j \rightarrow \epsilon^{-1/2} \partial_j$, $c^j \rightarrow 
\epsilon^{-1/2} E^j - \epsilon^{1/2} I^j$, $E^\alpha \rightarrow
\epsilon^{-1/2} E^\alpha$, $\widehat{c}^j \rightarrow \widehat{c}^j$ and
$\partial^{\cal E} \rightarrow \epsilon^{-1/2} \partial^{\cal E}$. One
finds from (\ref{Lichnomess}) that as $\epsilon \rightarrow 0$, 
in adapted coordinates the rescaling of $\epsilon B_4^2(u)$ approaches
\begin{align} 
& - 4 u (1-u) \left( \partial_j - \frac14 \: R^{TZ}_{jk} x^k \right)^2 
+ 4 u (1-u) {\cal R} \\
& +  E^j 
\left( \nabla^{{\cal E},sa}_{e_j} \partial^{\cal E} \right)
+ E^\alpha \left( \nabla^{{\cal E},sa}_{e_\alpha} \partial^{\cal E} \right)
\notag \\
& + \left( \frac12 - u \right) E^j [\psi_j, \partial^{\cal E}] 
+ \left( \frac12 - u \right) E^\alpha [\psi_\alpha, \partial^{\cal E}]
+ \left( \partial^{\cal E} \right)^2. \notag
\end{align}
Using local index methods as in the proof of
\cite[Theorem 3.15]{Bismut-Lott (1995)}, one finds
\begin{align}
\lim_{t \rightarrow 0} \TR_s \left( e^{-B_t^2(u)} \right) & =
\int_Z 
\left( 4u (1-u) \right)^{-n/2} \Pf \left[ 
\frac{4u(1-u) R^{TZ}}{2\pi} \right] \wedge \\
& \hspace{.5in} \Tr_s e^{- \left[   
 \left( \nabla^{{\cal E},flat,sa} \partial^{\cal E} \right)
+ \left( \frac12 - u \right) [\psi, \partial^{\cal E}] 
+ \left( \partial^{\cal E} \right)^2 - u(1-u) \psi^2 \right]} \notag \\
& = \int_Z 
e \left( TZ, \nabla^{TZ} \right) \wedge \Tr_s
e^{- \left[ \left( \nabla^{{\cal E},flat,sa} \partial^{\cal E} \right)
+ \left( \frac12 - u \right) [\psi, \partial^{\cal E}] 
+ \left( \partial^{\cal E} \right)^2 - u(1-u) \psi^2 \right]}. \notag
\end{align}
On the other hand,
\begin{equation}
\nabla^{\cal E}(u) = \nabla^{{\cal E},flat,sa} + \left( \frac12 - u \right)
\psi + \partial^{\cal E}
\end{equation}
and so
\begin{equation}
\left( \nabla^{\cal E}(u) \right)^2 = \left( \nabla^{{\cal E},flat,sa} 
\partial^{\cal E} \right)
+ \left( \frac12 - u \right) [\psi, \partial^{\cal E}] 
+ \left( \partial^{\cal E} \right)^2 - u(1-u) \psi^2.
\end{equation}
This gives the $t \rightarrow 0$ limit of (\ref{ttt}).

We have error estimates as in \cite[Theorem 3.16]{Bismut-Lott (1995)}, 
from which the proposition follows.
\end{pf}

Define $CS \left( B_t^\prime, h^{\cal W} \right) \in 
\overline{\Omega}^{\prime \prime, odd}(B, {\frak B})$ and
${\cal T}(t) \in \overline{\Omega}^{\prime \prime, even}(B, {\frak B})$
 as in Proposition \ref{tt}. 
\begin{proposition} \label{cslimit}
As $t \rightarrow 0$,
\begin{equation}
CS \left( B_t^\prime, h^{\cal W} \right) =
\begin{cases}
\int_Z e \left( TZ, \nabla^{TZ} \right) \wedge 
CS \left( \nabla^{\cal E}, h^{\cal E} \right) + O(t) & \text{if $n$ is even,}\\
O(\sqrt{t}) & \text{if $n$ is odd}
\end{cases}
\end{equation}
uniformly on $B$.
\end{proposition}
\begin{pf}
Given $\alpha_1, \alpha_2 \in \overline{\Omega}^{\prime \prime, *}
(B, {\frak B})$, let us write
\begin{equation}
\partial_z \left( \alpha_1 + z \alpha_2 \right) = \alpha_2.
\end{equation}
Then
\begin{equation}
CS \left( B_t^\prime, h^{\cal W} \right) = \frac12 \: \partial_z \int_0^1
\frac{1}{u(1-u)} \: 
\Tr_s e^{- \left[ B^2_t(u) + 2 u (1-u) z (B^\prime_t - 
B^{\prime \prime}_t) \right]} du.
\end{equation}
Let us do a rescaling as in the proof of Proposition 
\ref{Chernlimit}, with $z \rightarrow \epsilon^{-1/2} z$ in addition.
One
finds from (\ref{Lichnomess}) that as $\epsilon \rightarrow 0$, 
in adapted coordinates the rescaling of $\epsilon 
\left(B_4^2(u) + 2 u (1-u) z (B^\prime_4 - B^{\prime \prime}_4)\right)$ 
approaches
\begin{align} 
& - 4 u (1-u) \left( \partial_j - \frac14 \: R^{TZ}_{jk} x^k \right)^2 
+ 4 u (1-u) {\cal R} - 2 u (1-u) z \psi \\
& +  E^j 
\left( \nabla^{{\cal E},sa}_{e_j} \partial^{\cal E} \right)
+ E^\alpha \left( \nabla^{{\cal E},sa}_{e_\alpha} \partial^{\cal E} \right)
\notag \\
& + \left( \frac12 - u \right) E^j [\psi_j, \partial^{\cal E}] 
+ \left( \frac12 - u \right) E^\alpha [\psi_\alpha, \partial^{\cal E}]
+ \left( \partial^{\cal E} \right)^2. \notag
\end{align}
Proceeding as in the proof of \cite[Theorem 3.16]{Bismut-Lott (1995)},
one obtains
\begin{align}
\lim_{t \rightarrow 0}
CS \left( B_t^\prime, h^{\cal W} \right) & = \frac12 \: \partial_z \int_0^1
\frac{1}{u(1-u)} \: \int_Z e \left( TZ, \nabla^{TZ} \right) \wedge \Tr_s
e^{- \left[ \left( \nabla^{\cal E}(u) \right)^2 - 2u (1-u) z \psi  \right]} du
\\
& = \int_0^1
\int_Z e \left( TZ, \nabla^{TZ} \right) \wedge \Tr_s \left[ \psi \:
e^{- \left( \nabla^{\cal E}(u) \right)^2} \right] du \notag \\
& = \int_Z e \left( TZ, \nabla^{TZ} \right) \wedge
CS \left( \nabla^{\cal E}, h^{\cal E} \right). \notag
\end{align}
Although we are integrating over $u$, there is no problem with 
the $t \rightarrow 0 $ limit 
as the effective time parameter is $u (1-u)t$, which
only improves the convergence.
\end{pf}
\begin{proposition} \label{smallt2} 
As $t \rightarrow 0$,
\begin{equation}
{\cal T}(t) =
\begin{cases}
O(1) & \text{if $n$ is even,}\\
O(t^{-1/2}) & \text{if $n$ is odd}
\end{cases}
\end{equation}
uniformly on $B$.
\end{proposition}
\begin{pf}
Using the method of proof of \cite[Theorem 3.21]{Bismut-Lott (1995)},
one finds
\begin{equation}
{\cal T}(t) =
\begin{cases}
- \:\frac{n}{2} \:  
\frac{1}{t} \int_0^1 \Tr_s \left( e^{-B_t^2(u)} \right) du  + O(1) & 
\text{if $n$ is even,}\\
O(t^{-1/2}) & \text{if $n$ is odd}.
\end{cases}
\end{equation}
By Proposition \ref{Chernlimit}, if $n$ is even then
\begin{equation}
\lim_{t \rightarrow 0} \int_0^1 \Tr_s \left( e^{-B_t^2(u)} \right) du =
\int_Z e \left(TZ, \nabla^{TZ} \right) \wedge \int_0^1 \ch \left(
\nabla^{\cal E}(u) \right) du
\in \overline{\Omega}^{\prime, even}(B, {\frak B}).
\end{equation}
(Again, as the effective time parameter is $u(1-u)t$, there is no problem
in switching the $t \rightarrow 0$ limit and the $u$-integration.) 
As we quotient by 
$\overline{\Omega}^{\prime, even}(B, {\frak B})$ in defining
$\overline{\Omega}^{\prime \prime, even}(B, {\frak B})$, the proposition
follows.
\end{pf}

\subsection{Index Theorems} \label{Index Theorems}
We continue with the setup of Subsection \ref{Induced Superconnections}.
For each $b \in B$, let $H \left(Z_b; {\cal E} \big|_{Z_b} \right)$ denote
the cohomology of the complex $\left( {\cal W}_b, d^Z \right)$.
Put $\triangle_b = d^Z \left( d^Z \right)^* + \left( d^Z \right)^* d^Z
\in \Psi^{2}_{\frak B} \left(Z_b;
\Lambda(T^*Z_b) \otimes {\cal E}\big|_{Z_b}, 
\Lambda(T^*Z_b) \otimes {\cal E}\big|_{Z_b} \right)$.
Put $\overline{\cal E} = \Lambda \otimes_{\frak B} {\cal E}$ and
$\overline{\cal W}_b =
\Omega \left( Z_b; \overline{\cal E} \big|_{Z_b}
\right)$.
Let $\overline{\triangle}_b \in \Psi^{2}_{\Lambda} \left(Z_b;
\Lambda(T^*Z_b) \otimes \overline{\cal E} \big|_{Z_b}, 
\Lambda(T^*Z_b) \otimes \overline{\cal E}\big|_{Z_b} \right)$ be the
corresponding Laplacian in the $\Lambda$-pseudodifferential operator
calculus.

\begin{hypothesis} \label{gaphypothesis}
For each $b \in B$, the operator $\overline{d^Z} \in \End_{\Lambda}
\left(\overline{\cal W}_b\right)$ has closed image.
\end{hypothesis}

\begin{proposition} \label{gap}
Hypothesis \ref{gaphypothesis} is satisfied if and only if $0$ is isolated
in $\sigma( \overline{\triangle}_b)$.
\end{proposition}
\begin{pf}
This follows from standard arguments.  We omit the details.
\end{pf}

Hereafter, we assume that Hypothesis \ref{gaphypothesis} is satisfied.
\begin{proposition} \label{projmodules}
For each $b \in B$, $\HH \left(Z_b; {\cal E} \big|_{Z_b} \right)$ is a
finitely-generated projective ${\cal B}$-module.  The
$\{\HH \left(Z_b; {\cal E} \big|_{Z_b} \right)\}_{b \in B}$ fit together
to form a $\Z$-graded ${\frak B}$-vector bundle 
$H \left(Z; {\cal E} \big|_{Z} \right)$ on $B$ with a flat
structure.
\end{proposition}
\begin{pf}
The proof is similar to the proofs of Propositions \ref{projmodule} and 
\ref{fittogether}, with Proposition \ref{samespectrum} replacing Lemma 
\ref{smalllemma}. 
\end{pf}

There is an induced Hermitian metric 
$h^{H (Z; {\cal E} \big|_{Z} )}$ on 
$H \left(Z; {\cal E} \big|_{Z} \right)$ and an induced partially flat
connection $\nabla^{H (Z; {\cal E} \big|_{Z})}$ as in
(\ref{inducedconn}). 

\begin{proposition} \label{largetime}
For all $u \in (0,1)$, as $t \rightarrow \infty$,
\begin{equation}
\ch \left(B_t(u) \right) = 
\ch \left( \nabla^{H (Z; {\cal E} \big|_{Z} )}(u)
\right) + O(t^{-1/2})
\end{equation}
uniformly on $B$. Also,
\begin{equation}
CS\left(B^\prime_t, h^{\cal W}\right) = 
CS\left( \nabla^{H (Z; {\cal E} \big|_{Z} )},
h^{H (Z; {\cal E} \big|_{Z} )} 
\right) + O(t^{-1/2})
\end{equation}
uniformly on $B$.
\end{proposition}
\begin{pf}
Let $\lambda_0 > 0$ be the infimum of the nonzero spectrum of 
$\overline{\triangle}$.
For $r > 0$, put $\alpha_r = P^{Im(\triangle)} e^{-r \triangle}
P^{Im(\triangle)}$. By Proposition \ref{samespectrum} and Corollary
\ref{decay}, for each $j \ge 0$ there is a constant
$C_j > 0$ such that for all $r > 1$, 
\begin{equation} \label{newend1}
\parallel  \alpha_r \parallel_j \: \le \: C_j 
\: e^{-r \lambda_0/2}.
\end{equation}
The proof of the proposition is now formally the same as that of
Proposition \ref{csconv}, with (\ref{newend1}) replacing (\ref{end1}).
\end{pf}

\begin{proposition} \label{indextheorem}
We have
\begin{equation}
\left[ \ch \left( \nabla^{H (Z; {\cal E} \big|_{Z} )}
\right) \right] = \int_Z e(TZ) \wedge \left[\ch
\left(\nabla^{\cal E} \right) \right]
\text{ in } H^{\prime,even}_{\frak B}(B)
\end{equation}
and
\begin{equation}
\left[ CS \left( \nabla^{H (Z; {\cal E} \big|_{Z} )},
h^{H (Z; {\cal E} \big|_{Z} )}
\right) \right] = \int_Z e(TZ) \wedge 
\left[CS \left(\nabla^{\cal E}, h^{\cal E} \right)\right]
\text{ in }
\bigoplus 
\begin{Sb}
p > q \\
p + q \: odd
\end{Sb}
\HH^p(B; \overline{H}_q({\frak B})). 
\end{equation}
\end{proposition}
\begin{pf}
As in the finite-dimensional setting, one can verify that
$\left[ \ch \left(B_t(u) \right)  \right]$ and $\left[ CS \left( 
B^\prime_t, h^{\cal W} \right) \right]$ are independent of $t$.
For all $u \in (0,1)$, 
$\left[ \ch \left( \nabla^{H (Z; {\cal E} \big|_{Z} )}(u)
\right) \right] = \left[ \ch \left( \nabla^{H (Z; {\cal E} \big|_{Z} )}
\right) \right]$ and
$\left[\ch \left(\nabla^{\cal E}(u) \right) \right] = \left[\ch
\left(\nabla^{\cal E} \right) \right]$.
The proposition now follows from Propositions \ref{Chernlimit},
\ref{cslimit} and \ref{largetime}.
\end{pf}
{\bf Remark 8 : }
Proposition \ref{indextheorem} is also a consequence of the
topological index theorem of \cite{Dwyer-Weiss-Williams (1995)}.
Namely, Proposition \ref{projmodules} ensures that we can apply
(0-3) of their paper, as given in (\ref{DWWthm}) of the present paper.
Proposition \ref{indextheorem} follows from (\ref{DWWthm}) by applying
$\ch$ and $CS$.
 
\subsection{The Analytic Torsion Form II}
We continue with the assumptions of Subsection \ref{Index Theorems}.
Let $N$ be the number operator on ${\cal W}$. For $t > 0$, define
${\cal T}(t)$ as in (\ref{torsiondef}).
\begin{proposition} \label{torsionconv2}
As $t \rightarrow \infty$,
\begin{equation} \label{torsionasymp2}
{\cal T}(t) = O(t^{-3/2})
\end{equation}
uniformly on $B$.
\end{proposition}
\begin{pf}
The proof is formally the same as that of Proposition \ref{torsionconv}. 
We omit the details.
\end{pf}

Again, we have
\begin{equation} \label{partialtcs2}
\partial_t CS\left(B^\prime_t, h^{\cal W} \right) =
- \: d {\cal T}(t).
\end{equation}
\begin{definition}
The analytic torsion form
${\cal T} \in \overline{\Omega}^{\prime \prime, even}(M, {\frak B})$
is given by
\begin{equation} \label{torsioneqn2}
{\cal T} = \int_0^\infty {\cal T}(t) \: dt.
\end{equation}
\end{definition}
By Propositions \ref{smallt2} and \ref{torsionconv2}, the integral in
(\ref{torsioneqn2}) makes sense.
\begin{proposition} \label{dtor}
We have
\begin{equation}
d {\cal T} = \int_Z e \left(TZ, \nabla^{TZ} \right) \wedge
CS \left( \nabla^{\cal E}, h^{\cal E} \right) -
CS \left( \nabla^{H (Z; {\cal E} \big|_{Z} )},
h^{H (Z; {\cal E} \big|_{Z} )} \right)
\text{ in } \overline{\Omega}^{\prime \prime, odd}(B, {\frak B}).
\end{equation}
\end{proposition}
\begin{pf}
This follows from Proposition \ref{cslimit}, Proposition \ref{largetime}
and (\ref{partialtcs2}).
\end{pf}

\begin{corollary} \label{indep}
If $Z$ is odd-dimensional and $H (Z; {\cal E} \big|_{Z}) = 0$ then
${\cal T}$ is closed and so represents a class
$[{\cal T}] \in H^{\prime \prime,even}(B,{\frak B})$.
\end{corollary}
\begin{pf}
If $Z$ is odd-dimensional then $e \left(TZ, \nabla^{TZ} \right) = 0$. The
corollary now follows from Proposition \ref{dtor}.
\end{pf}

\begin{proposition} \label{cclass}
If $Z$ is odd-dimensional and $H (Z; {\cal E} \big|_{Z} ) = 0$ then
$[{\cal T}] \in H^{\prime \prime,even}(B,{\frak B})$ 
is independent of $g^{TZ}$, $T^HM$, $h^{\cal E}$ and
$\partial^{\cal E}$. Thus it only depends on the (smooth) topological
fiber bundle $Z \rightarrow M \rightarrow B$ and the flat structure
on ${\cal E}$.
\end{proposition}
\begin{pf}
Put ${\cal F} = \{g^{TZ}, T^HM, h^{\cal E},\partial^{\cal E} \}$ and let
${\cal F}^\prime$ be another choice of such data.
We can find a smooth $1$-parameter family 
$\{{\cal F}(\epsilon)\}_{\epsilon \in \R}$
such that ${\cal F}(0) = {\cal F}$ and ${\cal F}(1) = {\cal F}^\prime$.
Put $\widetilde{Z} = Z$, $\widetilde{M} = \R \times M$ and 
$\widetilde{B} = \R \times B$.
Let $p : \widetilde{M} \rightarrow M$ be projection onto the second factor
and put $\widetilde{\cal E} = p^* {\cal E}$. Then the family 
$\{{\cal F}(\epsilon)\}_{\epsilon \in \R}$ provides the data
$g^{T\widetilde{Z}}$, $T^H\widetilde{M}$, $h^{\widetilde{\cal E}}$ and
$\partial^{\widetilde{\cal E}}$ on the fiber bundle
$\widetilde{Z} \rightarrow  
\widetilde{M} \stackrel{\widetilde{\pi}}{\rightarrow} 
\widetilde{B}$. Let $\widetilde{d} = d\epsilon \:
\partial_\epsilon + d$ denote the differential on
$\overline{\Omega}^{\prime \prime, *}(\widetilde{B}, 
{\frak B})$. By the preceding
constructions, there is
an analytic torsion form on $\widetilde{B}$ which we can write as
$\widetilde{\cal T} = {\cal T}(\epsilon) +
d\epsilon \wedge {\cal T}^\prime(\epsilon)$,
satisfying
\begin{equation} \label{productT}
\widetilde{d} \widetilde{\cal T} = \int_{\widetilde{Z}} 
e \left(T\widetilde{Z}, \nabla^{T\widetilde{Z}} \right) \wedge
CS \left( \nabla^{\widetilde{\cal E}}, h^{\widetilde{\cal E}} \right) -
CS \left( \nabla^{H (\widetilde{Z}; \widetilde{\cal E} 
\big|_{\widetilde{Z}} )}, h^{H (\widetilde{Z}; \widetilde{\cal E} 
\big|_{\widetilde{Z}} )} \right).
\end{equation}
By our assumptions, the right-hand-side of (\ref{productT}) vanishes.
Thus $\partial_{\epsilon} {\cal T}(\epsilon) = d {\cal T}^\prime
(\epsilon)$, from which the proposition follows.
\end{pf} 
\begin{proposition} \label{duality}
Suppose that\\
1. $Z$ is even-dimensional\\
2. $TZ$ is oriented\\
3. $\nabla^{\cal E}$ is self-adjoint with respect to $h^{\cal E}$.\\  
Then ${\cal T} = 0$. 
\end{proposition}
\begin{pf}
This follows from an argument using Hodge duality, as in
\cite[Theorem 3.26]{Bismut-Lott (1995)}. We omit the details.
\end{pf}

Let us look more closedly at ${\cal T}_{[0]}$, the component of ${\cal T}$ in 
$\overline{\Omega}^{\prime \prime, 0}(B, {\frak B})$.
Assume for simplicity that $B$ is connected. Then
\begin{equation}
\overline{\Omega}^{\prime \prime, 0}(B, {\frak B}) =
\left( C^\infty(B)/\C\right) \otimes 
\left({\frak B}/\overline{[{\frak B},{\frak B}]}\right).
\end{equation}
As in (\ref{torequiv}),
\begin{equation}
{\cal T}_{[0]}  \equiv  
 - \int_0^\infty  \int_0^1 \Tr_s \left( N
\left( 1 - 2 t u(1-u) \triangle \right)
e^{- tu(1-u)\triangle}  \right) du 
\frac{dt}{t}.
\end{equation}
Define $g$ as in (\ref{g(t)}). Then
a specific lifting of  ${\cal T}_{[0]}$ to 
$C^\infty(B) \otimes \left({\frak B}/\overline{[{\frak B},{\frak B}]}\right)$
is given by
\begin{align}
{\cal T}_{[0]} & =  \int_0^\infty \left[ \Tr_s \left(
N g \left(t \triangle\right) \right) - 
\left( \frac{n}{2} \: \chi(Z) \: \Tr_s \left( I \big|_{\cal E} \right) -
\Tr_s \left( N \big|_{H(Z; {\cal E} \big|_Z)} \right) \right) g(t) \right. \\
& \left. \hspace{1in} +
\Tr_s \left( N \big|_{H(Z; {\cal E} \big|_Z)} \right) 
\right] \frac{dt}{t}. \notag
\end{align}
{\bf Example 10 : } If ${\frak B} = \C$ then as in
\cite[Theorem 3.29]{Bismut-Lott (1995)},
${\cal T}_{[0]}$ is 
the usual Ray-Singer analytic torsion
\cite{Ray-Singer (1971)}, considered to be a function on $B$.\\ \\
{\bf Example 11 : } Suppose that $\Gamma$ is a finite group and ${\frak B} =
\C \Gamma$.
Then ${\cal T}_{[0]}$ is equivalent to the equivariant analytic torsion
of \cite{Lott-Rothenberg (1991)}. \\ \\
{\bf Example 12 : } Suppose that $\Gamma$ is a discrete group and
${\frak B} = C^*_r \Gamma$.
Let $\tau$ be the trace on ${\frak B}$ given
by $\tau(\sum_{\gamma \in \Gamma} c_\gamma \gamma) = c_e$. 
Then $\tau \left( {\cal T}_{[0]} \right)$
is the $L^2$-analytic torsion of 
\cite{Lott (1992b),Mathai (1992)}.  
(In the cited papers, the $L^2$-torsion is defined using the group von
Neumann algebra and without the assumption of a gap in the spectrum of
$\triangle$, but with the assumption of positive Novikov-Shubin invariants.)

\section{Diffeomorphism Groups} \label{Diffeomorphism Groups}

Let $Z$ be a connected closed manifold and let $\Diff(Z)$ be its diffeomorphism
group, endowed with the natural smooth topology \cite{Milnor (1984)}.
For $i > 1$, let $\alpha : (S^i, *) \rightarrow (\Diff(Z), \Id)$ be a 
smooth map. We want to find invariants of $[\alpha] \in \pi_i(\Diff(Z))$.
Put $M_1 = M_2 = D^{i+1} \times Z$. Glue $M_1$ and $M_2$ along
their common boundary $S^i \times Z$ by identifying $(\theta, z) \in
\partial M_1$ with $(\theta, (\alpha(\theta))(z)) \in \partial M_2$. Let
$M = M_1 \cup_{S^i \times Z} M_2$ be the resulting manifold. Then
$M$ is the total space of a fiber bundle with base $B = S^{i+1}$ and
fiber $Z$. Any (smooth) topological invariant of such fiber bundles 
gives an invariant of $[\alpha]$.

As mentioned in the introduction, 
we are interested in the case when $Z$ is a $K(\Gamma, 1)$-manifold.
Then $\pi_1(M) = \Gamma$. Suppose that $\Gamma$ satisfies Hypothesis
\ref{introhypo} of the introduction.  There is a ${\frak B}$-vector
bundle ${\cal E} = {\frak B} \times_\Gamma \widetilde{M}$ on $M$.
Choosing $h \in C^\infty_0(\widetilde{M})$ satisfying (\ref{heqn}),
Proposition \ref{hprop} gives 
a partially flat connection $\nabla^{\cal E}$ on ${\cal E}$.
Let us add vertical Riemannian metrics $g^{TZ}$ and a horizontal distribution
$T^HM$ to the fiber bundle.

We would like to use the formalism of Section \ref{Fiber Bundles} to define
the analytic torsion form. By Proposition \ref{sa}, 
$\nabla^{\cal E}$ is self-adjoint and so Proposition \ref{duality} implies that
the torsion form vanishes if $\dim(Z)$ is even. (As the analysis is
effectively done on the universal cover $\widetilde{Z}$, the orientation
assumption on $TZ$ is irrelevant.) So assume that
$\dim(Z)$ is odd. 
Let ${\tau} \in Z^q(\Gamma; \C)$ be
a group cocycle and let $Z_\tau \in ZC^q({\frak B})$ be the 
cyclic cocycle coming from (\ref{cycliccocycle}) (with $x=e$).
We want to use
Proposition \ref{cclass} to
define the analytic torsion class $$[{\cal T}] \in 
H_{\frak B}^{\prime \prime,even}(B) = 
\bigoplus
\begin{Sb}
p + q \: even\\
p > q
\end{Sb}
\HH^p(B; \overline{H}_q({\frak B})),$$ take its integral over $B$ to get
$$\int_B [{\cal T}] \in \bigoplus
\begin{Sb}
q \equiv i+1 \mod  2\\
q < i + 1
\end{Sb}\overline{H}_q({\frak B})$$ and pair the result with  
$Z_\tau$.

In order to satisfy the hypotheses of Proposition \ref{cclass}, we need to
know that $H(Z; {\cal E}\big|_Z) = 0$ and that Hypothesis \ref{gaphypothesis}
is satisfied. Equivalently, we need to know that 
the $p$-form Laplacian on $\widetilde{Z}_b$ is invertible for all
$0 \le p \le \dim(Z)$. This is a topological condition on $Z$, but it
seems likely that it
is never satisfied
\cite{Lott (1997)}. To 
understand the nature of the problem, let us look in detail at the case when
$\Gamma$ is a free abelian group.
\subsection{Free Abelian Fundamental Groups}
Suppose that $\Gamma = \Z^n$. Then $Z = T^n$. Let $B\Gamma$ and
$\widehat{\Gamma}$ denote the classifying space of $\Gamma$
and  the Pontryagin dual of
$\Gamma$, respectively.  They are again $n$-tori, but it will be convenient
to distinguish them from $Z$. 

Under Fourier transform, $C^*_r \Gamma \cong C(\widehat{\Gamma})$. Take
${\frak B} = C^\infty(\widehat{\Gamma})$. Instead of using the universal
GDA of ${\frak B}$, we will simplify things and use the GDA of
smooth differential forms on $\widehat{\Gamma}$. This allows us to use
ordinary ``commutative'' analysis. All of the relevant steps
of the paper go through with this replacement. We now summarize the
statements.

First, there is a natural Hermitian line bundle $H$ with Hermitian
connection $\nabla^H$ on 
$B\Gamma \times \widehat{\Gamma}$ \cite[Section 3.1.1]{Lott (1992a)}.
(The third line of that section should read $\HH_1(M; \Z)_{mod \: Tor} \subset
\HH_1(M; \R)$.)  For all $\theta \in \widehat{\Gamma}$, 
the restriction of $\nabla^H$ to $B\Gamma \times \{ \theta \}$ is the flat 
connection on $B\Gamma$ with holonomy specified by $\theta$. 

We assume that we are given a fiber bundle $Z \rightarrow M \rightarrow B$ as
above, endowed with a vertical Riemannian metric and a horizontal
distribution.
Consider the fiber bundle $Z \rightarrow M \times \widehat{\Gamma} \rightarrow
B \times \widehat{\Gamma}$. It inherits a vertical Riemannian metric and
a horizontal distribution.  Let $f : M \rightarrow B\Gamma$ be a 
classifying map for the universal cover
$\widetilde{M} \rightarrow M$. Put $E_0 = 
(f \times \Id)^* H$, a Hermitian line bundle on  
$M \times \widehat{\Gamma}$. The pulled-back connection $\nabla^{E_0}$ is
partially flat. 

Let $\triangle$ be the vertical Laplacian of the fiber bundle 
$Z \rightarrow M \times \widehat{\Gamma} \rightarrow
B \times \widehat{\Gamma}$, acting on 
$\Omega \left(Z; E_0 \big|_Z \right)$. Then $\triangle$ is invertible 
except on the fibers over $B \times \{1\} \subset B \times \widehat{\Gamma}$.
This lack of invertibility on $B \times \{1\}$ is responsible for the fact
that Hypothesis \ref{gaphypothesis} is not satisfied.  The effect is that 
the analytic torsion form may be singular on $B \times \widehat{\Gamma}$,
with singularity along $B \times \{1\}$.

In order to get around this problem, one approach is to just remove the
singular subspace from consideration. Let $U \subset \widehat{\Gamma}$ be a
small neighborhood of $1 \in \widehat{\Gamma}$. Consider the restriction of
the fiber bundle to $B \times (\widehat{\Gamma} - U)$. Then the vertical
Laplacian is invertible and we can define the analytic torsion class
$$[{\cal T}] \in 
\bigoplus
\begin{Sb}
p + q \: even\\
p > q
\end{Sb}
\HH^p(B; \C) \otimes \HH^q(\widehat{\Gamma} - U; \C).$$
Now $\HH^q(\widehat{\Gamma} - U ; \C) \cong \HH_q(\Gamma; \C)$ if
$0 \le q < n$, and $\HH^n(\widehat{\Gamma} - U ; \C) = 0$. Thus
there is a (smooth) topological invariant of the fiber bundle given by
\begin{equation} \label{toruscase}
\int_B [{\cal T}] \: \in \bigoplus
\begin{Sb}
q \equiv i+1 \mod  2\\
q < \min(i + 1, n)
\end{Sb} \HH_q(\Gamma; \C).
\end{equation}
In fact, an argument involving complex conjugation shows that
the component of $\int_B [{\cal T}]$ in $\HH_q(\Gamma; \C)$ vanishes
unless $q \equiv i+1 \mod  4$. 

Comparing (\ref{toruscase}) with (\ref{pii})
(in the case $\Gamma = \Z^n$)
shows that $\int_B [{\cal T}]$ potentially detects all of
$\pi_*(\Diff(T^n)) \otimes_{\Z} \C$ in the stable range.  By removing $U$ from
$\widehat{\Gamma}$, we have lost the component of $\int_B [{\cal T}]$
in $\HH_n(\Gamma; \C)$ if $i+1 > n$, but this lies outside of the 
stable range, anyway. \\ \\
{\bf Remark 9 : } 
Although Hypothesis \ref{gaphypothesis} is not satisfied for the bundle
$Z \rightarrow M \rightarrow B$, we have seen that it is nevertheless
possible to extract most of the information in $[{\cal T}]$, due to the
fact that $\triangle$ is noninvertible only on a high-codimension subset
of $B \times
\widehat{\Gamma}$. Although $[{\cal T}]$ is possibly singular
on $B \times \{1\}$, it may be that its singularity is
sufficiently mild to still define the component of $\int_B [{\cal T}]$
in $\HH_n(\Gamma; \C)$. We have not looked at this point in detail.\\

In summary, when $\Gamma = \Z^n$ then a certain part of the  analytic 
torsion form is well-defined directly.  
We do not know what the situation is for the analytic torsion form in the case
of general $\Gamma$. In the next subsection we will make use of the
homotopic triviality of the fiber bundle in order to define a relative 
analytic torsion class for general $\Gamma$. 

\subsection{General Fundamental Groups}
Let $Z \rightarrow M \stackrel{\pi}{\rightarrow} B$ and ${\cal E}$ be as
described at the beginning of
Section \ref{Diffeomorphism Groups}. We again assume that $\Gamma$ satisfies
Hypothesis \ref{introhypo}, that $Z$ is a $K(\Gamma, 1)$-manifold 
and that $\dim(Z)$ is odd. 
Let $M^\prime = Z \times B$ be the
product bundle over $B$. Let ${\cal E}^\prime$ be the corresponding
${\frak B}$-vector bundle on $M^\prime$. From homotopy theory, we know that 
$M$ and $M^\prime$ are fiber-homotopy equivalent by some smooth map
$f : M \rightarrow M^\prime$. Furthermore, $f$ is unique
up to homotopy. It induces an isomorphism between the local systems
${\cal E}$ and ${\cal E}^\prime$.
We will show that the problems with invertibility of
the Laplacian cancel out when we consider $M$ relative to $M^\prime$.

Consider the restriction of $f$ to a single fiber $Z_b$. 
It acts by pullback on differential forms. 
However, this need not be a bounded, or even closable, operator.  
To get around this problem, we use the trick of 
\cite{Hilsum-Skandalis (1992)}, which involves  modifying $f$ to make it a
submersion.

Let $i : Z \rightarrow \R^N$ be an embedding of $Z$ in Euclidean space.
For $\epsilon > 0$ sufficiently small, let $U$ be an $\epsilon$-tubular 
neighborhood of $i(Z)$, with projection
$P : U \rightarrow Z$. Let $p_1 : Z \times B \rightarrow Z$ be projection
on the first factor. Let $\B^N$ denote the unit ball in $\R^N$. Consider the
fiber bundle $\B^N \times Z \rightarrow \B^N \times M \rightarrow B$.
Define $F : \B^N \times M \rightarrow M^\prime$ by
\begin{equation}
F(\vec{x}, m) = \left(P\left(\frac{\epsilon}{2} \: \vec{x} + 
(i \circ p_1 \circ f)(m) \right), \pi(m)\right).
\end{equation}
Then $F$ is a fiber-homotopy equivalence which is a fiberwise submersion.
Choose $\nu \in \Omega^N(\B^N)$ with support near $0 \in \B^N$ and total
integral $1$. Define
${\cal W}$ as in Subsection \ref{Induced Superconnections} and let
$\widetilde{\cal W}$ be the analogous object for the fiber bundle $M^\prime$.
Define a cochain map $T : \widetilde{\cal W} \rightarrow {\cal W}$ by
\begin{equation}
T(s^\prime) = \int_{\B^N} \nu \wedge F^* s^\prime,
\end{equation}
where $F^*$ acts fiberwise. Then $T$ is bounded.

Put $\widehat{\cal W} = {\cal W} \oplus \widetilde{\cal W}$, with the
$\Z$-grading $\widehat{\cal W}^i = {\cal W}^i \oplus (
\widetilde{\cal W}^\prime)^{i+1}$.
Let $A^\prime$ be the superconnection on ${\cal W}$ defined in 
Subsection \ref{Induced Superconnections} and let $\widetilde{A}^\prime$ be
the analogous superconnection on $\widetilde{\cal W}$. For $r \in \R$, 
define a superconnection $\widehat{A}_r^\prime$ on $\widehat{\cal W}$ by
\begin{equation}
\widehat{A}_r^\prime = 
\begin{pmatrix}
A^\prime & r T (-1)^N \\
0 & \widetilde{A}^\prime
\end{pmatrix}.
\end{equation}
The cochain part $\widehat{A}_{r,0,0,1}^\prime$ of $\widehat{A}_r^\prime$ is
\begin{equation}
\widehat{d}_r^Z = 
\begin{pmatrix}
d^Z & r T (-1)^N \\
0 & \widetilde{d}^Z
\end{pmatrix}.
\end{equation}
\begin{proposition}
The superconnection $\widehat{A}_r^\prime$ is partially flat on
$\widehat{\cal W}$.
\end{proposition}
\begin{pf}
It is enough to show that $A^{\prime,flat} T (-1)^N + 
T (-1)^N \widetilde{A}^{\prime,flat} = 0$.
Now $A^{\prime,flat}$ acts on $\Omega(B; {\cal W}) = 
\Omega(M; {\cal E})$ by exterior
differentiation $d^M$, and similarly for $\widetilde{A}^{\prime,flat}$.
Taking into account that $T$ is an odd variable, in ungraded language we
must show that $d^M T = T d^{M^\prime}$. 
As $T$ acts on $\Omega(M^\prime; {\cal E}^\prime)$ by
\begin{equation}
T(\omega^\prime) = \int_{\B^N} \nu \wedge F^* \omega^\prime,
\end{equation}
the proposition follows. 
\end{pf}

As $T$ is a cochain homotopy equivalence, if $r \ne 0$ then 
$H( \widehat{\cal W}, \widehat{d}_r^Z ) = 0$, while if $r = 0$ then 
$H^*( \widehat{\cal W}, \widehat{d}_r^Z ) = H^*(Z; {\cal E}_Z) \oplus
H^{*+1}(Z; {\cal E}_Z)$.

We now want to define an analytic torsion form ${\cal T}$
using the superconnection
$\widehat{A}_r^\prime$. If $t$ is large, we want $r$ to be nonzero, in order
to get the gap in the spectrum necessary for large-$t$ convergence of the
integral for ${\cal T}$. If $t$ is small, we want $r$ to be zero, in order to
use the small-time estimates separately on $M$ and $M^\prime$.

Choose a $\phi \in C^\infty_0([0, \infty))$ which is identically one near
$t = 0$. Put $h(t) = \sqrt{t} \: (1 - \phi(t))$ and 
\begin{equation}
\widehat{B}^\prime_t = t^{N/2} \widehat{A}^\prime_{h(t)} t^{-N/2}.
\end{equation}
The cochain part $\widehat{B}^\prime_{t,0,0,1}$ of $\widehat{B}^\prime_t$ is
\begin{equation}
\widehat{B}^\prime_{t,0,0,1} = 
\sqrt{t} \begin{pmatrix}
d^Z & (1 - \phi(t)) T (-1)^N \\
0 &  \widetilde{d}^Z
\end{pmatrix}.
\end{equation}
Let $\widehat{B}^{\prime \prime}_t$ be the adjoint of $\widehat{B}^\prime_t$.
For $u \in (0,1)$, put $\widehat{B}_t(u) = u \widehat{B}^\prime_t + (1-u)
\widehat{B}^{\prime \prime}_t$.
\begin{definition}
For $t > 0$, define $CS(t) \in 
\overline{\Omega}^{\prime \prime, odd}(B, {\frak B})$ and 
${\cal T}(t) \in \overline{\Omega}^{\prime \prime, even}(B, {\frak B})$ by
\begin{equation}
CS (t) = - \int_0^1 \Tr_s \left( \left( \widehat{B}^\prime_t -
\widehat{B}^{\prime \prime}_t \right) e^{- \widehat{B}_t^2(u)} \right) du
\end{equation}
and
\begin{align}
{\cal T}(t) = 
\: &- \: \frac1t \: 
\left( \int_0^1 \Tr_s \left( N
e^{- \widehat{B}_t^2(u)}  \right) du
+ \int_0^1 \int_0^1  u(1-u)  \right. \\
& \left. \hspace{.5in} \Tr_s \left(
N \left[ \widehat{B}_t^\prime - \widehat{B}^{\prime \prime}_t ,
e^{-r \widehat{B}_t^2(u)} 
\left( \widehat{B}_t^\prime - \widehat{B}^{\prime \prime}_t \right) 
e^{-(1-r) \widehat{B}_t^2(u)} \right] \right) dr du
\right) \notag \\
& +  \: h^\prime(t) \int_0^1 \int_0^1 \Tr_s \left( 
\begin{pmatrix}
0 & u T (-1)^N \\
(1-u) (-1)^N T^* & 0
\end{pmatrix} \right. \notag \\
& \left. \hspace{.5in} e^{-r \widehat{B}_t^2(u)} 
\left( \widehat{B}_t^\prime - \widehat{B}^{\prime \prime}_t \right) 
e^{-(1-r) \widehat{B}_t^2(u)}
\right) dr du. \notag
\end{align}
\end{definition}

\begin{proposition}
We have $\partial_t CS(t) = - d {\cal T}(t)$.
\end{proposition}
\begin{pf}
There is a partially flat superconnection on $\R^+ \times B$ given by
$d\epsilon \: \partial_\epsilon + \widehat{A}^\prime_{h(\epsilon)}$.
One can then proceed as in the proofs of Propositions \ref{CSvarprop} and
\ref{tt}. We omit the details.
\end{pf}

\begin{definition}
Define ${\cal T} \in \overline{\Omega}^{\prime \prime, even}(B, {\frak B})$ by
\begin{equation} \label{torint}
{\cal T} = \int_0^\infty {\cal T}(t) dt.
\end{equation}
\end{definition}

The integrand in (\ref{torint}) is integrable. For small $t$, this follows
from the fact that $1 - \phi$ vanishes identically near $t = 0$, so one
effectively has the difference of the torsion integrands of $M$ and $M^\prime$.
The large-$t$ convergence comes from 
the fact that $T$ is a cochain homotopy equivalence, which implies
that the Laplacian
$$\begin{pmatrix}
d^Z & T (-1)^N\\
0 & \widetilde{d}^Z
\end{pmatrix}^* 
\begin{pmatrix}
d^Z &  T (-1)^N \\
0 & \widetilde{d}^Z
\end{pmatrix} + 
\begin{pmatrix}
d^Z &  T (-1)^N \\
0 & \widetilde{d}^Z
\end{pmatrix} 
\begin{pmatrix}
d^Z &  T (-1)^N \\
0 & \widetilde{d}^Z
\end{pmatrix}^*
$$
is invertible \cite[Lemma 2.5]{Lott-Lueck (1995)}.

\begin{proposition}
The form ${\cal T}$ is closed.  Its class $[{\cal T}] \in
\overline{H}^{\prime \prime, even}(B, {\frak B})$ only depends
on $[\alpha] \in \pi_i(\Diff(Z))$.
\end{proposition}
\begin{pf}
This follows as in Corollary \ref{indep} and Proposition \ref{cclass}.
The only point to check is that any two choices of the auxiliary data
are connected by a smooth $1$-parameter family.  
This is obvious except, perhaps, for the choice
of embedding $i : Z \rightarrow \R^N$. If $i^\prime : Z \rightarrow 
\R^{N^\prime}$ is another choice, we can find some $N^{\prime \prime}$ and 
isometric embeddings $I : \R^N \rightarrow \R^{N^{\prime \prime}}$,
$I^\prime : \R^{N^\prime} \rightarrow \R^{N^{\prime \prime}}$ such that
$I \circ i$ and $I^\prime \circ i^\prime$ are connected by a smooth
$1$-parameter family of embeddings.
\end{pf}
 
So we have an invariant
\begin{equation}
\int_B [{\cal T}] \: \in \bigoplus
\begin{Sb}
q \equiv i+1 \mod  2\\
q < i + 1
\end{Sb} \overline{H}_q({\frak B}). 
\end{equation}
Because of the underlying real structures of the vector bundles involved,
one can show that the $\overline{H}_q({\frak B})$-component of
$\int_B [{\cal T}]$ vanishes unless $q \equiv i+1 \mod 4$. By Hypothesis
\ref{introhypo}, for each $[\tau] \in \HH^q(\Gamma; \C)$, there is a
representative $\tau \in Z^q(\Gamma; \C)$ such that 
$Z_\tau \in HC^q(\C \Gamma)$ extends to
a continuous cyclic cocycle on ${\frak B}$. Then the pairing
$\langle Z_\tau, \int_B [{\cal T}] \rangle \in \C$ is a numerical
invariant of $[\alpha] \in \pi_i(\Diff(Z))$.

\end{document}